\documentclass[10pt]{amsart}

\usepackage[dvipsnames]{xcolor}
\usepackage{cancel}
\usepackage{leftidx}
\usepackage{amsmath}
\usepackage{amscd}
\usepackage{amssymb}
\usepackage{rotating}
\usepackage{amsfonts}
\usepackage{amsthm}
\usepackage{verbatim}
\usepackage{nicematrix}
\usepackage{pdflscape}
\usepackage{bbm}
\usepackage[overload]{empheq}
\usepackage{mathtools}
\usepackage{arydshln}
\usepackage{cases}

\usepackage{cases}
\usepackage[top=1in, bottom=1in, left=1.25in, right=1.25in]{geometry}
\usepackage[colorlinks,linkcolor=blue, anchorcolor=blue, urlcolor=blue, citecolor=blue]{hyperref}

 \usepackage[hyperpageref]{backref}

\input xy
\xyoption{all}

\allowdisplaybreaks[1]

\NiceMatrixOptions{cell-space-top-limit = 1pt}

\newcommand{\pqed}{\hfill\qedsymbol\\}
\newcommand{\nn}{\nonumber}
\newcommand{\sfh}{\mathsf{h}}
\newcommand{\sfp}{\mathsf{p}}
\newcommand{\tsfh}{\tilde{\mathsf{h}}}

\newcommand{\Spec}{\mathrm{Spec}}
\newcommand{\wcM}{\widetilde{\mathcal{M}}}

\newtheorem{theorem}{Theorem}[section]
\newtheorem{proposition}[theorem]{Proposition}
\newtheorem{lemma}[theorem]{Lemma}
\newtheorem{corollary}[theorem]{Corollary}

\newtheorem{theorem/definition}[theorem]{Theorem/Definition}

\theoremstyle{remark}
\newtheorem{remark}{Remark}[section]

\theoremstyle{definition}
 \newtheorem{example}[theorem]{Example}
\newtheorem{definition}[theorem]{Definition}
\newtheorem{construction}[theorem]{Construction}

\begin{document}

\title{MIRROR SYMMETRY FOR QUADRIC HYPERSURFACES}
\author{Xiaowen Hu}
\address{School of Sciences, Great Bay University, Dongguan, China}
\email{huxw06@gmail.com}


\begin{abstract}
We compute Przyjalkowski-Shramov's resolution of the Calabi-Yau compactification of Givental's mirror Landau-Ginzburg model of the quadric hypersurfaces. We deduce the Picard-Fuchs equation for the narrow periods, which mirror the ambient quantum cohomology of quadric hypersurfaces. Then by an indirect approach using the irreducibility of the narrow Picard-Fuchs operator we deduce the Picard-Fuchs equation of the broad period, which mirrors the quantum cohomology of quadric hypersurfaces involving primitive cohomology classes.
The result suggests a natural choice of the opposite space in Barannikov's construction of Frobenius manifolds.
Finally, we show an isomorphism between the Frobenius manifolds associated with the quantum cohomology of a quadric hypersurface and its mirror Landau-Ginzburg model.
\end{abstract}

\maketitle

\tableofcontents

\section{Introduction}
For a Fano manifold $X$, mirror symmetry predicts that there exists a Landau-Ginzburg model $(V,f)$, where $V$ is a smooth quasi-projective variety, $w$ a regular function on $V$ (called \emph{potential function}), such that the Frobenius manifold associated with the genus 0 Gromov-Witten invariants of $X$ (called \emph{A-model}) is isomorphic to the Frobenius manifold associated with  $(V,w)$ (called \emph{B-model}). The construction of the Frobenius manifold associated with a Landau-Ginzburg model dates back to K. Saito \cite{Sai83}, in the name of \emph{flat structures}. It is generalized to various types of Landau-Ginzburg models by Barannikov \cite{Bar00}, and Douai-Sabbah \cite{DS03}, just to mention a few. 

For a Fano complete intersection $X$ in projective spaces, there is a mirror LG model defined by Givental (\cite{Giv95}, see (\ref{eq-LGmodel-Givental}) in Section \ref{sec:compactification-desingularization}). This LG model does not give a full mirror symmetry in terms of Frobenius manifolds, unless $X$ is a projective space. For example, when $X$ is a quadric hypersurface, Givental's mirror LG model does not have the predicted number of critical points. It is expected that a fiberwise compactification (i.e. making the potential function $w$ proper)  will remedy this.

In \cite{PRW16}, Pech, Rietsch, and Williams constructed a partial fiberwise compactification of  Givental's LG model, from a Lie theoretic or representation theoretic view, and showed a version of mirror symmetry. More precisely, if $X$ is an odd dimensional (resp. even dimensional) quadric hypersurface, they showed an isomorphism (resp. injective homomorphism) from the $D$-module associated with the quantum cohomology of $X$ to the $D$-module associated with the Gauss-Manin connection in their LG model. Lam and Templier \cite{LT17} then showed  that the morphism of $D$-modules is an isomorphism, as a special case of  mirror symmetry for minuscule flag varieties, by using tools coming from geometric representation theory.

In this paper we adopt the approach \cite{PS15} of Przyjalkowski-Shramov   to compactify and desingularize Givental's Landau-Ginzburg model
\begin{equation}\label{eq-mirrorLGModel-quadric-Givental-intro}
      \frac{q(x_1+1)^2}{x_1  \prod_{i=1}^{n-1}y_i}+y_1+\dots+y_{n-1}:(\mathbb{C}^{\times})^n\rightarrow \mathbb{C}.
\end{equation}
 The resulted model $(\overline{Z}^n,W_0)$ has an affine open subset $Z^n\cong \Spec[z_1,\dots,z_n,\frac{1}{z_1\cdots z_n-1}]$, and the restriction of the potential function to $Z^n$ is
\begin{equation}
      W_0(z_1,\cdots,z_n,q)=\frac{q z_1^2 z_{n} }{\prod_{i=1}^{n}z_i-1}+z_2+\sum_{i=2}^{n-1}z_{i}z_{i+1},
\end{equation}
where $q\in \mathbb{C}$ is regarded  as parameterizing a family of such models. The volume form $\frac{\mathrm{d}x_1 \mathrm{d}y_1\cdots \mathrm{d}y_{n-1}}{x_1y_1\cdot y_{n-1}}$ extends to be
\begin{equation}\label{eq-volumeForm-compatifiedLGModel-mirrorQuadric-intro}
            \Omega=\frac{dz_1\cdots dz_n}{z_1\cdots z_n-1}.
\end{equation}
We take $(Z^n,W_0)$  as our partially compactified LG model.  It is birationally equivalent to the model of Pech-Rietsch-Williams, and in dimension 3 birationally equivalent to the model of Gorbounov-Smirnov \cite{GS15}. It has the expected property that all critical points of $W_0$ lie in $Z^n$.
In our computation of the ingredients for the construction of B-model Frobenius manifolds, we do not explicitly use the proper model $(\overline{Z}^n,W_0)$. Nevertheless, the existence of such a model with this property plays a crucial role 
(see Lemma \ref{lem-pairing-oscilatoryIntegral}, \ref{lem-algebraicity-restrictionOfcalE}, and Corollary \ref{cor-APrioriPicardFuchs-broad}). The perspective that Landau-Ginzburg models need to be appropriately compactified first appeared in \cite{GKR12}, and our main references are \cite[\S 3.5]{Gro11} and \cite{KKP17}.

We adopt Barannikov's construction \cite{Bar00} of Frobenius manifolds associated with Landau-Ginzburg models, as  detailed and clarified in \cite{Gro11}. 
We denote by $\mathrm{FM}_B(Z^n,W_0)$ the resulted formal Frobenius manifold, and by $\mathrm{FM}_A(Q^n)$ the formal Frobenius manifold associated with the quantum cohomology of $Q^n$.
The following is our main theorem.

\begin{theorem}[= Theorem \ref{thm-mirrorSym-FrobManifolds-quadric}]\label{thm-mirrorSym-FrobManifolds-quadric-intro}
There exists an isomorphism of formal Frobenius manifolds from $\mathrm{FM}_A(Q^n)$ to $\mathrm{FM}_B(Z^n,W_0)$.
\end{theorem}
As $\mathrm{FM}_A(Q^n)$ is the germ of an analytic Frobenius manifold (\cite[Theorem 5.1]{HK21}), so is $\mathrm{FM}_B(Z^n,W_0)$.

For the construction of $\mathrm{FM}_B(Z^n,W_0)$ we find first a universal unfolding of $W_0$:
\begin{eqnarray}\label{eq-unfolding-LGMirrorOfQuadric-intro}
W=t_0+W_0+\sum_{i=2}^{n}t_i W_0^i-2e^{t_1}t_n+t_{n+1} e^{t_1}z_1.
\end{eqnarray}
In Section \ref{sec:universalUnfolding} we explain the requirement for $W$ and why we choose this one. This immediately yields a \emph{semi-infinite variation of Hodge structures}. To obtain a Frobenius manifold, we need to compute the Gauss-Manin connection, and  choose an \emph{opposite subspace} of the space of flat sections. We compute the Gauss-Manin connection on differential forms, and simultaneously the periods of the forms over the Lefschetz thimbles. For the  Lefschetz thimble $\Xi$ at an arbitrary critical point (an explicit Lefschetz thimble can be given by choosing a riemannian metric on $\overline{Z}^n$) we call 
\begin{equation*}
      \psi=\int_{\Xi}e^{\frac{W_0}{\hbar}}\Omega
\end{equation*}
a \emph{narrow} period, as it corresponds to quantum cohomology of $Q^n$ with only ambient cohomology classes as insertions (i.e. those pulled back from the ambient projective space). When the dimension $n$ is odd, the quadric hypersurface $Q^n$ has only ambient cohomology, and the Picard-Fuchs equation for narrow periods suffices to give the B-model Frobenius manifold. Suppose now that $n$ is even and $\geq 4$. We call 
\begin{equation*}
      \theta=\int_{\Xi}z_1e^{\frac{W_0}{\hbar}}\Omega.
\end{equation*}
a \emph{broad} period, as it corresponds to quantum cohomology of $Q^n$ with both primitive and ambient cohomology classes as insertions.
The computation of the Picard-Fuchs equations of the narrow periods is parallel to that in the case of projective spaces, while that of the broad periods is the main technical point of this paper. We do this by an indirect method. First we show, in terms of twisted de Rham cohomology classes and Gauss-Manin connections,
\begin{eqnarray}\label{eq-GaussManinConnection-broadToNarrow-intro}
(\nabla_{\hbar\partial_{\hbar}}^{\mathrm{GM}}-\frac{n}{2})^{\frac{n}{2}}[z_1\Omega]
=\frac{\hbar^{\frac{n}{2}}}{2(-n)^{\frac{n}{2}}q}(\nabla_{\hbar\partial_{\hbar}}^{\mathrm{GM}})^{n}[\Omega].
\end{eqnarray}
Secondly, by considering monodromy of homology basis $\Xi$, we show that a Picard-Fuchs operator for $[z_1\Omega]$ should satisfy certain constraints.
Then by an irreducibility result \cite{BBH88} of the Picard-Fuchs operator for $[\Omega]$, we obtain: let $D=\hbar\frac{\partial}{\partial \hbar}$, then for any Lefschetz thimble $\Xi$,
\begin{eqnarray}\label{eq-PicardFuchsOperator-broad-intro}
\big((D^2-\frac{n^2}{4})^{\frac{n}{2}+1}-4n^n q\hbar^{-n}(D-\frac{n}{2})(D-n)\big)\int_{\Xi}z_1e^{\frac{W_0}{\hbar}}\Omega=0.
\end{eqnarray}
The space of periods  is given in Corollary \ref{cor-periods-final}, from which there follows a seemingly natural choice of opposite space $\mathcal{H}_{-}$ of the space of flat sections (see Section \ref{sec:opposite space}). The final results of Section \ref{sec:Picard-Fuchs} might also be obtained from results in terms of $D$-modules in \cite{PRW16} and \cite{LT17}, but I did not study this in detail.

Finally we show the isomorphism of Frobenius manifolds by a reconstruction result similar to \cite[Section 5]{HK21}. A difference to loc. cit. is that in the present paper we show the vanishing of correlators with an odd number of broad insertions by the WDVV equation, and not using the deformation invariance of Gromov-Witten invariants.

There are difficulties in showing  mirror symmetry, in terms of Frobenius manifolds, for other Fano complete intersections in projective spaces. A crucial difficulty lies in the current constructions of Frobenius manifolds for Landau-Ginzburg models: they work only for those potential functions with only isolated critical points. For example, as we showed in \cite{Hu21}, the Frobenius manifold associated with the quantum cohomology of an even dimensional intersection of multidegree $\mathbf{d}=(2,2)$ is generically semisimple\footnote{This is predicted by a conjecture of Dubrovin. Indeed, the only smooth complete intersections in projective spaces having full exceptional collections are cubic surfaces, quadric hypersurfaces, and even dimensional intersections of multidegree $\mathbf{d}=(2,2)$.}. Establishing mirror symmetry for this class of complete intersections is desirable from the view of quantum cohomology (see the introduction of \cite{Hu21}). However, computations as in Section \ref{sec:compactification-desingularization} show that the corresponding compactified and desingularized mirror LG model has a non-isolated critical locus, which is similar to 3-dim intersections of multidegree $\mathbf{d}=(2,2)$ in \cite[Section 5]{KKOY09}.  We refer the reader to \cite{LW22} for an attempt towards Frobenius manifolds associated with LG models with non-isolated critical locus,.

The paper is organized as follows. In Section \ref{sec:compactification-desingularization} we recall Przyjalkowski-Shramov's compactification of desingularization of Givental's LG model. Then we compute the case of the mirror of $Q^n$. 
In Section \ref{sec:criticalPoints} we compute the critical points of our LG model and show the  isomorphism from the Milnor ring to the small quantum cohomology of $Q^n$. Then we explain our choice of universal unfolding $W$. 
 In Section \ref{sec:preparation-BModel-FrobeniusManifold} we recall Barannikov's construction of Frobenius manifolds from \cite[Chap. 2]{Gro11}. In Section \ref{sec:Picard-Fuchs} we compute the Picard-Fuchs equations in the way we outlined above. Then in Section \ref{sec:mirrorSymmetryFrobeniusManifolds} we collect the ingredients of Barannikov's construction and show the mirror isomorphism of Frobenius manifolds.
In Section \ref{sec:preparation-BModel-FrobeniusManifold} and Section \ref{sec:mirrorSymmetryFrobeniusManifolds}, the proofs that are verbatim as the case of the mirror LG model for $\mathbb{P}^n$ in \cite[Chap. 2]{Gro11}  are only briefly sketched or omitted.
 \\

\emph{Acknowledgement}:
I am grateful to Weiqiang He, Hua-Zhong Ke, Changzheng Li, Toshio Oshima, Helge Ruddat, and especially Hao Wen for very helpful discussions. 
This work is supported by NSFC 11701579 and NSFC 11831017.

\section{Przyjalkowski-Shramov desingularization}\label{sec:compactification-desingularization}

\begin{definition}
A \emph{quasiprojective (resp. projective) Landau-Ginzburg model} is a triple $(V,f,\mathrm{vol}_V)$ where $V$ is a smooth quasiprojective variety over $\mathbb{C}$, $f:V\rightarrow \mathbb{A}^1$ a dominant (resp. projective) morphism, and $\mathrm{vol}_V$ a nowhere vanishing section (resp. a meromorphic section) of the canonical sheaf $K_V$.
\end{definition}

Let $V_{\mathbf{d}}\subset \mathbb{P}^{n+r}$ be a smooth complete intersection of multidegree $\mathbf{d}=(d_1,\dots,d_r)$ in $\mathbb{P}^{n+r}$. Let
\begin{equation*}
      l=n+r-\sum_{i=1}^r d_i=\mbox{Fano index}-1.
\end{equation*}
The family of Givental's mirror Landau-Ginzburg model is $({Y}_{\mathbf{d}},f_{\mathbf{d}})$
\begin{gather}\label{eq-LGmodel-Givental}
   f_{\mathbf{d}}: {Y}_{\mathbf{d}}=\mathbb{C}^{n}=\prod_{i=1}^r(\mathbb{C}^{\times})^{ d_i-1}\times \mathbb{C}^{l}\rightarrow \mathbb{C},\nn\\ 
   f_{\mathbf{d}}(x_{1,1},\dots,x_{1,d_1-1},\dots,x_{r,1},\dots,x_{r,d_r-1},y_1,\dots,y_l)\\
    =  \frac{q\prod_{i=1}^r (x_{i,1}+\dots+x_{i,d_i-1}+1)^{d_i}}{\prod_{i=1}^r\prod_{j=1}^{d_i-1}x_{i,j} \prod_{j=1}^{l}y_j}+y_1+\dots+y_{l}\nn
\end{gather}
parametrized by $q\in \mathbb{C}^{\times}$.

The graph of $f_{\mathbf{d}}$ has a singular fiberwise Calabi-Yau compactification, which is the family of hypersurfaces $\overline{Y}_{\mathbf{d}}$ in $\prod_{i=1}^r\mathbb{P}^{d_i-1}\times \mathbb{P}^{l}\times \mathbb{A}^{1}$ defined by
\begin{equation}\label{eq-sing-GiventalMirror-CYCompactifiation}
    qy_0^{l+1}\prod_{i=1}^r(x_{i,1}+\dots+x_{i,d_i})^{d_i}=(\lambda y_0-\sum_{i=1}^{l}y_i)\prod_{i=1}^r\prod_{j=1}^{d_i}x_{i,j}\prod_{i=1}^{l}y_i
\end{equation}
 parametrized by $q$, where $\lambda=f_{\mathbf{d}}$ on $Y_{\mathbf{d}}$.

In \cite[Resolution Process 4.4]{PS15} Przyjalkowski and Shramov constructed crepant resolutions of  (\ref{eq-sing-GiventalMirror-CYCompactifiation}), and regard $\lambda$ as the Landau-Ginzburg potential function on the resolved space. We recall their construction. Locally this hypersurface is defined by an equation of the form
\begin{equation}\label{eq-PS-hypersurface}
      a_1^{d_1}\cdots a_k^{d_k}= \sigma x_1\cdots x_s
\end{equation}
in the affine space $\mathbb{A}(a_1,\dots,a_k,x_1,\dots,x_s, \mu)=\Spec\ \mathbb{C}[a_1,\dots,a_k,x_1,\dots,x_s, \sigma]$, where $d_i>0$, and $\sigma\in \mathbb{C}$ is regarded a parameter. We say that this hypersurface has weight $(s,\sum_{i=1}^k d_i)$. 

When $d_1\geq s$, we blowup $\{a_1=x_1=\dots=x_s=0\}$. Then in the  chart $a_1\neq 0$, we get a hypersurface of weight $(s,\sum d_i-s)$, and in the chart $x_j\neq 0$, we get a hypersurface of weight $(s-1,\sum d_i+d_1-s)$. Both weights are lexicographically smaller than $(s,\sum d_i)$.

When $d_1< s$, we blowup $\{a_1=x_1=\dots=x_{d_1}=0\}$. Then in the  chart $a_1\neq 0$, we get a hypersurface of weight $(s,\sum d_i-d_1)$, and in the chart $x_j\neq 0$, we get a hypersurface of weight $(s-1,\sum d_i)$. Both weights are again lexicographically smaller than $(s,\sum d_i)$.

The process terminates at a hypersurface of the form (\ref{eq-PS-hypersurface}) with weight $(0,\sum d_i)$ or $(s,0)$, which is smooth.

There are choices in expressing the local structure of the singular hypersurface into the form (\ref{eq-PS-hypersurface}). So the resolution is not unique. It is shown in \cite{PS15} that one can make a consistent choice of centers in the process of resolution to obtain a global crepant resolution of (\ref{eq-sing-GiventalMirror-CYCompactifiation}).

\subsection{Mirror LG model of quadric hypersurfaces}
Let $Q^n\subset \mathbb{P}^{n+1}$ be a smooth quadric hypersurface of dimension $n$. 
The potential function of the Givental mirror is
\begin{equation}\label{eq-mirrorLGModel-quadric-Givental}
      \frac{q(x_1+1)^2}{x_1  \prod_{i=1}^{n-1}y_i}+y_1+\dots+y_{n-1}
\end{equation}
and its singular Calabi-Yau compactification is the family of hypersurfaces $\overline{Y}$ in $\mathbb{P}^{1}\times \mathbb{P}^{n-1}\times \mathbb{A}^{1}$ defined by
\begin{equation}
    q y_0^{n}(x_1+x_2)^2=(\lambda y_0- \sum_{i=1}^{n-1}y_i)x_1 x_2 \prod_{i=1}^{n-1}y_i
\end{equation}
 parametrized by $q\in \mathbb{C}^{\times}$. The singular locus  of $\overline{Y}$ is
 \begin{eqnarray*}
 && \bigcup_{1\leq i_1< i_2\leq n-1}\{ x_1+x_2=0, y_{i_1}=y_{i_2}=0 \}
  \cup \bigcup_{1\leq i\leq n-1} \{ x_1+x_2=0, y_i=0, \sum_{j=1}^{n-1}y_j=\lambda y_0
 \}\\
 && \cup \{y_0=0,\mbox{and at least two of } x_1,x_2,y_1,\dots,y_{n-1},\sum_{i=1}^{n-1}y_i\ \mbox{vanish}
 \}.
 \end{eqnarray*}
 We are going to compute the Przyjalkowski-Shramov resolution of $\overline{Y}$.

Along the inverse image of the  locus $\{y_0=0\}$, $\lambda$ will be a coordinate on the local charts after a Przyjalkowski-Shramov desingularization. To obtain the critical locus of the function $\lambda$ it suffices to consider the affine open subset
\begin{equation}
      Y_0=\{y_0\neq 0,x_2\neq 0\} \cap \overline{Y}.
\end{equation}
 In the process of Przyjalkowski-Shramov resolution recalled above, in every step up we blow up along $V\big((x_1+1)^{i},y_a^{(j)},y_b^{(k)}\big)$ in a suitable affine chart, where $1\leq a<b\leq n-1$, $(x_1+1)^{i}$ is the defininig equation of the proper transform of $V(x_1+1)$ in this chart, and $y_a^{(j)}$ is the defining equation of the proper transform of $V(y_a)$ in this chart, and similarly for $y_b^{(k)}$. In the affine chart $U(y_b^{(k)}\neq 0)$ of the blowing-up, the defining equation of the proper transform of $V(y_a)$ is denoted by $y_a^{(j+1)}$, while in the affine chart $U(y_a^{(j)}\neq 0)$ of the blowup, the defining equation of the proper transform of $V(y_b)$ is denoted by $y_b^{(k+1)}$; in both charts the defining equation of the proper transform of $V(x_1+1)$ is denoted by $(x_1+1)^{(i+1)}$. Note that the open subset $U\big((x_1+1)^{(i+1)}\neq 0\big)$ is contained in the domain of the Givental LG model.   So a chart of the final resolution that may contain new critical points corresponds to a \emph{decision sequence}, which chooses the affine chart after each step of blowing-up that one needs to further blow up. We illustrate the process of Przyjalkowski-Shramov resolution by a \emph{decision tree}, as in Figure \ref{fig:decisionTree}. 

\newgeometry{left=1in,right=1in,top=1in,bottom=1in,nohead}
\begin{landscape}
(\emph{bla} is short for \emph{blow up along}.)
\begin{figure}[htbp]
\centering
\caption{Decision Tree}\label{fig:decisionTree}
\[
\scalebox{0.8}{%
\xymatrix{
&&&& &&& Y_0 \ar[d]^{\mbox{blow up along } V(x_1+1,y_1,y_2)}  &&& &&&& \\
&&&& &&& \ar[dllll]_{y_1^{(1)}\neq 0}  \ar[drrrr]^{y_2^{(1)}\neq 0} &&& &&&& \\
&&& Y_1  \ar[d]^{\mbox{bla } V\big((x_1+1)^{(1)},y_2^{(1)},y_3\big)}  
&&&&&&&& Y_2  \ar[d]^{\mbox{bla } V\big((x_1+1)^{(1)},y_1^{(1)},y_3\big)} &&\\
&&&  \ar[dll]_{y_2^{(2)}\neq 0}  \ar[drr]^{y_3^{(1)}\neq 0} && &&&& &&
 \ar[dll]_{y_1^{(2)}\neq 0}  \ar[drr]^{y_3^{(1)}\neq 0} && \\
&Y_{12} \ar[d]_{\mbox{bla } V\big((x_1+1)^{(2)},y_3^{(1)},y_4\big)} &&&& 
Y_{13} \ar[d]_{\mbox{bla } V\big((x_1+1)^{(2)},y_2^{(2)},y_4\big)} &&&&
Y_{21} \ar[d]_{\mbox{bla }  V\big((x_1+1)^{(2)},y_3^{(1)},y_4\big)} &&&& 
Y_{23} \ar[d]_{\mbox{bla } V\big((x_1+1)^{(2)},y_1^{(2)},y_4\big)} & &\\
& \ar[dl]_{y_3^{(2)}\neq 0}  \ar[dr]^{y_4^{(1)}\neq 0} &&&&
 \ar[dl]_{y_2^{(3)}\neq 0}  \ar[dr]^{y_4^{(1)}\neq 0}&&&&
 \ar[dl]_{y_3^{(2)}\neq 0}  \ar[dr]^{y_4^{(1)}\neq 0} &&&&
  \ar[dl]_{y_1^{(3)}\neq 0}  \ar[dr]^{y_4^{(1)}\neq 0}&& \\
Y_{123} \ar[d] && Y_{124} \ar[d] && Y_{132}\ar[d] && Y_{134}\ar[d] &&  Y_{213} \ar[d]&& Y_{214}\ar[d] 
&& Y_{231}\ar[d] && Y_{234}\ar[d] \\
\dots && \dots && \dots && \dots && \dots && \dots && \dots && \dots
}
}
\]
\end{figure}
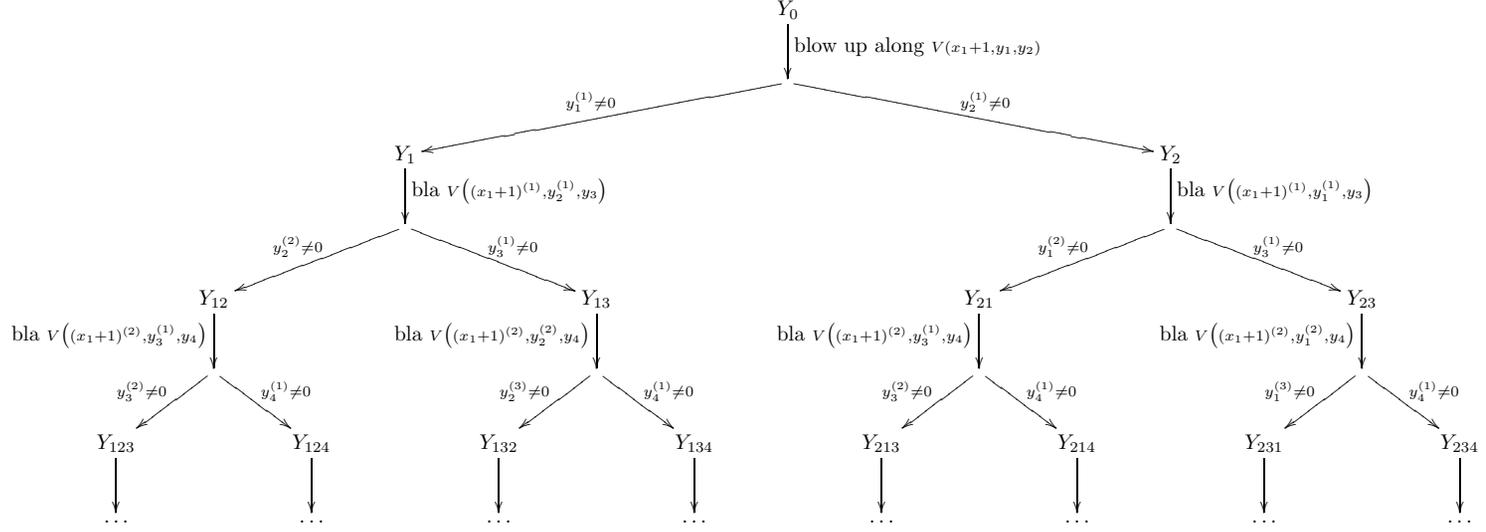

\normalsize{
where, for example,
\begin{eqnarray*}
&&Y_0=\{qy_0^{n}(x_1+1)^2=(\lambda y_0-y_1-\dots-y_{n-1})x_1 y_1 \dots y_{n-1}\},\\
&&Y_1=\{qy_0^n\big((x_1+1)^{(1)}\big)^2=(\lambda y_0-y_1-y_2^{(1)}y_1-\prod_{i=3}^{n-1}y_i)\big((x_1+1)^{(1)}y_1-1\big)y_2^{(1)}\prod_{i=3}^{n-1}y_i\},\\
&&Y_2=\{qy_0^n\big((x_1+1)^{(1)}\big)^2=(\lambda y_0-y_1^{(1)}y_2-\prod_{i=2}^{n-1}y_i)\big((x_1+1)^{(1)}y_2-1\big)y_1^{(1)}\prod_{i=3}^{n-1}y_i\},\\
&& Y_{23}=\{qy_0^n\big((x_1+1)^{(2)}\big)^2 =(\lambda y_0-y_1^{(2)}y_2y_3-\sum_{i=2}^{n-1}y_i)\big((x_1+1)^{(2)}y_2y_3-1\big)y_1^{(2)}\prod_{i=4}^{n-1}y_i\},\\
&& Y_{231}=\{qy_0^n\big((x_1+1)^{(3)}\big)^2 =(\lambda y_0-y_1^{(2)}y_2y_3-y_2-y_3-y_1^{(2)}y_4^{(1)}-\sum_{i=5}^{n-1}y_i)\big((x_1+1)^{(3)}y_1^{(2)}y_2y_3-1\big)y_4^{(1)}\prod_{i=5}^{n-1}y_i\}.
\end{eqnarray*}
}
\end{landscape}
\restoregeometry

\begin{definition}
Let $N$ be a natural number.
Suppose $1= b_1<b_2<\dots<b_r\leq N$. Let $l_i=b_{i+1}-b_i$ for $1\leq i\leq r-1$, and $l_r=N+1-b_r$. For $1\leq i\leq r$, let $c_i$ be the sequence $(b_i,b_i+1,\dots,b_i+l_i-1)$. Then the concatenation $c_1 \cdots c_r$ is equal to the sequence $S=(1,2,\dots,N)$ and we call the ordered collection of subsets $c_1,\dots,c_r$ a \emph{subdivision} of $S$.
\end{definition}

\begin{definition}

With a subdivision $c_1,\dots,c_r$ of $S=\{1,2,\dots,n-1\}$, we associate an affine chart in the final resolution as follows. 
\begin{enumerate}
       \item  First we blow up $U_0=Y_0$ along $V(x_1+1,y_1,y_2)$, and denote the proper transform in the chart $y_2^{(1)}\neq 0$ of the blowing-up by $U_1$.   For $k$ from 2 to $l_1-1$, we blowup $U_{k-1}$ along $V\big( (x_1+1)^{(k-1)},y_1^{(k-1)},y_{k+1})$, and denote the proper transform in  the chart $y_{k+1}^{(1)}\neq 0$ of the blowing-up by $U_k$. Then blowup $U_{l_1-1}$ along $V\big( (x_1+1)^{(l_1-1)},y_1^{(l_1-1)},y_{l_1+1})$, and denote the proper transform in  the chart $y_{1}^{(l_1)}\neq 0$ of the blowing-up by $U_{c_1}$. 
       \item We repeat the first step inductively. More precisely, suppose that we have obtained $U_{c_1\dots c_{e-1}}$. Let $K=\sum_{i=1}^{e-1} l_i$. We blow up $U_0=U_{c_1\dots c_{e-1}}$ along $V\big((x_1+1)^{K},y_{b_{e}}^{(1)},y_{b_{e}+1})$, and denote the proper transform in  the chart $y_{b_e+1}^{(1)}\neq 0$ of the blowing-up by $U_1$.   For $k$ from 2 to $l_e-1$, we blowup $U_{k-1}$ along $V\big( (x_1+1)^{(K+k-1)},y_{b_{e}}^{(k-1)},y_{b_{e}+k})$, and denote the proper transform in the chart 
       $y_{b_{e}+k}^{(1)}\neq 0$ of the blowing-up by $U_k$. Then blowup $U_{l_e-1}$ along $V\big( (x_1+1)^{(K+l_e-1)},y_{b_{e}}^{(l_e-1)},y_{b_{e}+l_e})$, and denote the proper transform in  the chart $y_{b_{e}}^{(l_e)}\neq 0$ of the blowing-up by $U_{c_1\dots c_e}$. 
 \end{enumerate} 
\end{definition}

\begin{example}
Suppose $n\geq 4$. Then $Y_{231}$ in Figure \ref{fig:decisionTree} is identified with $U_{(123)}$.
\pqed
\end{example}

\begin{example}
Suppose $n=7$. Then $U_{(123)(456)}$ is defined by
\begin{eqnarray*}
qy_0^n\big((x_1+1)^{(6)}\big)^2y_4^{(3)}& =&(\lambda y_0-y_1^{(2)}y_2y_3-y_2-y_3-y_1^{(2)}y_4^{(3)}y_5y_6-y_5-y_6)\\
&&\cdot\big((x_1+1)^{(6)}y_1^{(2)}y_2y_3y_4^{(3)}y_5y_6-1\big).
\end{eqnarray*}
So on the chart $y_0\neq 0$,
\begin{equation*}
      \lambda=\frac{q\big((x_1+1)^{(6)}\big)^2y_4^{(3)}}{(x_1+1)^{(6)}y_1^{(2)}y_2y_3y_4^{(3)}y_5y_6-1}
      +y_1^{(2)}y_2y_3+y_2+y_3+y_1^{(2)}y_4^{(3)}y_5y_6+y_5+y_6.
\end{equation*}
Tracing back the changes of coordinates  we have 
\[
x_1=(x_1+1)^{(6)}y_1^{(2)}y_2y_3y_4^{(3)}y_5y_6-1,\
y_4=y_1^{(2)}y_4^{(3)}y_5y_6,\
y_1=y_1^{(2)}y_2y_3.
\]
\pqed
\end{example}

If the decision sequence is $\mathbf{c}=(c_1,c_2,\dots,c_k)$, where $c_i=(b_{i},\dots, b_{i}+l_i-1)$, the final proper transform is defined by
\begin{eqnarray}\label{eq-definitionEquation-Uc}
     qy_0^n\big((x_1+1)^{(n-1)}\big)^2y_{b_{k}}^{(l_k-\delta_{k,1})}& =&(\lambda y_0-\sum_{i=1}^k(y_{b_{i-1}}^{(l_{i-1}-\delta_{i-1,1})}y_{b_{i}}^{(l_i-\delta_{i,1})}y_{b_{i}+1}\cdots y_{b_{i}+l_i-1}+y_{b_{i}+1}+\dots+y_{b_{i}+l_i-1})\nn\\
     &&\times\big((x_1+1)^{(n-1)} \prod_{i=1}^k(y_{b_{i}}^{(l_i-\delta_{i,1})}y_{b_{i}+1}\cdots y_{b_{i}+l_i-1})-1\big),
\end{eqnarray}
where $y_{b_{i-1}}^{(l_{i-1}-\delta_{i-1,1})}=1$ when $i=1$, as a convention. So on the chart $y_0\neq 0$,
\begin{eqnarray}\label{eq-potentialFunctionLambda-on-Uc}
\lambda&=&\frac{q\big((x_1+1)^{(n-1)}\big)^2y_{b_{k}}^{(l_k-\delta_{k,1})}}{(x_1+1)^{(n-1)} \prod_{i=1}^k(y_{b_{i}}^{(l_i-\delta_{i,1})}y_{b_{i}+1}
\cdots y_{b_{i}+l_i-1})-1}\nn\\
&&+\sum_{i=1}^k\big(y_{b_{i-1}}^{(l_{i-1}-\delta_{i-1,1})}y_{b_{i}}^{(l_i-\delta_{i,1})}y_{b_{i}+1}\cdots y_{b_{i}+l_i-1}+y_{b_{i}+1}+\dots+y_{b_{i}+l_i-1}\big).
\end{eqnarray}
The relation of the new coordinates to the original ones is
\begin{eqnarray}\label{eq-mirrorLG-quadric-UcChart-relationOfCoordinates}
\begin{cases}
      x_1=(x_1+1)^{(n-1)} \prod_{i=1}^k(y_{b_{i}}^{(l_i-\delta_{i,1})}y_{b_{i}+1}\cdots y_{b_{i}+l_i-1})-1,& \\
      y_{b_{i}}=y_{b_{i-1}}^{(l_{i-1}-\delta_{i-1,1})}y_{b_{i}}^{(l_i-\delta_{i,1})}y_{b_{i}+1}\cdots y_{b_{i}+l_i-1},& \mbox{for}\ 1\leq i\leq k.
\end{cases}      
\end{eqnarray}

\begin{proposition}\label{eq-criticalPoints-charts}
If $\mathbf{c}\neq (1)(2)\dots(n-1)$, then the critical points of $W^{\mathrm{PS}}$ in $U_{\mathbf{c}}$ are critical points of the Givental model.
\end{proposition}
\begin{proof}
We compute the partial derivatives as follows.
\begin{eqnarray}\label{eq-LGMirror-quadric-chartUc-gradient-1}
\frac{\partial\lambda}{\partial (x_1+1)^{(n-1)}}&=&
\frac{q(x_1+1)^{(n-1)} y_{b_{k}}^{(l_k-\delta_{k,1})} }{\big((x_1+1)^{(n-1)} \prod_{i=1}^k(y_{b_{i}}^{(l_i-\delta_{i,1})}y_{b_{i}+1}\cdots y_{b_{i}+l_i-1})-1\big)^2}\nn\\
&&\times\big((x_1+1)^{(n-1)} \prod_{i=1}^k(y_{b_{i}}^{(l_i-\delta_{i,1})}y_{b_{i}+1}\cdots y_{b_{i}+l_i-1})-2\big).
\end{eqnarray}
For $1\leq j<k$,
\begin{eqnarray}\label{eq-LGMirror-quadric-chartUc-gradient-bj1}
\frac{\partial\lambda}{\partial y_{b_{j}}^{(l_j-\delta_{j,1})} }&=&-\frac{q(x_1+1)^{(n-1)} y_{b_{k}}^{(l_k-\delta_{k,1})}  \prod_{1\leq i\leq k, i\neq j}(y_{b_{i}}^{(l_i-\delta_{i,1})}y_{b_{i}+1}\cdots y_{b_{i}+l_i-1}) \cdot (y_{b_{j}+1}\cdots y_{b_{j}+l_j-1}) }{\big((x_1+1)^{(n-1)} \prod_{i=1}^k(y_{b_{i}}^{(l_i-\delta_{i,1})}y_{b_{i}+1}\cdots y_{b_{i}+l_i-1})-1\big)^2}\nn\\
&&+y_{b_{j-1}}^{(l_{j-1}-\delta_{j-1,1})}y_{b_{j}+1}\cdots y_{b_{j}+l_j-1}+y_{b_{j+1}}^{(l_{j+1})}y_{b_{j+1}+1}\cdots y_{b_{j+1}+l_{j+1}-1},
\end{eqnarray}
and when $j=k$,
\begin{eqnarray}\label{eq-LGMirror-quadric-chartUc-gradient-bk1}
\frac{\partial\lambda}{\partial y_{b_{k}}^{(l_k-\delta_{k,1})} }&=&
\frac{q\big((x_1+1)^{(n-1)}\big)^2}{(x_1+1)^{(n-1)} \prod_{i=1}^k(y_{b_{i}}^{(l_i-\delta_{i,1})}y_{b_{i}+1}\cdots y_{b_{i}+l_i-1})-1}\nn\\
&&-\frac{q(x_1+1)^{(n-1)}  \prod_{1\leq i\leq k}(y_{b_{i}}^{(l_i-\delta_{i,1})}y_{b_{i}+1}\cdots y_{b_{i}+l_i-1})  }{\big((x_1+1)^{(n-1)} \prod_{i=1}^k(y_{b_{i}}^{(l_i-\delta_{i,1})}y_{b_{i}+1}\cdots y_{b_{i}+l_i-1})-1\big)^2}\nn\\
&&+y_{b_{k-1}}^{(l_{k-1}-\delta_{k-1,1})}y_{b_{k}+1}\cdots y_{b_{k}+l_k-1}.
\end{eqnarray}
For $1\leq j\leq k$ and  $1\leq e\leq l_j-1$,
\begin{align}\label{eq-LGMirror-quadric-chartUc-gradient-bje}
\frac{\partial\lambda}{\partial y_{b_{j}+e} }=&-\frac{q(x_1+1)^{(n-1)} y_{b_{k}}^{(l_k-\delta_{k,1})} }{\big((x_1+1)^{(n-1)} \prod_{i=1}^k(y_{b_{i}}^{(l_i-\delta_{i,1})}y_{b_{i}+1}\cdots y_{b_{i}+l_i-1})-1\big)^2}\\
&\cdot  \prod_{1\leq i\leq k, i\neq j}(y_{b_{i}}^{(l_i-\delta_{i,1})}y_{b_{i}+1}\cdots y_{b_{i}+l_i-1}) \cdot (y_{b_{j}}^{(l_j-\delta_{j,1})}y_{b_{j}+1}\cdots  \widehat{y_{b_{j}+e}}\cdots y_{b_{j}+l_j-1}) \nn\\
&+y_{b_{j-1}}^{(l_{j-1}-\delta_{j-1,1})}y_{b_{j}}^{(l_j-\delta_{j,1})}y_{b_{j}+1}\cdots  \widehat{y_{b_{j}+e}}\cdots y_{b_{j}+l_j-1}+1,\nn
\end{align}
where $\widehat{(\cdot)}$ means deleting this factor. Then we compute the critical points in $U_{\mathbf{c}}$.  By (\ref{eq-LGMirror-quadric-chartUc-gradient-1}),
\[
(x_1+1)^{(n-1)} y_{b_{k}}^{(l_k-\delta_{k,1})}
\big((x_1+1)^{(n-1)} \prod_{i=1}^k(y_{b_{i}}^{(l_i-\delta_{i,1})}y_{b_{i}+1}\cdots y_{b_{i}+l_i-1})-2\big)=0.
\]
If $(x_1+1)^{(n-1)} \prod_{i=1}^k(y_{b_{i}}^{(l_i-\delta_{i,1})}y_{b_{i}+1}\cdots y_{b_{i}+l_i-1})-2$, then by (\ref{eq-mirrorLG-quadric-UcChart-relationOfCoordinates}) the critical point lies in the Givental chart and is a critical point of Givental model.

Assume $(x_1+1)^{(n-1)}(p)=0$. Then by (\ref{eq-LGMirror-quadric-chartUc-gradient-bje}) we have
\begin{equation*}
      y_{b_{j}+e}=y_{b_{j}+1},\ \mbox{for}\ 1\leq j\leq k\ \mbox{and}\ 1\leq e\leq l_j-1.
\end{equation*}
By (\ref{eq-LGMirror-quadric-chartUc-gradient-bk1}),
\[
y_{b_{k-1}}^{(l_{k-1}-\delta_{k-1,1})}y_{b_{k}+1}\cdots y_{b_{k}+l_k-1}=0.
\]
In view of (\ref{eq-LGMirror-quadric-chartUc-gradient-bje}) with $j=k$, this is impossible  unless $l_k=1$. Hence $l_k=1$ and $y_{b_{k-1}}^{(l_{k-1}-\delta_{k-1,1})}=0$, and thus by (\ref{eq-LGMirror-quadric-chartUc-gradient-bje}) with $j=k-1$ we get $l_{k-1}=0$. Then using (\ref{eq-LGMirror-quadric-chartUc-gradient-bj1}) in place of (\ref{eq-LGMirror-quadric-chartUc-gradient-bk1}) and an induction on $r$, the same argument yields $l_{k-2r}=l_{k-2r-1}=1$ and $y_{b_{k-2r-1}}=0$  for $r\geq 0$. This is the chart $U_{(1)(2)\dots(n-1)}$ we dealt with above. 

Assume $y_{b_{k}}^{(l_k-\delta_{k,1})}=0$. Then (\ref{eq-LGMirror-quadric-chartUc-gradient-bje}) with $j=k$  is impossible unless $l_k=1$. The $j=k-1$ case of (\ref{eq-LGMirror-quadric-chartUc-gradient-bj1}) conflicts with those case of (\ref{eq-LGMirror-quadric-chartUc-gradient-bje}), unless $l_{k-1}=1$. Hence $l_{k-1}=l_k=1$ and $y_{b_{k-2}}^{(l_{k-2}-\delta_{k-2,1})}=0$. Then an induction on $r$ yields $l_{k-2r}=l_{k-2r-1}=1$ and $y_{b_{k-2r}}=0$  for $r\geq 0$. This again is the chart $U_{(1)(2)\dots(n-1)}$ we dealt with above.
\end{proof}

\begin{corollary}\label{cor-compatifiedLGModel-mirrorQuadric}
There is a projective Landau-Ginzburg model $(\overline{Z}^n,W_0)$, such that 
\begin{enumerate}
      \item[(i)] it extends the Givental model (\ref{eq-mirrorLGModel-quadric-Givental});
      \item[(ii)] it contains an open subset $Z^n$, which contains the critical points of $W_0$, and $\overline{Z}^n\backslash Z^n$ is a normal crossing divisor;
      \item[(iii)] there is an isomorphism $Z^n\cong\{z_1\cdots z_n\neq 1\}\subset \mathbb{C}^n$, and under this isomorphism we have
      \begin{equation*}
           W_0|_{Z^n}(z_1,\cdots,z_n,q)=\frac{q z_1^2 z_{n} }{\prod_{i=1}^{n}z_i-1}+z_2+\sum_{i=2}^{n-1}z_{i}z_{i+1}.
      \end{equation*}
      Moreover, the standard $n$-form
      \begin{equation*}
            \frac{dx_1dy_1\cdots dy_{n-1}}{x_1 y_1\cdots y_{n-1}}
      \end{equation*}
      on the torus $(\mathbb{C}^{*})^n$ extends to $Z^n$:
      \begin{equation}\label{eq-volumeForm-compatifiedLGModel-mirrorQuadric}
            \Omega=\frac{dz_1\cdots dz_n}{z_1\cdots z_n-1}.
      \end{equation}
\end{enumerate}
\end{corollary}
\begin{proof}
Let $\tilde{Z}^n$ be the Przyjalkowski-Shramov resolution of $\overline{Y}$, and $Z^n=U_{(1)(2)\dots(n-1)}$. Let $W_0=\lambda$. By Proposition \ref{eq-criticalPoints-charts}, $Z^n$ contains all critical points of $W_0$. By blowing up smooth centers contained in the reduced divisor $\tilde{Z}^n\backslash Z^n$ if necessary, one gets $\overline{Z}^n$ such that $\overline{Z}^n\backslash Z^n$ is a normal crossing divisor.
By (\ref{eq-definitionEquation-Uc}), (\ref{eq-potentialFunctionLambda-on-Uc}), and (\ref{eq-mirrorLG-quadric-UcChart-relationOfCoordinates}), 
on $U_{(1)(2)\dots(n-1)}$, 
$$(x_1+1)^{(n-1)},y_1,y_2^{(1)},\dots,y_{n-1}^{(1)}$$ form a system of coordinates, and 
\begin{equation}\label{eq-LGmodel-partialCompactification-desing-0}
      \lambda=\frac{q\big((x_1+1)^{(n-1)}\big)^2 y_{n-1}^{(1)} }{(x_1+1)^{(n-1)}y_1\prod_{i=2}^{n-1}y_i^{(1)}-1}+y_1+y_2^{(1)}y_1+\sum_{i=2}^{n-2}y_{i}^{(1)}y_{i+1}^{(1)}.
\end{equation}
We rewrite the potential (\ref{eq-LGmodel-partialCompactification-desing-0}) as
\begin{equation}\label{eq-LGmodel-partialCompactification-desing}
      \lambda(z_1,\cdots,z_n)=\frac{q z_1^2 z_{n} }{\prod_{i=1}^{n}z_i-1}+z_2+\sum_{i=2}^{n-1}z_{i}z_{i+1}
\end{equation}
defined on the open subset $\{\prod_{i=1}^n z_i\neq 1\}$ of $\mathbb{C}^n$, where
\begin{eqnarray}\label{eq-transform-xyToZ}
\begin{cases}
z_1=\begin{cases} 
\frac{x_1+1}{y_2 y_4 \cdots y_{n-1}},& \mbox{if}\ n\ \mbox{odd};\\
\frac{x_1+1}{y_1 y_3 \cdots y_{n-1}},& \mbox{if}\ n\ \mbox{even},
\end{cases}\\
z_{2i}= \frac{y_1 y_3\cdots y_{2i-1}}{y_2 y_4\cdots y_{2i-2}},& \mbox{for}\ 2\leq 2i\leq n,\\
z_{2i+1}=\frac{y_2 y_4\cdots y_{2i}}{y_1 y_3\cdots y_{2i-1}},&  \mbox{for}\ 3\leq 2i+1\leq n.
\end{cases}
\end{eqnarray}
Conversely,
\begin{equation}\label{eq-transform-zToXy}
\begin{cases}
x_1=\prod_{i=1}^n z_i-1,\\
y_1=z_2,\\
 y_i=z_{i}z_{i+1},& \mbox{for}\ 2\leq i\leq n-1.
\end{cases}
\end{equation}
So we get (\ref{eq-volumeForm-compatifiedLGModel-mirrorQuadric}).
\end{proof}

\begin{definition}
We define the quasiprojective Landau-Ginzburg model mirror to $Q^n$ to be $(Z^n,W_0,\Omega)$ in Corollary \ref{cor-compatifiedLGModel-mirrorQuadric}. Here we denote the restriction of $W_0$ on $Z^n$ still by $W_0$, as an abuse of notation.
\end{definition}

One can further blowup $\overline{Z}^n$ to extend $\Omega$ to be a meromorphic section, so that one obtains a projective Landau-Ginzburg model. For the purpose of this paper, it suffices to know that $Z^n$ has a fiberwise compactification with no critical points in the complement.

\section{The critical points and the universal unfolding}\label{sec:criticalPoints}
In this section we compute the critical points of $W_0$, and show that the Milnor ring of $W_0$ is isomorphic to the small quantum cohomology ring of $Q^n$. The structure of the Milnor ring suggests our choice of the universal unfolding $W$, as will be explained in Remark \ref{rem:choice-universalUnfolding}. We present the result for even dimension $n$, and leave the case of odd dimensions to the reader.

\subsection{Critical points and the Milnor ring}
Suppose $n\geq 4$ even.
We regard $W_0$ as a map from $Z^n\times \mathbb{C}^{\times}\rightarrow \mathbb{C}$.
Define the Milnor ring to be
\begin{equation*}
      \mathrm{MR}_n=\mathbb{C}[z_1,\dots,z_n,\frac{1}{z_1\cdots z_n-1},q]/(\frac{\partial W_0}{\partial z_1},\dots,\frac{\partial W_0}{\partial z_n}).
\end{equation*}

\begin{proposition}
The critical locus of $W_0:Z^n\times \mathbb{C}^{\times}\rightarrow \mathbb{C}$ is $V(I_{1})\sqcup V(I_{2})$ where
\begin{equation*}
      I_{1}=\sum_{i=1}^{\frac{n}{2}-1}\big(z_{2i+1}-1\big)+\sum_{i=1}^{\frac{n}{2}}(z_{2i}-z_n)+(qz_1^2-1)+(z_1 z_n^{\frac{n}{2}}-2),
\end{equation*}
and
\begin{equation*}
      I_{2}=\sum_{i=1}^{\frac{n}{2}-1}\big(z_{2i+1}-(-1)^i\big)+\sum_{i=1}^{\frac{n}{2}}(z_{2i})+\big(qz_1^2-(-1)^{\frac{n-2}{2}}\big).
\end{equation*}
In particular, for any fixed $q\in \mathbb{C}^{\times}$, $W_0$ has $n+2$ critical points, which are all nondegenerate.
Moreover, the critical values are
\begin{equation*}
      W_0(p_j)=\begin{dcases}
      n \zeta^{j}(4q)^{\frac{1}{n}},& \mbox{for}\ 1\leq j\leq n,\\
      0,& \mbox{for}\ j=n+1,n+2.
      \end{dcases}
\end{equation*}

\end{proposition}
\begin{proof}
We have
\begin{equation}\label{eq-derivatives-W0-zj}
\frac{\partial W_0}{\partial z_j}=
\begin{dcases}
\frac{q z_1 z_n(\prod_{i=1}^{n}z_i-2)}{(\prod_{i=1}^{n}z_i-1)^2},& j=1 \\
-\frac{q z_1^2 z_{n}\cdot z_1\prod_{i=3}^n z_i }{(\prod_{i=1}^{n}z_i-1)^2}+1+z_3,& j=2\\
-\frac{q z_1^2 z_{n}\cdot \prod_{i=1}^{j-1}z_i\cdot \prod_{i=j+1}^n z_i }{(\prod_{i=1}^{n}z_i-1)^2}
+z_{j-1}+z_{j+1},& 3\leq j\leq n-1,\\
-\frac{q z_1^2 }{(\prod_{i=1}^{n}z_i-1)^2}+z_{n-1},& j=n.
\end{dcases}
\end{equation}

Fix a (primitive) $n$-th root of unity $\zeta$.
The critical points are $p_1,\dots,p_{n},p_{n+1},p_{n+2}$, whose coordinates are
\begin{equation}\label{eq-zCoordinates-criticalPoints}
      \begin{cases}
 \mbox{for $1\leq j\leq n$}:\
z_1(p_j)=\zeta^{-\frac{jn}{2}}q^{-\frac{1}{2}},
z_{2i+1}(p_j)=1\ \mbox{for}\ 1\leq i\leq\frac{n}{2}-1,\\
\hspace{2.5cm} z_{2i}(p_j)=\zeta^{j}(4q)^{\frac{1}{n}},\ \mbox{for}\ 1\leq i\leq\frac{n}{2},\\
z_1(p_{n+1})= (\sqrt{-1})^{\frac{n-2}{2}}q^{-\frac{1}{2}},
z_{2i+1}(p_{n+1})=(-1)^i\ \mbox{for}\ 1\leq i\leq\frac{n}{2}-1,\\
\hspace{2.5cm} z_{2i}(p_{n+1})=0,\ \mbox{for}\ 1\leq i\leq\frac{n}{2},\\
z_1(p_{n+2})= -(\sqrt{-1})^{\frac{n-2}{2}}q^{-\frac{1}{2}},
z_{2i+1}(p_{n+2})=(-1)^i\ \mbox{for}\ 1\leq i\leq\frac{n}{2}-1,\\
\hspace{2.5cm} z_{2i}(p_{n+2})=0,\ \mbox{for}\ 1\leq i\leq\frac{n}{2}.
\end{cases}
\end{equation}
Note that in $\mathrm{MR}_n$ we have
\begin{equation}\label{eq-W0=nz2}
      W_0=n z_2.
\end{equation}
So the critical values follow.
\end{proof}

By evaluating the functions in $\mathbb{C}[z_1,\dots,z_n,\frac{1}{z_1\cdots z_n-1},q]$ at the critical points, one can find their relations in $\mathrm{MR}_n$. 
By Vandermonde, $\{z_2^i\}_{0\leq i\leq n}\cup\{z_1\}$ is a basis of $\mathrm{MR}_n$, and
\begin{equation}\label{eq-quadric-mirror-jac-relations}
      z_2^{n+1}=4q z_2,\ z_1 z_2=\frac{z_2^{\frac{n}{2}+1}}{2q},\ 
      z_1^2=\begin{cases}
      q^{-1},\ \mbox{if}\ n\equiv 2 \mod 4,\\
      -q^{-1}+\frac{z_2^n}{2q^2},\ \mbox{if}\ n\equiv 0 \mod 4.
      \end{cases}
\end{equation}
Let
\begin{equation}\label{eq-mirrorLGofQuadric-mirrorBasisToSmallQH}
      f_1=z_2,\  f_2=(\sqrt{-1})^{\frac{n}{2}}(2q z_1-z_2^{\frac{n}{2}}).
\end{equation}
Then in $\mathrm{MR}_n$,
\begin{equation}
      f_1^{n+1}=4q f_1,\ f_1f_2=0,\ f_2^2=-4q+f_1^n.
\end{equation}

\begin{corollary}\label{cor-milnorRing-smallQuantumCohOfQn}
Let $\mathbbm1$ be the identity element of the cohomology of $Q^n$.
Let $\sfh\in H^2(Q^n;\mathbb{Z})$ be the hyperplane class of the quadric hypersurface $Q^n$.
Let $\sfp\in H^{n}(Q^n;\mathbb{Z})$ be a primitive cohomology class of $Q^n$ satisfying $\int_{Q^n}\sfp\cup\sfp=2$.
Then the map $1\mapsto \mathbbm1$, $f_1\mapsto \sfh$, $f_2\mapsto \sfp$  gives an isomorphism of $\mathbb{C}[q]$-algebras from $\mathrm{MR}_n$ to the small quantum cohomology ring of the quadric $Q^n$.
\end{corollary}
\begin{proof}
The small quantum multiplication of $Q^n$ is given by \cite[Prop. 1]{Bea95}.
One easily checks that the given map matches  the small quantum products with the products in the Milnor ring.
\end{proof}

For later use, we recall the 3-points genus 0 Gromov-Witten invariants of $Q^n$ computed from \cite[Prop. 1]{Bea95} and the Poincaré pairings $\int_{Q^n}\sfh^i\cup\sfh^j=2\delta_{i+j,n}$, $\int_{Q^n}\sfh^i\cup\sfp=0$, and $\int_{Q^n}\sfp\cup\sfp=2$.
\begin{lemma}[{\cite[Lemma 5.7]{HK21}}]\label{lem-nonzerothreepointinvariants-quadric}
The non-zero 3-point genus 0 Gromov-Witten invariants of $Q^n$ are the following:
      \begin{equation}\label{eq-nonzerothreepointinvariants-quadric}
      \begin{cases}
      \langle \sfp,\sfp, \mathbbm1\rangle=2,\ \langle \sfp,\sfp, \sfh_{n}\rangle=-4,\\
      \langle \sfh_{i},\sfh_{j},\sfh_{k}\rangle=2,\ \mbox{if}\ i+j+k=n,\\
      \langle \sfh_{i},\sfh_{n-i},\sfh_{n}\rangle=4,\ \mbox{if}\ 1\leq i\leq n-1,\\
      \langle \sfh_{i},\sfh_{j},\sfh_{k}\rangle=8,\ \mbox{if}\ i+j+k=2n,\ \mbox{and}\ i,j,k\neq n,\\
      \langle \sfh_{n},\sfh_{n},\sfh_{n}\rangle=8.
      \end{cases}
      \end{equation}
Here we omit the (stable map) degree of a Gromov-Witten invariant because it is  determined by the dimension constraint of Gromov-Witten invariants and that $Q^n$ has a positive Fano index.\pqed
\end{lemma}

Recall that (\cite[Prop. 2.35]{Gro11}) $\mathrm{Hess}(W_0)(p)$ is defined to be the determinant of the Hessian matrix $(\frac{\partial^2 W_0}{\partial u_i\partial u_j})_{1\leq i,j\leq n}$ at $p$, where $u_1,\dots,u_n$ are local coordinates near $p$ such that $\Omega=\mathrm{d}u_1\dots \mathrm{d}u_n$.
For later use we need to compute the Hessian of $W_0$ at the critical points.
\begin{lemma}\label{lem-hessianValue-quadric}
\begin{equation}\label{eq-hessianValue-quadric}
      \mathrm{Hess}(W_0)(p_i)=\begin{cases}
      2n q,& \mbox{if}\ 1\leq i\leq n,\\
      -4q,& \mbox{if}\ i=n+1,n+2.
      \end{cases}
\end{equation}
\end{lemma}
\begin{proof}
Near $p_k$, for $1\leq k\leq n$, let $x=e^{\theta_1}$, $y_i=e^{\theta_{i+1}}$ for $1\leq i\leq n-1$, so that $\Omega=\prod_{i=1}^{n}\mathrm{d}\theta_i$. Then
\begin{eqnarray*}
      &z_i=\frac{e^{\theta_1}+1}{e^{\theta_2+\theta_4+\dots+\theta_n}},\
      z_2=e^{\theta_2},\ z_3=e^{\theta_3-\theta_2},\ z_4=e^{\theta_4+\theta_2-\theta_3},\dots
\end{eqnarray*}
We have thus
\begin{equation*}
\frac{\partial^2 W_0}{\partial \theta_1\partial \theta_1}(p_k)=\frac{2q}{(\zeta^{k}(4q)^{\frac{1}{n}})^{n-1}}
=\frac{2\zeta^{k}q}{(4q)^{\frac{n-1}{n}}}=\frac{\zeta^k (4q)^{\frac{1}{n}}}{2},
\end{equation*}
\begin{equation*}
\frac{\partial^2 W_0}{\partial \theta_1\partial \theta_j}(p_k)=0,\ 2\leq j\leq n,
\end{equation*}
\begin{equation*}
\frac{\partial^2 W_0}{\partial \theta_j\partial \theta_l}(p_k)=(1+\delta_{j,l})\zeta^{k}(4q)^{\frac{1}{n}},\ 2\leq j\neq l\leq n.
\end{equation*}
Then we get $\mathrm{Hess}(W_0)(p_i)=2n q$ for $1\leq i\leq n$.
Near $p_{n+1}$, first we set
\[
\tilde{w}_1=z_1-(\sqrt{-1})^{\frac{n-2}{2}}q^{-\frac{1}{2}},\
w_{2i+1}=z_{2i+1}-(-1)^i\ \mbox{for}\ 1\leq i\leq\frac{n}{2}-1,
\]
\[
w_{2i}=0,\ \mbox{for}\ 1\leq i\leq\frac{n}{2}.
\]
Thus
\[
\Omega=\frac{\mathrm{d}\tilde{w}_1\cdots \mathrm{d}w_{n}}{\big(\tilde{w}_1+(\sqrt{-1})^{\frac{n-2}{2}}q^{-\frac{1}{2}}\big)
\prod_{i=1}^{n/2-1}\big(w_{2i+1}+(-1)^i\big)\prod_{i=1}^{n/2}w_{2i}-1}.
\]
Then we define
\begin{equation*}
w_1=\frac{\log\Big(1+\frac{\prod_{i=1}^{n/2-1}\big(w_{2i+1}+(-1)^i\big)\prod_{i=1}^{n/2}w_{2i}}
{(\sqrt{-1})^{\frac{n-2}{2}}q^{-\frac{1}{2}}\prod_{i=1}^{n/2-1}\big(w_{2i+1}+(-1)^i\big)\prod_{i=1}^{n/2}w_{2i}-1}\tilde{w}_1\Big)}{\prod_{i=1}^{n/2-1}\big(w_{2i+1}+(-1)^i\big)\prod_{i=1}^{n/2}w_{2i}},
\end{equation*}
so that 
\[
\Omega=\mathrm{d}w_1\cdots \mathrm{d}w_{n}.
\]
We have
\[
\frac{\partial^2 W_0}{\partial z_r \partial z_s}(p_{n+1})=\begin{cases}
-2(\sqrt{-1})^{\frac{n-2}{2}}q^{\frac{1}{2}},& \mbox{if}\ (r,s)=(1,n)\ \mbox{or}\ (n,1);\\
1,& \mbox{if}\ 2\leq r,s\leq n\ \mbox{and}\ |r-s|=1;\\
0,& \mbox{otherwise}.
\end{cases}
\]
Since
\begin{multline}\label{eq-def-w1}
w_1=\frac{\log\big(1+\frac{\prod_{i=2}^{n}z_i}{(\sqrt{-1})^{\frac{n-2}{2}}q^{-\frac{1}{2}}\prod_{i=2}^{n}z_i-1}(z_1-(\sqrt{-1})^{\frac{n-2}{2}}q^{-\frac{1}{2}})\big)}{\prod_{i=2}^{n}z_i}\\
= \frac{z_1-(\sqrt{-1})^{\frac{n-2}{2}}q^{-\frac{1}{2}}}{(\sqrt{-1})^{\frac{n-2}{2}}q^{-\frac{1}{2}}\prod_{i=2}^{n}z_i-1}
-\frac{\prod_{i=2}^{n}z_i (z_1-(\sqrt{-1})^{\frac{n-2}{2}}q^{-\frac{1}{2}})}{2\big((\sqrt{-1})^{\frac{n-2}{2}}q^{-\frac{1}{2}}\prod_{i=2}^{n}z_i-1\big)^2}
+O\big(\sum_{i=1}^n|z_i-z_i(p_{n+1})|^2\big),
\end{multline}
we have
\[
\frac{\partial w_i}{\partial z_j}(p_{n+1})=\begin{cases}
-1,& \mbox{if}\ i=j=1,\\
1,& \mbox{if}\ 2\leq i=j\leq n,\\
0,& \mbox{otherwise}.
\end{cases}
\]
Hence
\[
\frac{\partial^2 W_0}{\partial w_j \partial w_l}(p_{n+1})=\begin{cases}
2(\sqrt{-1})^{\frac{n-2}{2}}q^{\frac{1}{2}},& \mbox{if}\ (j,l)=(1,n)\ \mbox{or}\ (n,1);\\
1,& \mbox{if}\ 2\leq j,l\leq n\ \mbox{and}\ |j-l|=1;\\
0,& \mbox{otherwise}.
\end{cases}
\]
It follows that $\mathrm{Hess}(W_0)(p_{n+1})= -4q$. Similarly we have $\mathrm{Hess}(W_0)(p_{n+2})= -4q$.
\end{proof}

\subsection{Universal unfolding and Euler field}\label{sec:universalUnfolding}

Following Mark Gross, we put $q=e^{t_1}$. Then we set
\begin{eqnarray}\label{eq-unfolding-LGMirrorOfQuadric-final}
W=t_0+W_0+\sum_{i=2}^{n}t_i W_0^i-2e^{t_1}t_n+t_{n+1} e^{t_1}z_1: Z^n\times \mathbb{C}^{n+2}\rightarrow \mathbb{C}.
\end{eqnarray}
We define the Euler vector field on the moduli space $\mathbb{C}^{n+2}$ to be
\begin{equation}\label{eq-EulerField-final}
      E=n\partial_{t_1}+\sum_{i=0}^{n}(1-i)t_i\partial_{t_i}+(1-\frac{n}{2})t_{n+1}\partial_{t_{n+1}}.
\end{equation}
Then $E$ has a lift to the total space $Z^n\times \mathbb{C}^{n+2}$
\begin{equation*}
      \tilde{E}=n\partial_{t_1}+\sum_{i=0}^{n}(1-i)t_i\partial_{t_i}+(1-\frac{n}{2})t_{n+1}\partial_{t_{n+1}}
      -\frac{n}{2}z_{1}\partial_{z_1}
      +\sum_{i=1}^{\frac{n}{2}} z_{2i} \partial_{z_{2i}},
\end{equation*}
satisfying
\begin{equation}\label{eq-EulerField-liftingProperty-final}
      \tilde{E}(W)=W.
\end{equation}

\begin{remark}\label{rem:choice-universalUnfolding}
A priori, we have the following requirements  (see \cite[proof of Prop. 2.38]{Gro11}, and Proposition \ref{eq-semiInfiniteHodge-mirrorQuadric}) for a universal unfolding $W$:
\begin{enumerate}
      \item The evaluation of the derivatives $\frac{\partial W}{\partial t_i}$, for $1\leq i\leq n$, at $\mathbf{t}=0$, span the Milnor ring $\mathrm{MR}_n$.
      \item The exists a lift $\tilde{E}$ satisfying (\ref{eq-EulerField-liftingProperty-final}).
\end{enumerate}
Such $W$ is not unique. For example, for $C_1,C_2\in \mathbb{C}$, and $a_0,a_1\in \mathbb{Z}_{\geq 0}$ such that $a_1\geq 1$ and $a_1\geq a_0$, we  take
\begin{eqnarray*}
W=t_0+W_0+\sum_{i=2}^{n}t_i W_0^i+C_1e^{t_1}t_n+C_2t_{n+1}e^{a_0 t_1}z_1^{2a_1-1} z_2^{n(a_1-a_0)}.
\end{eqnarray*}
Then (i) is satisfied, and the same form of $\tilde{E}$ as above makes (ii) satisfied. We choose $a_0=a_1=C_2=1$ from our partiality. Our choice $C_2=-2$ arises from that
the $n$-th small quantum power of the hyperplane class $\sfh$ of $Q^n$ is equal to $\sfh^n-2$ (from e.g. Lemma \ref{lem-nonzerothreepointinvariants-quadric}), and thus we have a slightly simpler form of the identification of the A-side and B-side moduli spaces than the choice $C_2=0$.
\end{remark}

\section{Preparations for the B-model Frobenius manifolds}\label{sec:preparation-BModel-FrobeniusManifold}
\subsection{Construction of Frobenius manifolds}
We recall the Barannikov(-Kontsevich)   construction  (\cite{BK98}, \cite{Bar00}, \cite{Bar01}) of the Frobenius manifold associated with an LG model, which can be regarded as a generalization of K. Saito \cite{Sai83}. We follow a slight modification of it by M. Gross and his terminology \cite[Chap. 2]{Gro11}; see also \cite{Iri09}.
\begin{construction}\label{cons-Barannikov}\phantom\\
\begin{enumerate}
      \item[(i)] First we need a \emph{semi-infinite variation of Hodge structure} (SVHS for short) with a \emph{grading}. We refer the reader to \cite[Definition 2.20]{Gro11} for the definition. We only recall that a SVHS consists of a space $\mathcal{M}$, a locally free $\mathcal{O}_{\mathcal{M}}\{\hbar\}$-module $\mathcal{E}$ of finite rank, a flat connection $\nabla:\mathcal{E}\rightarrow \Omega_{\mathcal{M}}^1\otimes \hbar^{-1}\mathcal{E}$, and a pairing $(\cdot,\cdot)_{\mathcal{E}}:\mathcal{E}\times \mathcal{E}\rightarrow \mathcal{O}_{\mathcal{M}}\{\hbar\}$ satisfying certain conditions. A grading is a $\mathbb{C}$-linear endomorphism $\mathrm{Gr}$ of $\mathcal{E}$ satisfying certain conditions.
      \item[(ii)] Suppose that $\mathcal{M}$ is simply connected, and $\big(\mathcal{E},\nabla,(\cdot,\cdot)_{\mathcal{E}}\big)$ is a SVHS over $\mathcal{M}$. Let 
      \[
      \mathcal{H}:=\{s \in \Gamma(\mathcal{M}, \mathcal{E}\otimes_{\mathcal{O}_{\mathcal{M}}\{\hbar\}}\mathcal{O}_{\mathcal{M}}\{\hbar,\hbar^{-1}\})| \nabla s = 0\}.
      \]
      By parallel transport we get a family of subspace $\mathcal{E}_x$ of $\mathcal{H}$ parametrized by $\mathcal{M}$. Note that $\mathcal{H}$ is a $\mathbb{C}\{\hbar,\hbar^{-1}\}$-module, and every $\mathcal{E}_x$ is a $\mathbb{C}\{\hbar\}$-submodule.
      \item[(iii)] Fix a point $0\in \mathcal{M}$. We choose an $\mathcal{O}(\mathbb{P}^1\backslash\{0\})$-submodule $\mathcal{H}_{-}$ of $\mathcal{H}$ satisfying that $\mathcal{H}_{-}\oplus \mathcal{E}_0\rightarrow \mathcal{H}$ is an isomorphism. One can show that the projections 
      $\mathcal{E}_0\cap \hbar \mathcal{H}_{-}\rightarrow \mathcal{E}_0/\hbar \mathcal{E}_0$ and 
      $\mathcal{E}_0\cap \hbar \mathcal{H}_{-}\rightarrow \hbar\mathcal{H}_{-}/\mathcal{H}_{-}$ are isomorphisms.
      By (ii) we regard $\mathcal{E}$ as a subbundle of the trivial bundle $\mathcal{H}_{\mathcal{M}}:=\mathcal{H}\times \mathcal{M}$. The projection
      $\mathcal{E}\cap \hbar \mathcal{H}_{-,\mathcal{M}}\rightarrow (\hbar\mathcal{H}_{-}/\mathcal{H}_{-})_{\mathcal{M}}$ is an isomorphism at 0, and thus is an isomorphism in an open neighborhood $\mathcal{U}$ of 0. We denote by $\tau$ the induced trivialization 
      \begin{equation*}
            \tau:(\hbar\mathcal{H}_{-}/\mathcal{H}_{-})\otimes_{\mathbb{C}}\mathcal{O}_{\mathcal{U}}\{\hbar\}\rightarrow \mathcal{E}|_{\mathcal{U}}.
      \end{equation*}
      The connection $\nabla$ induces a connection $\nabla$ on $:(\hbar\mathcal{H}_{-}/\mathcal{H}_{-})\otimes_{\mathbb{C}}\mathcal{O}_{\mathcal{U}}\{\hbar\}$.
      Then we  write this connection in the trivialization $\tau$ as $\nabla=d+\hbar^{-1}A$, where A is an $\mathrm{End}(\hbar \mathcal{H}_{-}/\mathcal{H}_{-})$-valued 1-form on $\mathcal{U}$. For a vector field $X$ on $\mathcal{U}$, let $A_X:=X\lrcorner A\in \mathrm{End}(\hbar \mathcal{H}_{-}/\mathcal{H}_{-})_{\mathcal{U}}$.
      \item[(iv)] Define a symplectic form $\overline{\Omega}(\cdot,\cdot)$ on $\mathcal{H}$ by $\overline{\Omega}(s_1, s_2)=$  the coefficient of $\hbar^{-1}$ in $(s_1, s_2)_{\mathcal{E}}$. Suppose given a grading $\mathrm{Gr}$ of the SVHS $\big(\mathcal{E},\nabla,(\cdot,\cdot)_{\mathcal{E}}\big)$. We furthermore require that the opposite space $\mathcal{H}_{-}$ is isotropic with respect to $\overline{\Omega}$, and is preserved by $\mathrm{Gr}$. Then $\mathrm{Gr}$ induces an endomorphism $\mathrm{Gr}_0$ of $\hbar\mathcal{H}_{-}/\mathcal{H}_{-}$.
      \item[(v)] Suppose given $\Omega_0\in \hbar \mathcal{H}_{-}$, satisfying 
            \begin{itemize}
                  \item $[\Omega_0]$ is an eigenvector of $\mathrm{Gr}_0$;
                  \item the section $s_0':=\tau([\Omega_0]\otimes 1)$ satisfies \emph{miniversality}: the map
                  \begin{equation*}
                        \mathcal{T}_{\mathcal{U}}\rightarrow \mathcal{E}/\hbar \mathcal{E},\
                        X\mapsto \hbar\nabla_{X}s'_0
                  \end{equation*}
                  is an isomorphism.
            \end{itemize}
            Then we have an induced \emph{period map},
            \begin{equation*}
                  \Psi:\mathcal{U}\rightarrow \mathcal{H},\
                  q\mapsto (\Omega_0+\mathcal{H}_{-})\cap \mathcal{E}_q,
            \end{equation*}
            and a map 
            \begin{equation*}
                  \psi:\mathcal{U}\rightarrow \hbar \mathcal{H}_{-}/\mathcal{H}_{-},\
                  q\mapsto [\hbar(\Psi(q)-\Omega_0)].
            \end{equation*}
            Miniversality of $s'_0$ implies that $\psi$ is an local isomorphism. So $\psi$ induces an isomorphism $\mathcal{T}_{\mathcal{U}}\xrightarrow{\sim} (\hbar \mathcal{H}_{-}/\mathcal{H}_{-})_{\mathcal{U}}$, and we have $\psi_*(X)= A_X([\Omega_0])$. The trivial flat connection on 
            $(\hbar\mathcal{H}_{-}/\mathcal{H}_{-})_{\mathcal{U}}$ induces a flat connection $\nabla^{\mathcal{U}}$ on $\mathcal{T}_{\mathcal{U}}$.
      \item[(vi)] That $\mathcal{H}_{-}$ is isotropic with respect to $\overline{\Omega}$ enables us to define a non-degenerate symmetric bilinear pairing on $\hbar\mathcal{H}_{-}/\mathcal{H}_{-}$ by $(s_1,s_2)_{\hbar\mathcal{H}_{-}/\mathcal{H}_{-}}:=(s_1,s_2)_{\mathcal{E}}|_{\hbar=\infty}$. It induced a pairing $g$ on $\mathcal{T}_{\mathcal{U}}$.
      \item[(vii)] We define a product on $\mathcal{T}_{\mathcal{U}}$ by $A_{X_1\circ X_2}[\Omega_0]=A_{X_1} A_{X_2}[\Omega_0]$. We define
      \begin{equation*}
            \mathcal{A}:S^3 \mathcal{T}_{\mathcal{U}}\rightarrow \mathcal{O}_{\mathcal{U}},\
            \mathcal{A}(X_1,X_2,X_3)=g(X_1\circ X_2,X_3).
      \end{equation*}
      We define $e$ to be the vector field on $\mathcal{U}$ by requiring $A_{e}[\Omega_0]=[\Omega_0]$.
      Then $(\nabla^{\mathcal{U}},g,\mathcal{A})$ gives a Frobenius structure on $\mathcal{U}$ such that $e$ is the flat identity.
\end{enumerate}
\end{construction}

\subsection{Semi-infinite variation of Hodge structures}
\begin{definition}
Let $\wcM$ be the ringed space $(\mathbb{C},\mathcal{O}_{\wcM})$, where $\mathbb{C}$ is the complex manifold $\mathbb{C}$ with coordinate $t_1$, and for any open subset $U\subset \mathbb{C}$, 
\[
\mathcal{O}_{\wcM}(U)=\big\{\mbox{formal power series}\ \sum_{i_0,i_2,i_3,\dots,i_{n+1}\in \mathbb{Z}_{\geq 0}}f_{i_0,i_2,i_3,\dots,i_{n+1}}t_0^{i_0}t_2^{i_2}t_3^{i_3}\cdots t_{n+1}^{i_{n+1}}| f_{i_0,i_2,i_3,\dots,i_{n+1}}\in \mathcal{O}_{\mathbb{C}}(U)
\big\}.
\]
\end{definition}
We can regard $\wcM$ as the formal completion of  the analytic space $\mathbb{C}^{n+2}$ with coordinates $t_0,\dots,t_{n+1}$ along $\{t_0=t_2=t_3=\dots=t_{n+1}=0\}$.
Recall that  $Z^n=\{(z_1,\dots,z_n)\in \mathbb{C}^n| z_1\cdots z_n\neq 1\}$. Define $\mathcal{Z}^n=\wcM\times Z^n$, and $\pi: \mathcal{Z}^n\rightarrow \wcM$ the obvious projection. According Section \ref{sec:universalUnfolding}, we define $W\in \Gamma(\mathcal{Z}^n,\mathcal{O}_{\mathcal{Z}^n})$ by
\begin{eqnarray*}
W(t_0,t_1,\dots,t_{n+1},z_1,\dots,z_n)
=t_0+W_0+\sum_{i=2}^{n}t_i W_0^i-2e^{t_1}t_n+t_{n+1} e^{t_1}z_1.
\end{eqnarray*}
On $\mathbb{C}^{\times}$ we fix a usual coordinate $\hbar$.
We define a local system $R$ on  $\wcM\times \mathbb{C}^{\times}$   such that the fiber over $(t_1,\hbar)$ is 
\begin{equation*}
      H_n\big(\pi^{-1}(t_1),\mathrm{Re}(\frac{W|_{\pi^{-1}(t_1)}}{\hbar})\ll 0\big)
      :=\varprojlim_{c\rightarrow -\infty}H_n\big(\pi^{-1}(t_1),\mathrm{Re}(\frac{W|_{\pi^{-1}(t_1)}}{\hbar})<c\big).
\end{equation*}
Here $W|_{\pi^{-1}(t_1)}$ is understood as $W(0,t_1,0,\dots,0,z_1,\dots,z_n)=W_0(z_1,\dots,z_n,e^{t_1})$. Then we define $\mathcal{R}=R\otimes_{\mathbb{C}}\mathcal{O}_{\wcM\times \mathbb{C}^{\times}}$, and $\mathcal{R}^{\vee}$ the dual of $\mathcal{R}$. The locally free sheaf $\mathcal{R}$ carries a flat connection $\nabla^{\mathrm{GM}}$ whose flat sections are local sections of  $R$. The dual connection on $\mathcal{R}^{\vee}$ is also denoted by $\nabla^{\mathrm{GM}}$.


We define a sheaf $\Omega_{\wcM}\{\hbar\}$ exactly as in \cite[Definition 2.17]{Gro11} on $\wcM$, and the locally free $\Omega_{\wcM}\{\hbar\}$-module $\mathcal{E}$ exactly as in \cite[Definition 2.37]{Gro11}. Then $\mathcal{E}$ carries a connection $\nabla$ induced by $\nabla^{\mathrm{GM}}$.  Denote by $(-)$ the map $\wcM\times \mathbb{C}^{\times}\rightarrow \wcM\times \mathbb{C}^{\times}$, $(u,\hbar)\mapsto (u,-\hbar)$. There is a natural perfect intersection pairing (\cite[\S 5]{Pha85})
\[
H_n\big(\pi^{-1}(t_1),\mathrm{Re}(\frac{W|_{\pi^{-1}(t_1)}}{\hbar})\ll 0\big)\times 
H_n\big(\pi^{-1}(t_1),\mathrm{Re}(-\frac{W|_{\pi^{-1}(t_1)}}{\hbar})\ll 0\big)\rightarrow \mathbb{C},
\]
which induces an intersection pairing $(-)^* \mathcal{R}^{\vee}\times \mathcal{R}^{\vee}\rightarrow \mathcal{O}_{\wcM\times \mathbb{C}^{\times}}$, denoted by $(\cdot,\cdot)$. Then we define the pairing
\[
(s_1,s_2)_{\mathcal{E}}(\hbar)=\frac{(-1)^{\frac{n(n+1)}{2}}}{(2\pi\sqrt{-1}\hbar)^n}\big((-)^* s_1,s_2\big).
\]
We define the Euler vector field on $\wcM$ to be
\begin{equation*}
      E=n\partial_{t_1}+\sum_{i=0}^{n}(1-i)t_i\partial_{t_i}+(1-\frac{n}{2})t_{n+1}\partial_{t_{n+1}}.
\end{equation*}
We refer the reader to \cite[Definition 2.20]{Gro11} for the definition of semi-infinite Hodge structure and the grading operators.
\begin{proposition}\label{eq-semiInfiniteHodge-mirrorQuadric}
The data $\mathcal{E},\nabla$ and $(\cdot,\cdot)_{\mathcal{E}}$ yield a semi-infinite Hodge structure parametrized by $\wcM$. The operator $\mathrm{Gr}$, defined by $\mathrm{Gr}(s)=\nabla^{\mathrm{GM}}_{\hbar \partial_\hbar+E}(s)-s$, is a grading operator with $D=2-n$.
\end{proposition}
\begin{proof}
Note that by our construction in Section \ref{sec:universalUnfolding} we can lift the vector field $E$ to a vector field $\tilde{E}$
on $\mathcal{Z}^n$, satisfying $\tilde{E}(W)=W$. Then the proof of \cite[Prop. 2.38]{Gro11} carries over almost verbatim. 
\end{proof}

For computations of the Gauss-Manin connection, we need the following description of fibers of $\mathcal{R}^{\vee}$. Fix $(t_1,\hbar)$. Let $f$ be an algebraic function on $Z^n$, then for $\Gamma\in H_n\big(\pi^{-1}(t_1),\mathrm{Re}(\frac{W|_{\pi^{-1}(t_1)}}{\hbar})\ll 0\big)$, the integration
\begin{equation}\label{eq-oscilatoryIntegral}
\int_\Gamma e^{\frac{W|_{\pi^{-1}(t_1)}}{\hbar}}f \Omega
\end{equation}
is well-defined. We get a pairing
\begin{equation}\label{eq-pairing-oscilatoryIntegral}
      H_n\big(\pi^{-1}(t_1),\mathrm{Re}(\frac{W|_{\pi^{-1}(t_1)}}{\hbar})\ll 0\big)
      \times H^n\big( \Gamma(\Omega_{\pi^{-1}(t_1)}^{\bullet},\mathrm{d}+\frac{\mathrm{d}W|_{\pi^{-1}(t_1)}}{\hbar}\wedge)\big)\rightarrow \mathbb{C}.
\end{equation}
\begin{lemma}\label{lem-pairing-oscilatoryIntegral}
The pairing (\ref{eq-pairing-oscilatoryIntegral}) is a perfect pairing.
\end{lemma}
\begin{proof}
By the stationary phase approximation (\cite[Prop. 2.35]{Gro11}) and computations in Section \ref{sec:pairings}, the induced map 
\[
H_n\big(\pi^{-1}(t_1),\mathrm{Re}(\frac{W|_{\pi^{-1}(t_1)}}{\hbar})\ll 0\big)
\rightarrow H^n\big( \Gamma(\Omega_{\pi^{-1}(t_1)}^{\bullet},\mathrm{d}+\frac{\mathrm{d}W|_{\pi^{-1}(t_1)}}{\hbar}\wedge)\big)^{\vee}
\]
is injective. It thus remains to show that both sides have the same rank. Since LHS has rank $n+2$, which is equal to the rank of the Jacobi ring of $W_0$, it suffices to show 
\begin{equation}\label{eq-rank-twistedDeRham}
h^n\big( \Gamma(\Omega_{\pi^{-1}(t_1)}^{\bullet},\mathrm{d}+\frac{\mathrm{d}W|_{\pi^{-1}(t_1)}}{\hbar}\wedge)\big)
=h^n\big( \Gamma(\Omega_{\pi^{-1}(t_1)}^{\bullet},\mathrm{d}W|_{\pi^{-1}(t_1)}\wedge)\big).
\end{equation}
Since $\pi^{-1}(t_1)\cong Z^n$ is affine, one has
\[
\mathbb{H}^n\big(\pi^{-1}(t_1),(\Omega_{\pi^{-1}(t_1)}^{\bullet},\mathrm{d}+\frac{\mathrm{d}W|_{\pi^{-1}(t_1)}}{\hbar}\wedge)\big)
\cong H^n\big( \Gamma(\Omega_{\pi^{-1}(t_1)}^{\bullet},\mathrm{d}+\frac{\mathrm{d}W|_{\pi^{-1}(t_1)}}{\hbar}\wedge)\big).
\]

By Corollary \ref{cor-compatifiedLGModel-mirrorQuadric} we have a projective Landau-Ginzburg model $(\overline{Z},W_0)$ such that the critical locus is contained in $Z^n$, and $V=\overline{Z}^n\backslash Z^n$ is a normal crossing divisor. 
In other words $\mathrm{d}W_0$ is locally non-vanishing in $\Omega_{\overline{Z}^n}^1(\log V)$. 
So (\ref{eq-rank-twistedDeRham}) follows from its log version in \cite{Sab99} (see the arguments in \cite[Page 129]{Gro11}).
\end{proof}
In our indirect computation of broad periods we need the following.
\begin{lemma}\label{lem-algebraicity-restrictionOfcalE}
Let $\wcM'=\Spec(\mathbb{C}[t_1])\times \Spec(\mathbb{C}[\hbar,\hbar^{-1}])$, and $\mathcal{Y}=X\times\wcM'$. Let $\mathrm{d}=\mathrm{d}_{\mathcal{Y}/\wcM'}$ denote the relative differential. Then the sheaf associated with the cohomology module 
\begin{equation*}
      H^n\big(\Gamma(\Omega_{\mathcal{Y}/\wcM',\mathrm{Zar}}^{\bullet},\mathrm{d}+\frac{\mathrm{d}W|_{\wcM'}}{\hbar}\wedge)\big)
\end{equation*}
is a locally free coherent sheaf on $\Spec(\mathbb{C}[t_1])\times \Spec(\mathbb{C}[\hbar,\hbar^{-1}])$.
\end{lemma}
\begin{proof}
 The complex $ \Gamma(\Omega_{\mathcal{Y}/\wcM',\mathrm{Zar}}^{\bullet},\mathrm{d}+\frac{\mathrm{d}W|_{\wcM'}}{\hbar}\wedge)$ is $\mathcal{O}_{\wcM'}$-linear. The sheaf in the statement is thus a quasi-coherent sheaf. It then suffices to note that it has a constant finite rank, which is a consequence of Lemma \ref{lem-pairing-oscilatoryIntegral} and a comparison theorem between analytic and algebraic twisted de Rham cohomology (\cite{Sab99}, \cite[Theorem 2.31]{Gro11}).
\end{proof}

\section{Picard-Fuchs equations for the periods}\label{sec:Picard-Fuchs}
To obtain a Frobenius manifold from a semi-infinite variation of Hodge structures, we need to choose an $\mathcal{O}(\mathbb{P}^1\backslash\{0\})$-submodule $\mathcal{H}_{-}$ of $\mathcal{H}$ satisfying the requirements in Construction \ref{cons-Barannikov} (iii) and (iv). There is no canonical choice of $\mathcal{H}_{-}$ in general. In this section we compute the space of periods (\ref{eq-oscilatoryIntegral}) to find a natural one.

For any $t_1\in \mathbb{C}$, we identify $\pi^{-1}(t_1)$ with $Z^n$, and $W|_{\pi^{-1}(t_1)}$ with $W_0(z_1,\dots,z_n,q=e^{t_1})$.
Let $\Xi\in H_n\big(Z^n,\mathrm{Re}(\frac{W_0}{\hbar})\ll 0\big)$ be an arbitrary Lefschetz thimble. In this section we seek the Picard-Fuchs equation satisfied by the \emph{narrow period}
\begin{equation*}
      \psi=\int_{\Xi}e^{\frac{W_0}{\hbar}}\Omega
\end{equation*}
and the \emph{broad period}
\begin{equation*}
      \theta=\int_{\Xi}z_1e^{\frac{W_0}{\hbar}}\Omega.
\end{equation*}
Denote by $[f\Omega]$ the section of $\mathcal{R}^{\vee}$ defined by the integration $\int_\Xi e^{\frac{W}{\hbar}}f\Omega$. 
Recall that the Gauss-Manin connection satisfies
\begin{equation}\label{eq-GM-connection-hbar}
      \nabla_{\hbar\partial_\hbar}^{\mathrm{GM}}[f \Omega]=[(\hbar\partial_\hbar f-\hbar^{-1}W_0f)\Omega].
\end{equation}

\subsection{Picard-Fuchs equation for narrow periods}
\begin{proposition}\label{prop-PicardFuchs-Quadric}
\begin{equation}\label{eq-GaussManinConn-Quadric-narror}
(-\frac{1}{n}\nabla_{\hbar\partial_{\hbar}}^{\mathrm{GM}})^{n+1}[\Omega]
=   2q\hbar^{-n}(-\frac{2}{n}\nabla_{\hbar\partial_{\hbar}}^{\mathrm{GM}}[\Omega]+[\Omega]).
\end{equation}
Equivalently,
\begin{equation}\label{eq-PicardFuchs-Quadric}
      (-\frac{1}{n}\hbar\frac{\partial}{\partial \hbar})^{n+1}\psi=2q\hbar^{-n}(-\frac{2\hbar}{n}\frac{\partial}{\partial \hbar}+1)\psi.
\end{equation}
\end{proposition}

\begin{remark}\label{rem:proof-PicardFuchsEq-narrorPeriods}
If we restrict $[\Omega]$ to the domain $\{x_1 \prod_{i=1}^{n-1}y_i\neq 0\}$, the Picard-Fuchs equation is known (see e.g \cite{Iri11}). In dimension $n=3$, Gorbounov and Smirnov showed that one can use results in $D$-modules to reduce the computation to the domain $\{x_1 \prod_{i=1}^{n-1}y_i\neq 0\}$ (see \cite[Lemma 4.1.3]{GS15}. In the following proof, we adapt the proof of \cite[Lemma 4.2.3]{GS15} to all dimensions $n$, and replace the use of $D$-modules  with an elementary observation.
\end{remark}
 
We begin the proof by computing the Gauss-Manin connection in the domain $U=\{x_1 \prod_{i=1}^{n-1}y_i\neq 0\}$.  Recall 
 \begin{equation}\label{eq-restriction-W0-Omega-GiventalDomain}
       W_0|_{U}=\frac{q(x_1+1)^2}{x_1  \prod_{i=1}^{n-1}y_i}+y_1+\dots+y_{n-1},\
       \Omega|_{U}=\frac{dx_1dy_1\cdots dy_{n-1}}{x_1 y_1\cdots y_{n-1}}.
 \end{equation}

\begin{lemma}\label{lem-GaussManinConnection-yk}
Restricting on $U$, we have
\begin{eqnarray}\label{eq-GaussManinConnection-yk}
\nabla_{\hbar^2\partial_{\hbar}}^{\mathrm{GM}}[\prod_{k=1}^l y_k\Omega]=
\begin{cases}
 l\hbar[\prod_{k=1}^l y_k\Omega]-n[\prod_{k=1}^{l+1} y_k\Omega],& \mbox{if}\ 0\leq l<n-1,\\
 (n-1)\hbar[\prod_{k=1}^{n-1} y_k\Omega]-2nq[x_1\Omega]-2nq[\Omega],& \mbox{if}\ l=n-1.
 \end{cases}
\end{eqnarray}
\end{lemma}
\begin{proof}
 By (\ref{eq-GM-connection-hbar}) and (\ref{eq-restriction-W0-Omega-GiventalDomain}) we have
\begin{eqnarray*}
\nabla_{\hbar^2\partial_{\hbar}}^{\mathrm{GM}}[\prod_{k=1}^l y_k\Omega]=-[(\prod_{k=1}^l y_k) W_0\Omega]
=-[\big(\frac{q(x_1+1)^2}{x_1 \prod_{k=l+1}^{n-1} y_k}+(y_1+\dots+y_{n-1})\prod_{k=1}^l y_k\big)\Omega].
\end{eqnarray*}
For $1\leq j\leq l$,
\begin{eqnarray*}
&&(d+\frac{d W_0}{\hbar})\big((\prod_{k=1}^l y_k)y_j\iota_{\partial_{y_j}}\Omega\big)=
(\prod_{k=1}^l y_k)\Omega+\hbar^{-1} (\prod_{k=1}^l y_k) y_j \frac{\partial W_0}{\partial y_j}\Omega\nn\\
&=& \big(\prod_{k=1}^l y_k+\hbar^{-1}(-\frac{q(x_1+1)^2}{x_1  \prod_{i=l+1}^{n-1}y_i}+ y_j\prod_{k=1}^l y_k)\big)\Omega.
\end{eqnarray*}
For $l+1\leq j\leq n-1$,
\begin{eqnarray*}
&&(d+\frac{d W_0}{\hbar})\big((\prod_{k=1}^l y_k)y_j\iota_{\partial_{y_j}}\Omega\big)=
\hbar^{-1} (\prod_{k=1}^l y_k) y_j \frac{\partial W_0}{\partial y_j}\Omega\nn\\
&=& \hbar^{-1}(-\frac{q(x_1+1)^2}{x_1  \prod_{i=l+1}^{n-1}y_i}+ y_j\prod_{k=1}^l y_k)\Omega.
\end{eqnarray*}
It follows that for $1\leq l\leq n-1$,
\begin{eqnarray}\label{eq-proof-prop-PicardFuchs-Quadric-4}
\nabla_{\hbar^2\partial_{\hbar}}^{\mathrm{GM}}[\prod_{k=1}^l y_k\Omega]
&=& l\hbar[\prod_{k=1}^l y_k\Omega]
-[\big((l+1)\frac{q(x_1+1)^2}{x_1 \prod_{k=l+1}^{n-1} y_k}+(y_{l+1}+\dots+y_{n-1})\prod_{k=1}^l y_k\big)\Omega]\nn\\
&=& l\hbar[\prod_{k=1}^l y_k\Omega]
-n[\frac{q(x_1+1)^2}{x_1 \prod_{k=l+1}^{n-1} y_k}\Omega].
\end{eqnarray}
For  $1\leq l<n-1$ we get (\ref{eq-GaussManinConnection-yk}).
When $l=n-1$ in (\ref{eq-proof-prop-PicardFuchs-Quadric-4}), using
\begin{equation}\label{eq-proof-prop-PicardFuchs-Quadric-5}
(d+\frac{d W_0}{\hbar})\big((\prod_{k=1}^{n-1}y_k)x_1\iota_{\partial_{x_1}}\Omega\big)=\hbar^{-1}\frac{q (x_1^2-1)}{x_1}\Omega,
\end{equation}
we get
\begin{eqnarray*}
&&\nabla_{\hbar^2\partial_{\hbar}}^{\mathrm{GM}}[\prod_{k=1}^{n-1}y_k\Omega]= (n-1)\hbar[\prod_{k=1}^{n-1} y_k\Omega]
-n[\frac{q(x_1+1)^2}{x_1 }\Omega]\\
&=& (n-1)\hbar[\prod_{k=1}^{n-1} y_k\Omega]
-n[\frac{q(2x_1^2+2x_1)}{x_1 }\Omega]\\
&=& (n-1)\hbar[\prod_{k=1}^{n-1} y_k\Omega]
-2nq[x_1\Omega]-2nq[\Omega].
\end{eqnarray*}
\end{proof}

\begin{lemma}\label{lem-GaussManinConnection-x1}
Restricting on $U$, we have
\begin{equation}\label{eq-GaussManinConnection-x1}
\nabla_{\hbar^2\partial_{\hbar}}^{\mathrm{GM}}[x_1\Omega]= -n[y_1 \Omega]+n\hbar[x_1\Omega].
\end{equation}
\end{lemma}
\begin{proof}
Since
\begin{eqnarray*}
(d+\frac{d W_0}{\hbar})(x_1 y_j\iota_{\partial_{y_j}}\Omega)=\hbar^{-1}x_1 y_j\frac{\partial W_0}{\partial y_j}\Omega
=\hbar^{-1}(-\frac{q(x_1+1)^2}{\prod_{i=1}^{n-1}y_i}+x_1 y_j)\Omega,
\end{eqnarray*}
we have, for $1\leq j\leq n-1$,
\begin{equation*}
      [x_1 y_j\Omega]= [x_1 y_1\Omega].
\end{equation*}
 By (\ref{eq-GM-connection-hbar}) and (\ref{eq-restriction-W0-Omega-GiventalDomain}) we have
\begin{eqnarray*}
\nabla_{\hbar^2\partial_{\hbar}}^{\mathrm{GM}}[x_1\Omega]=-[x_1 W_0\Omega]
= -[(\frac{q(x_1+1)^2}{\prod_{i=1}^{n-1}y_i}+x_1(y_1+\dots+y_{n-1})\Omega].
\end{eqnarray*}
It follows that
\begin{eqnarray}\label{eq-proof-lem-GaussManinConnection-x1-intermediate}
\nabla_{\hbar^2\partial_{\hbar}}^{\mathrm{GM}}[x_1\Omega]=-[x_1 W_0\Omega]= -n[x_1 y_1 \Omega].
\end{eqnarray}
Note that (I learned this identity from \cite[Proof of Lemma 4.2.3(ii)]{GS15})
\begin{eqnarray*}
x_1 y_1=y_1+x_1(y_1 \frac{\partial W_0}{\partial y_1}+x_1 \frac{\partial W_0}{\partial x_1})
-(y_1 \frac{\partial W_0}{\partial y_1}-x_1 \frac{\partial W_0}{\partial x_1}).
\end{eqnarray*}
Then (\ref{eq-GaussManinConnection-x1}) follows from (\ref{eq-proof-lem-GaussManinConnection-x1-intermediate}) and
\begin{equation*}
(d+\frac{d W_0}{\hbar})(x_1 y_1\iota_{\partial_{y_1}} \Omega)=\hbar^{-1}x_1 y_1\frac{\partial W_0}{\partial y_1}\Omega,
\end{equation*}
\begin{equation*}
(d+\frac{d W_0}{\hbar})(y_1\iota_{\partial_{y_1}} \Omega)=\hbar^{-1} y_1\frac{\partial W_0}{\partial y_1}\Omega,
\end{equation*}
\begin{equation}\label{eq-lem-GaussManinConnection-x1-proof-last}
(d+\frac{d W_0}{\hbar})\big((x_1^2+x_1)\iota_{\partial_{x_1}} \Omega\big)=x_1\Omega+\hbar^{-1}(x_1^2+x_1)\frac{\partial W_0}{\partial x_1}\Omega.
\end{equation}
\end{proof}

\begin{lemma}\label{lem-narrowPicardFuchs-intermediate}
\begin{subequations}\label{eq-narrowPicardFuchs-intermediate}
\begin{align}[left ={(\nabla_{\hbar\partial_{\hbar}}^{\mathrm{GM}})^{l}[\Omega]=\empheqlbrace}]
&(-n)^l \hbar^{-l}[y_1\cdots y_l \Omega],\ \mbox{if}\ 1\leq l\leq n-1,\label{eq-narrowPicardFuchs-intermediate-l<=n-1}\\
&2(-n)^n q \hbar^{-n}[(x_1+1)\Omega],\  \mbox{if}\ l=n.\label{eq-narrowPicardFuchs-intermediate-l=n}
\end{align}
\end{subequations}
\end{lemma}
\begin{proof}
We show this by induction on $l$. The case $l=1$ is the case $l=0$ of (\ref{eq-GaussManinConnection-yk}). Suppose (\ref{eq-narrowPicardFuchs-intermediate-l<=n-1}) holds for $l$, for some $l<n-1$. Then by (\ref{eq-GaussManinConnection-yk}) we have
\begin{eqnarray*}
&& (\nabla_{\hbar\partial_{\hbar}}^{\mathrm{GM}})^{l+1}[\Omega]
=(-n)^l(-l\hbar^{-l}[y_1\cdots y_l \Omega]
+l\hbar^{-l}[y_1\cdots y_l \Omega]-n\hbar^{-l-1}[y_1\cdots y_{l+1} \Omega])\\
&=& (-n)^{l+1}\hbar^{-l-1}[y_1\cdots y_{l+1} \Omega].
\end{eqnarray*}
Finally by (\ref{eq-GaussManinConnection-yk}) we have
\begin{multline*}
(\nabla_{\hbar\partial_{\hbar}}^{\mathrm{GM}})^{n}[\Omega]= (-n)^{n-1}\Big((1-n)\hbar^{-n+1}[y_1\cdots y_{n-1} \Omega]\\
+\hbar^{-n}\big((n-1)\hbar[\prod_{k=1}^l y_k\Omega]-2nq[x_1\Omega]-2nq[\Omega]\big)
\Big)
= 2(-n)^n q\hbar^{-n}([x_1\Omega]+[\Omega]).
\end{multline*}
\end{proof}

\begin{proof}[Proof of Proposition \ref{prop-PicardFuchs-Quadric}]
We return to the domain $V=\{z_1\cdots z_n\neq 1\}$ in our assertion. 
By (\ref{eq-transform-xyToZ}) and (\ref{eq-transform-zToXy}) we have
\begin{eqnarray*}
(x_1^2+x_1)\frac{\partial}{\partial x_1}
= (\prod_{i=1}^n z_i-1)z_1\frac{\partial}{\partial z_1},
\end{eqnarray*}
and
\begin{eqnarray*}
(\prod_{k=1}^{n-1}y_k)x_1\frac{\partial}{\partial x_1}=\prod_{i=2}^{n-1}z_i \cdot (\prod_{i=1}^n z_i-1)\frac{\partial}{\partial z_1},
\end{eqnarray*}
the form $(\prod_{k=1}^{n-1}y_k)x_1\frac{\partial}{\partial x_1}$ on RHS of (\ref{eq-proof-prop-PicardFuchs-Quadric-5}), and $(x_1^2+x_1)\iota_{\partial_{x_1}} \Omega$ on RHS of (\ref{eq-lem-GaussManinConnection-x1-proof-last}), are analytic on $V$. Similarly one finds that  $y_i\iota_{\partial_{y_i}}\Omega$ is analytic on $V$ for $1\leq i\leq n-1$. Then by checking the proofs of Lemma \ref{lem-GaussManinConnection-yk}-\ref{lem-narrowPicardFuchs-intermediate}, one sees that (\ref{eq-GaussManinConnection-yk}), (\ref{eq-GaussManinConnection-x1}), and (\ref{eq-narrowPicardFuchs-intermediate}) hold on $V$. In particular, on $V$, by (\ref{eq-narrowPicardFuchs-intermediate-l=n}), (\ref{eq-GaussManinConnection-yk}), and (\ref{eq-GaussManinConnection-x1}), we have
\begin{eqnarray*}
&& (\nabla_{\hbar\partial_{\hbar}}^{\mathrm{GM}})^{n+1}[\Omega]
= 2(-n)^n q\nabla_{\hbar\partial_{\hbar}}^{\mathrm{GM}}([\hbar^{-n}x_1\Omega]+[\hbar^{-n}\Omega])\\
&=& 2(-n)^{n+1}\hbar^{-n} q(2 \hbar^{-1}[y_1 \Omega]+[\Omega])\\
&=&  2(-n)^{n+1}q\hbar^{-n}(-\frac{2}{n}\nabla_{\hbar\partial_{\hbar}}^{\mathrm{GM}}[\Omega]+[\Omega]).
\end{eqnarray*}
\end{proof}

\begin{proposition}\label{prop-PicardFuchs-Quadric-solution}
The coefficients of $1,\alpha,\dots,\alpha^n$ in
\begin{equation}\label{eq-function-xi}
      \xi(\hbar,\alpha)=q^{\alpha}\hbar^{-n \alpha}\sum_{d=0}^{+\infty}\hbar^{-nd}\frac{(2d)!}{(d!)^{n+2}}\frac{q^d\binom{2 \alpha+2d}{2d}}{\binom{\alpha+d}{d}^2}
\end{equation}
form a system of fundamental solutions to (\ref{eq-PicardFuchs-Quadric}).
\end{proposition}
\begin{proof}
Let $\psi=\hbar^{-n \alpha}\varphi $. Then $\psi$ satisfies (\ref{eq-PicardFuchs-Quadric}) iff
\begin{equation*}
      (\alpha-\frac{1}{n}\hbar\frac{\partial}{\partial \hbar})^{n+1}\varphi
      =2q \hbar^{-n}(2 \alpha+1-\frac{2}{n}\hbar\frac{\partial}{\partial \hbar})\varphi.
\end{equation*}
Write
\begin{equation*}
      \varphi=\sum_{k=0}^{+\infty}\hbar^{-k}g_k.
\end{equation*}
Then
\begin{equation*}
      (\alpha+\frac{n+k}{n})^{n+1}g_{n+k}=2q (2\alpha+1+\frac{2k}{n})g_k.
\end{equation*}
So $g_k=0$ if $n\nmid k$, and if $g_0=1$ then
\begin{equation*}
      g_{nd}=\frac{(4q)^d (\alpha+d-\frac{1}{2})(\alpha+d-\frac{3}{2})\cdots(\alpha+\frac{1}{2})}{\big((\alpha+d)(\alpha+d-1)\cdots(\alpha+1)\big)^{n+1}}=\frac{(4q)^d (\alpha+\frac{1}{2})_d}{\big((\alpha+1)_d\big)^{n+1}}=\frac{(2d)!}{(d!)^{n+2}}
      \frac{q^d\binom{2 \alpha+2d}{2d}}{\binom{\alpha+d}{d}^2}.
\end{equation*}
Thus
\begin{equation*}
      \varphi=\sum_{d=0}^{+\infty}\frac{(4q)^d (\alpha+\frac{1}{2})_d}{\big((\alpha+1)_d\big)^{n+1}}\hbar^{-nd}
      =\sum_{d=0}^{+\infty}\frac{(2d)!}{(d!)^{n+2}}\frac{q^d\binom{2 \alpha+2d}{2d}}{\binom{\alpha+d}{d}^2}\hbar^{-nd}.
\end{equation*}
\end{proof}

\begin{corollary}\label{cor-localHomologyBasis-Xi}
Let $\xi=\sum_{i=0}^n \alpha^i \xi_i$.
There is a local basis of sections  $\Xi_0,\dots,\Xi_n,\Xi_{n+1}$ of $R$  such that the integrals
$\int_{\Xi_i}e^{W_0/\hbar}\Omega=\xi_i$, $i=0,\dots,n$, and $\int_{\Xi_{n+1}}e^{W_0/\hbar}\Omega=0$.
\end{corollary}
\begin{proof}
The $\hbar$-direction is determined by the Picard-Fuchs equation (\ref{eq-PicardFuchs-Quadric}). We are left to determine the variation in the $q$-direction.
Recall
\begin{equation*}
      W_0(z_1,\cdots,z_n,q)=\frac{q z_1^2 z_{n} }{\prod_{i=1}^{n}z_i-1}+z_2+\sum_{i=2}^{n-1}z_{i}z_{i+1}.
\end{equation*}
Substitute $z_i$ by
\begin{equation}\label{eq-changeOfVar-lambda}
      z_i\mapsto \begin{cases}
      \lambda^{-\frac{n}{2}} z_i,& i=1;\\
      \lambda z_i,& i\ \mbox{is even};\\
      z_i,& i\geq 3\ \mbox{and is odd},
      \end{cases}
\end{equation}
one gets
\begin{equation*}
      W_0(\lambda^{-\frac{n}{2}}z_1,\lambda z_2,z_3,\dots,z_{n-1},\lambda z_n,q)=\lambda W_0(z_1,\dots,z_n,\lambda^{-n}q),
\end{equation*}
and 
\begin{equation*}
      \Omega=\frac{dz_1\cdots dz_n}{z_1\cdots z_n-1}\mapsto \Omega.
\end{equation*}
It follows that 
\begin{equation}\label{eq-changeOfVar-narrowPeriod}
      \int_{\Xi_i}e^{\frac{W_0}{\hbar}}\Omega=\int_{\Xi_i}e^{\frac{q^{\frac{1}{n}}W_0(q=1)}{\hbar}}\Omega,
\end{equation}
and thus
\begin{equation*}
\xi(\hbar,\alpha,q)=\sum_{i=0}^n \alpha^i\int_{\Xi_i}e^{\frac{W_0}{\hbar}}\Omega=C q^{\alpha}\hbar^{-n \alpha}\sum_{d=0}^{+\infty}\hbar^{-nd}\frac{(2d)!}{(d!)^{n+2}}\frac{q^d\binom{2 \alpha+2d}{2d}}{\binom{\alpha+d}{d}^2}.
\end{equation*}
where $C$ is a constant.
\end{proof}

\subsection{Irreducibility of the narrow Picard-Fuchs equation}
In this section we show the irreducibility of the Picard-Fuchs equation (\ref{eq-PicardFuchs-Quadric}). By a change of variables $u=4q \hbar^{-n}$ this is equivalent to showing the irreducibility of the ODE
\begin{equation}\label{eq-PicardFuchs-Quadric-inTemrsOfLambda}
      (u\frac{\partial}{\partial u})^{n+1}\psi=u(u\frac{\partial}{\partial u}+\frac{1}{2})\psi.
\end{equation}
Namely, for a nonzero solution $\psi$, the functions $\psi,\psi',\dots,\psi^{(n)}$ are linear independent over $\mathbb{C}[u,u^{-1}]$. Recall that (\ref{eq-PicardFuchs-Quadric-inTemrsOfLambda}) has a solution 
\begin{equation*}
      \psi(u)=\sum_{d\geq 0}\frac{(2d)!}{(d!)^{n+2}} (\frac{u}{4})^d. 
\end{equation*}
For our use in the previous section, we need only to show

\begin{lemma}[Beukers-Brownawell-Heckman]\label{lem-irreducibility}
The series
$\psi,\psi',\dots,\psi^{(n)}$ are linear independent over $\mathbb{C}(u)$.
\end{lemma}
\begin{proof} This is shown in the proof of \cite[Lemma 4.2]{BBH88}. For the reader's convenience, we recall their argument.
 Suppose the existence of $N\in \mathbb{N}$, and $f_i(u)=\sum_{j=0}^N f_{i,j}u^j\in \mathbb{C}[u]$, for $0\leq i\leq n$, which satisfy
\begin{equation}\label{eq-irreducibility-coeff}
      \sum_{i=0}^{n}f_i D^i \psi=0.
\end{equation}
Then
\begin{eqnarray*}
&&\sum_{i=0}^{n}f_i D^i \psi\\
&=& \sum_{i=0}^{n}\sum_{j=0}^{N}\big(f_{i,j}u^j \sum_{d\geq 0}\frac{(2d)!d^j}{(d!)^{n+2}} (\frac{u}{4})^d\big)
\\
&=&\sum_{l\geq 0}u^{l}\left( \sum_{i=0}^n \sum_{j=0}^N \frac{(2(l-j))!(l-j)^i }{4^{l-j}((l-j)!)^{n+2}}f_{i,j}
      \right).
\end{eqnarray*}
So (\ref{eq-irreducibility-coeff}) is equivalent to the collection of  infinitely many linear equations for $f_{k,j}\in \mathbb{C}$
\begin{equation*}
       \sum_{i=0}^n \sum_{j=0}^N \frac{(2(l-j))!(l-j)^i }{4^{l-j}((l-j)!)^{n+2}}f_{i,j}=0.
\end{equation*}
This system of linear equations depends on $n$ and $N$. 
We need to show that its rank  is equal to $(n+1)(N+1)$.

Using the Pochhammer symbol, we rewrite $\psi$ as
\begin{equation*}
      \psi(u)=\sum_{d\geq 0}\frac{(\frac{1}{2})_d}{((1)_d)^{n+1}} u^d,
\end{equation*}
and the linear equations as
\begin{equation*}
       \sum_{i=0}^n \sum_{j=0}^N \frac{(\frac{1}{2})_{l-j}(l-j)^i }{\big((1)_{l-j}\big)^{n+1}}f_{i,j}=0.
\end{equation*}
Define $g_{j}(w)=\sum_{i=0}^{n}f_{i,j}w^{i}$. Then
\begin{equation*}
      \sum_{j=0}^N \frac{(\frac{1}{2})_{l-j}}{\big((1)_{l-j}\big)^{n+1}}g_j(l-j)=0
\end{equation*}
for any $l\geq 0$. Suppose $l\geq N$. Then
\[
\frac{(\frac{1}{2})_{l-j}}{\big((1)_{l-j}\big)^{n+1}}\Big/\frac{(\frac{1}{2})_{l-N}}{\big((1)_{l-N}\big)^{n+1}}
=\frac{\prod_{k=0}^{N-j-1}(\frac{1}{2}+l-N+k)}{\prod_{k=0}^{N-j-1}(1+l-N+k)^{n+1}}.
\]
So for all $l\geq N$,
\[
\sum_{j=0}^N \frac{\prod_{k=0}^{N-j-1}(\frac{1}{2}+l-N+k)}{\prod_{k=0}^{N-j-1}(1+l-N+k)^{n+1}}g_j(l-j)=0.
\]
But this is a rational function for $l$. We get thus an identity for the indeterminate $w$
\[
\sum_{j=0}^N \frac{\prod_{k=0}^{N-j-1}(w+\frac{1}{2}-N+k)}{\prod_{k=0}^{N-j-1}(w+1-N+k)^{n+1}}g_j(w-j)=0.
\]
This is impossible because only when $j=0$, $\frac{\prod_{k=0}^{N-j-1}(w+\frac{1}{2}-N+k)}{\prod_{k=0}^{N-j-1}(w+1-N+k)^{n+1}}$ has a pole  at $w=0$, and moreover it is a pole of order $n+1$, while $\mathrm{deg}(g_j)\leq n$.
\end{proof}

\begin{corollary}\label{cor-irreducibility-narrowPicardFuchsOperator}
The Picard-Fuchs operator
\begin{equation}\label{eq-narrowPicardFuchsOperator}
      \mathcal{P}_{\mathrm{na}}:=(-\frac{\hbar}{n}\frac{\partial}{\partial \hbar})^{n+1}-2q\hbar^{-n}(-\frac{2\hbar}{n}\frac{\partial}{\partial \hbar}+1)
\end{equation}
is irreducible.
\end{corollary}

\subsection{Picard-Fuchs equation for broad periods}
In this section, $n$ is an even integer and $n\geq 4$.
\begin{proposition}\label{prop-GaussManinConnection-broad-0}
\begin{eqnarray}\label{eq-GaussManinConnection-broad-0}
\nabla_{\hbar^2\partial_{\hbar}}^{\mathrm{GM}}[z_1\Omega]=\frac{n}{2}\hbar[z_1\Omega]-n[z_1 z_2\Omega].
\end{eqnarray}
\end{proposition}
\begin{proof}
Using (\ref{eq-derivatives-W0-zj}) we have
\begin{eqnarray}\label{eq-prop-GaussManinConnection-broad-0-proof-1}
&&(d+\frac{d W_0}{\hbar})\big(z_1^2 \iota_{\partial_{z_1}}\Omega
-z_1 z_n \iota_{\partial_{z_n}}\Omega\big)\nn\\
&=& \big(z_1+\hbar^{-1}\big(\frac{q z_1^3 z_n}{\prod_{i=1}^n z_i-1}-z_1 z_{n-1} z_n)\big)\Omega,
\end{eqnarray}
\begin{eqnarray}\label{eq-prop-GaussManinConnection-broad-0-proof-2}
&&(d+\frac{d W_0}{\hbar})\big(z_1z_2 \iota_{\partial_{z_2}}\Omega
-z_1 z_n \iota_{\partial_{z_n}}\Omega\big)\nn\\
&=& \hbar^{-1}\big(-\frac{q z_1^3 z_n}{\prod_{i=1}^n z_i-1}+z_1z_2+z_1z_2z_3-z_1 z_{n-1} z_n)\big)\Omega,
\end{eqnarray}
and for $3\leq j\leq n-1$,
\begin{eqnarray}\label{eq-prop-GaussManinConnection-broad-0-proof-3}
&&(d+\frac{d W_0}{\hbar})\big(z_1z_j \iota_{\partial_{z_j}}\Omega
-z_1 z_n \iota_{\partial_{z_n}}\Omega\big)\nn\\
&=& \hbar^{-1}\big(-\frac{q z_1^3 z_n}{\prod_{i=1}^n z_i-1}+z_1z_{j-1}z_j+z_1z_jz_{j+1}-z_1 z_{n-1} z_n)\big)\Omega.
\end{eqnarray}
Thus (\ref{eq-prop-GaussManinConnection-broad-0-proof-2}), (\ref{eq-prop-GaussManinConnection-broad-0-proof-3}) and induction on $j$ yield, for $3\leq j\leq n$
\begin{equation*}
   [z_1z_{j-1}z_{j}\Omega]= 
   \begin{cases}
   [z_1 z_2\Omega],& \mbox{if $j$ is even},\\
      [z_1 z_2 z_3\Omega],& \mbox{if $j$ is odd}.
   \end{cases}  
\end{equation*}
Then (\ref{eq-prop-GaussManinConnection-broad-0-proof-1}) and (\ref{eq-prop-GaussManinConnection-broad-0-proof-2}) yield
\[
[\frac{q z_1^3 z_n}{\prod_{i=1}^n z_i-1}\Omega]=
[z_1z_{2}z_{3}\Omega]=[z_1 z_2\Omega]-\hbar[z_1\Omega].
\]
It follows that
\begin{eqnarray*}
&&\nabla_{\hbar^2\partial_{\hbar}}^{\mathrm{GM}}[z_1\Omega]=
 -[(\frac{q z_1^3 z_{n} }{\prod_{i=1}^{n}z_i-1}+z_1z_2+z_1\sum_{i=2}^{n-1}z_{i}z_{i+1}) \Omega]\\
 &=&\frac{n}{2}\hbar[z_1\Omega]-n[z_1 z_2\Omega].
\end{eqnarray*}
\end{proof}

\begin{proposition}\label{prop-GaussManinConnection-broad-k}
For $1\leq k\leq \frac{n}{2}-1$,
\begin{eqnarray}\label{eq-GaussManinConnection-broad-k}
\nabla_{\hbar^2\partial_{\hbar}}^{\mathrm{GM}}[z_1\cdots z_{2k}\Omega]
=(\frac{n}{2}+k)\hbar[z_1\cdots z_{2k}\Omega]-n[z_1\cdots z_{2k+2}\Omega].
\end{eqnarray}
\end{proposition}
\begin{proof}
Using (\ref{eq-derivatives-W0-zj}) we have
\begin{eqnarray*}
&&(d+\frac{d W_0}{\hbar})\big(z_1^2 z_2\cdots z_{2k} \iota_{\partial_{z_1}}\Omega
-z_1\cdots z_{2k}\cdot z_n \iota_{\partial_{z_n}}\Omega\big)\\
&=& \big(z_1\cdots z_{2k}+\hbar^{-1}(\frac{q z_1^3z_2\cdots z_{2k} z_n}{\prod_{i=1}^n z_i-1}-z_1\cdots z_{2k}\cdot z_{n-1} z_n)\big)\Omega,
\end{eqnarray*}
\begin{eqnarray*}
&&(d+\frac{d W_0}{\hbar})\big(z_1z_2^2 z_3\cdots z_{2k} \iota_{\partial_{z_2}}\Omega
-z_1\cdots z_{2k}\cdot z_n \iota_{\partial_{z_n}}\Omega\big)\\
&=& \Big(z_1\cdots z_{2k}+\hbar^{-1}\big(-\frac{q z_1^3z_2\cdots z_{2k} z_n}{\prod_{i=1}^n z_i-1}
+z_1\cdots z_{2k}(z_2+z_2 z_3-z_{n-1} z_n)\big)\Big)\Omega,
\end{eqnarray*}
and for $3\leq j\leq 2k$,
\begin{eqnarray*}
&&(d+\frac{d W_0}{\hbar})\big(z_1\cdots z_{2k}\cdot z_j \iota_{\partial_{z_j}}\Omega
-z_1\cdots z_{2k}\cdot z_n \iota_{\partial_{z_n}}\Omega\big)\\
&=&\Big(z_1\cdots z_{2k}+\hbar^{-1}\big(-\frac{q z_1^3 z_2\cdots z_{2k} z_n}{\prod_{i=1}^n z_i-1}
+z_1\cdots z_{2k}(z_{j-1}z_j+z_jz_{j+1}-z_{n-1} z_n)\big)\Big)\Omega,
\end{eqnarray*}
and for $2k+1\leq j\leq n-1$,
\begin{eqnarray*}
&&(d+\frac{d W_0}{\hbar})\big(z_1\cdots z_{2k}\cdot z_j \iota_{\partial_{z_j}}\Omega
-z_1\cdots z_{2k}\cdot z_n \iota_{\partial_{z_n}}\Omega\big)\\
&=&\hbar^{-1}\big(-\frac{q z_1^3 z_2\cdots z_{2k} z_n}{\prod_{i=1}^n z_i-1}
+z_1\cdots z_{2k}(z_{j-1}z_j+z_jz_{j+1}-z_{n-1} z_n)\big)\Omega.
\end{eqnarray*}

It follows that
\begin{subequations}
\begin{eqnarray}
&& [\frac{q z_1^3 z_2\cdots z_{2k} z_n}{\prod_{i=1}^n z_i-1}\Omega]\nn\\
&=& [z_1\cdots z_{2k}\cdot z_{n-1}z_n\Omega]-\hbar[z_1\cdots z_{2k}\Omega]  \label{eq-prop-GaussManinConnection-broad-k-proor-1}\\
&=& [z_1\cdots z_{2k}\cdot z_2\Omega]+[z_1\cdots z_{2k}\cdot z_2z_3\Omega]-[z_1\cdots z_{2k}\cdot z_{n-1}z_n\Omega]
+\hbar[z_1\cdots z_{2k}\Omega]   \label{eq-prop-GaussManinConnection-broad-k-proor-2}\\
&=& \left\{\begin{gathered}
[z_1\cdots z_{2k}\cdot z_{j-1}z_j\Omega]+[z_1\cdots z_{2k}\cdot z_jz_{j+1}\Omega]\\
-[z_1\cdots z_{2k}\cdot z_{n-1}z_n\Omega]
+\hbar[z_1\cdots z_{2k}\Omega],\ \mbox{if}\ 3\leq j\leq 2k;   \label{eq-prop-GaussManinConnection-broad-k-proor-3}\\
[z_1\cdots z_{2k}\cdot z_{j-1}z_j\Omega]+[z_1\cdots z_{2k}\cdot z_jz_{j+1}\Omega]\\
-[z_1\cdots z_{2k}\cdot z_{n-1}z_n\Omega],\ \mbox{if}\ 2k+1\leq j\leq n-1.\label{eq-prop-GaussManinConnection-broad-k-proor-4}
\end{gathered}\right.  
\end{eqnarray}
\end{subequations}
For $3\leq j\leq 2k$, by (\ref{eq-prop-GaussManinConnection-broad-k-proor-2})=(\ref{eq-prop-GaussManinConnection-broad-k-proor-3}) and induction on $j$ we get
\begin{equation}\label{eq-prop-GaussManinConnection-broad-k-proor-5}
      [z_1\cdots z_{2k}\cdot z_{2i-1}z_{2i}\Omega]=[z_1\cdots z_{2k}\cdot z_2 \Omega],\ \mbox{for}\ 2\leq i\leq k,
\end{equation}
and
\begin{equation}\label{eq-prop-GaussManinConnection-broad-k-proor-6}
      [z_1\cdots z_{2k}\cdot z_{2i}z_{2i+1}\Omega]=[z_1\cdots z_{2k}\cdot z_{2}z_{3} \Omega],\ \mbox{for}\ 2\leq i\leq k.
\end{equation}
Then for $2k+1\leq j\leq n-1$, by  (\ref{eq-prop-GaussManinConnection-broad-k-proor-2})=(\ref{eq-prop-GaussManinConnection-broad-k-proor-3}) and induction on $j$ we get
\begin{equation}\label{eq-prop-GaussManinConnection-broad-k-proor-7}
      [z_1\cdots z_{2k}\cdot z_{2i-1}z_{2i}\Omega]=[z_1\cdots z_{2k}\cdot z_2 \Omega]+\hbar[z_1\cdots z_{2k}\Omega],\ \mbox{for}\ k+1\leq i\leq \frac{n}{2},
\end{equation}
and
\begin{equation}\label{eq-prop-GaussManinConnection-broad-k-proor-8}
      [z_1\cdots z_{2k}\cdot z_{2i}z_{2i+1}\Omega]=[z_1\cdots z_{2k}\cdot z_{2}z_{3} \Omega],\ \mbox{for}\ k+1\leq i\leq \frac{n}{2}-1.
\end{equation}
The assumption $k\leq \frac{n}{2}-1$ allows us to take $i=\frac{n}{2}$ in (\ref{eq-prop-GaussManinConnection-broad-k-proor-7}). Then (\ref{eq-prop-GaussManinConnection-broad-k-proor-1})=(\ref{eq-prop-GaussManinConnection-broad-k-proor-2}) yields
\begin{equation*}
      [z_1\cdots z_{2k}\cdot z_2\Omega]=[z_1\cdots z_{2k}\cdot z_2z_3\Omega].
\end{equation*}
Thus
\begin{equation}\label{eq-prop-GaussManinConnection-broad-k-proor-9}
    [z_1\cdots z_{2k}\cdot z_{j-1}z_{j}\Omega]=
    \begin{cases}
    [z_1\cdots z_{2k}\cdot z_2 \Omega]+\hbar[z_1\cdots z_{2k}\Omega],& \mbox{if $j$ is even and}\ j\geq 2k+2,\\
     [z_1\cdots z_{2k}\cdot z_2 \Omega],& \mbox{for other $j\geq 3$},
    \end{cases}
\end{equation}
and (\ref{eq-prop-GaussManinConnection-broad-k-proor-1}) yields
\begin{eqnarray}\label{eq-prop-GaussManinConnection-broad-k-proor-10}
[\frac{q z_1^3 z_2\cdots z_{2k}\cdot z_n}{\prod_{i=1}^n z_i-1}\Omega]=[z_1\cdots z_{2k+2} \Omega]-\hbar[z_1\cdots z_{2k}\Omega].
\end{eqnarray}
Hence by (\ref{eq-GM-connection-hbar}), (\ref{eq-prop-GaussManinConnection-broad-k-proor-9}), and (\ref{eq-prop-GaussManinConnection-broad-k-proor-10})
 we have
\begin{eqnarray*}
&&\nabla_{\hbar^2\partial_{\hbar}}^{\mathrm{GM}}[z_1\cdots z_{2k}\Omega]
= -[(\frac{q z_1^3 z_2\cdots z_{2k}z_{n} }{\prod_{i=1}^{n}z_i-1}+z_1\cdots z_{2k}(z_2+\sum_{i=2}^{n-1}z_{i}z_{i+1}) \Omega]\\
&=&(\frac{n}{2}+k)\hbar[z_1\cdots z_{2k}\Omega]-n[z_1\cdots z_{2k+2}\Omega].
\end{eqnarray*}
\end{proof}

\begin{corollary}\label{cor-PicardFuchs-broadPeriod}
For $1\leq k\leq \frac{n}{2}$,
\begin{eqnarray}\label{eq-PicardFuchs-broadPeriod}
(\nabla_{\hbar\partial_{\hbar}}^{\mathrm{GM}}-\frac{n}{2})^{k}[z_1\Omega]
=(-n)^k\hbar^{-k}[\prod_{i=1}^{2k}z_i\Omega].
\end{eqnarray}
\end{corollary}
\begin{proof}
We prove by induction on $k$. The case $k=1$ is (\ref{eq-GaussManinConnection-broad-0}). 
Suppose that our assertion is true for $k$. 
Then (\ref{eq-GaussManinConnection-broad-k}) yields
\begin{eqnarray*}
&&(\nabla_{\hbar\partial_{\hbar}}^{\mathrm{GM}}-\frac{n}{2})^{k+1}[z_1\Omega]
=(-n)^k(\nabla_{\hbar\partial_{\hbar}}^{\mathrm{GM}}-\frac{n}{2})\hbar^{-k}[\prod_{i=1}^{2k}z_i\Omega]\\
&=&(-n)^k\big(-k\hbar^{-k}[\prod_{i=1}^{2k}z_i\Omega]
+(\frac{n}{2}+k)\hbar^{-k}[\prod_{i=1}^{2k}z_i\Omega]-n\hbar^{-k-1}[\prod_{i=1}^{2k+2}z_i\Omega]
-\frac{n}{2}\hbar^{-k}[\prod_{i=1}^{2k}z_i\Omega]\big)\\
&=&(-n)^{k+1}\hbar^{-k-1}[\prod_{i=1}^{2k+2}z_i\Omega].
\end{eqnarray*}
\end{proof}

\begin{corollary}
\begin{eqnarray}\label{eq-GaussManinConnection-broadToNarrow}
(\nabla_{\hbar\partial_{\hbar}}^{\mathrm{GM}}-\frac{n}{2})^{\frac{n}{2}}[z_1\Omega]
=\frac{\hbar^{\frac{n}{2}}}{2(-n)^{\frac{n}{2}}q}(\nabla_{\hbar\partial_{\hbar}}^{\mathrm{GM}})^{n}[\Omega].
\end{eqnarray}
\end{corollary}
\begin{proof}
By the proof of Proposition \ref{prop-PicardFuchs-Quadric}, (\ref{eq-narrowPicardFuchs-intermediate-l=n}) holds on $V$. Namely we have
\begin{eqnarray*}
(-n)^{\frac{n}{2}}\hbar^{-\frac{n}{2}}[\prod_{i=1}^{n}z_i\Omega]
=\frac{\hbar^{\frac{n}{2}}}{2(-n)^{\frac{n}{2}}q}(\nabla_{\hbar\partial_{\hbar}}^{\mathrm{GM}})^{n}[\Omega].
\end{eqnarray*}
Then (\ref{eq-GaussManinConnection-broadToNarrow}) follows from the $k=\frac{n}{2}$ case of (\ref{eq-PicardFuchs-broadPeriod}).
\end{proof}

\begin{proposition}\label{prop-fix-homologyBasis-preliminary}
Suppose $q=1$. 
Then locally there exists a unique basis  $\Xi_0,\dots,\Xi_{n+1}$ of the local system $R$ such that
\begin{equation*}
      \int_{\Xi_i}e^{\frac{W_0}{\hbar}} \Omega=\xi_i,\  0\leq i\leq n,
\end{equation*}
\begin{equation*}
      \int_{\Xi_{n+1}}e^{\frac{W_0}{\hbar}} \Omega=0,
\end{equation*}
\begin{equation}\label{eq-broadPeriod-Xin+1}
      \int_{\Xi_{n+1}}e^{\frac{W_0}{\hbar}}z_1 \Omega=
      \hbar^{\frac{n}{2}},
\end{equation}
and for $0\leq k\leq n$,
\begin{equation}\label{eq-thetak-xik-relation-monodromy}
      \int_{\Xi_k}e^{\frac{W_0}{\hbar}}z_1\Omega=\hbar^{\frac{n}{2}}\sum_{i=1}^{\frac{n}{2}-1}c_i(\Xi_k) \log^i\hbar
      +\frac{\hbar^{\frac{n}{2}}}{2(-n)^{\frac{n}{2}}}(\hbar\frac{\partial}{\partial \hbar})^{\frac{n}{2}}\int_{\Xi_k}e^{\frac{W_0}{\hbar}} \Omega,
\end{equation}
where$c_i(\Xi_k)\in \mathbb{C}$.
\end{proposition}
\begin{proof}
A  solution to $(\nabla_{\hbar\partial_{\hbar}}^{\mathrm{GM}}-\frac{n}{2})^{\frac{n}{2}}f=0$ has the form
\begin{equation*}
      f(\hbar)=\hbar^{\frac{n}{2}}\sum_{i=0}^{\frac{n}{2}-1}a_i \log^i\hbar,
\end{equation*}
where $a_i\in \mathbb{C}$. Thus (\ref{eq-GaussManinConnection-broadToNarrow}) yields
\begin{equation}\label{eq-GaussManinConnection-broadToNarrow-universalSolution}
      \int_{\Xi}e^{\frac{W_0}{\hbar}}z_1\Omega=\hbar^{\frac{n}{2}}\sum_{i=0}^{\frac{n}{2}-1}c_i(\Xi) \log^i\hbar
      +\frac{\hbar^{\frac{n}{2}}}{2(-n)^{\frac{n}{2}}}(\hbar\frac{\partial}{\partial \hbar})^{\frac{n}{2}}\int_{\Xi}e^{\frac{W_0}{\hbar}} \Omega,
\end{equation}
where $c_i(\Xi)\in \mathbb{C}$ depends on $\Xi$.
Suppose that $\Xi_0,\dots,\Xi_n,\Xi_{n+1}$ is a local basis of $R$ in Corollary \ref{cor-localHomologyBasis-Xi} so that
\begin{equation*}
      \xi_i=\int_{\Xi_i}e^{\frac{W_0}{\hbar}} \Omega
\end{equation*}
and $\int_{\Xi_{n+1}}e^{W_0/\hbar}\Omega=0$. 
 Let $\hbar$ undergo a counterclockwise loop in the $\mathbb{C}^{\times}_{\hbar}$ plane around $0$. Then $\hbar^{-n\alpha}$ is replaced by $\hbar^{-n\alpha}\exp(-n\times 2\pi \sqrt{-1}\alpha)$. Then from (\ref{eq-function-xi}) it follows that there exist $b_0,\dots,b_{n+1}\in \mathbb{C}$ such that
\begin{equation}\label{eq-monodromy-Xik}
      \Xi_k\mapsto \sum_{j=0}^k\frac{(-2n\pi\sqrt{-1})^{k-j}}{(k-j)!}\Xi_j+ b_k \Xi_{n+1},\ \mbox{for}\ 0\leq k\leq n,
\end{equation}
and
\begin{equation*}
      \Xi_{n+1}\mapsto b_{n+1} \Xi_{n+1},\ \mbox{for}\ 0\leq k\leq n.
\end{equation*}
Compare whats happen on both sides of (\ref{eq-GaussManinConnection-broadToNarrow-universalSolution}), in the case $\Xi=\Xi_{n+1}$:
\begin{eqnarray*}
\int_{\Xi_{n+1}}e^{\frac{W_0}{\hbar}}z_1\Omega &\mapsto& b_{n+1}\int_{\Xi_{n+1}}e^{\frac{W_0}{\hbar}}z_1\Omega,\\
\hbar^{\frac{n}{2}}\sum_{i=0}^{\frac{n}{2}-1}c_i(\Xi_{n+1}) \log^i\hbar &\mapsto&
\hbar^{\frac{n}{2}}\sum_{i=0}^{\frac{n}{2}-1}c_i(\Xi_{n+1}) (\log\hbar+2\pi\sqrt{-1})^i.
\end{eqnarray*}
It follows that $c_i(\Xi_{n+1})=0$ for $i>0$, and $b_{n+1}=1$. Therefore we can fix $\Xi_{n+1}$ by demanding $c_0(\Xi_{n+1})=1$, and we do. Moreover, we fix $\Xi_0,\dots,\Xi_{n+1}$ by demanding that
\begin{equation}\label{eq-fix-Xi}
      c_0(\Xi_i)=\delta_{i,n+1}.
\end{equation}
\end{proof}

\begin{lemma}\label{lem-vanishing-ciXik-preliminary}
If $i\geq k-\frac{n}{2}$, $c_i(\Xi_k)=0$.
\end{lemma}
\begin{proof}
Similar to the above proof, let $\hbar$ undergo a counterclockwise loop in the $\mathbb{C}^{\times}_{\hbar}$ plane around $0$.  We 
compare what happens on both sides of (\ref{eq-GaussManinConnection-broadToNarrow-universalSolution}), in the case $\Xi=\Xi_{k}$, for $0\leq k\leq n+1$. We have
\begin{eqnarray*}
&&\hbar^{\frac{n}{2}}\sum_{i=1}^{\frac{n}{2}-1}c_i(\Xi_k) \log^i\hbar
      +\frac{\hbar^{\frac{n}{2}}}{2(-n)^{\frac{n}{2}}}(\hbar\frac{\partial}{\partial \hbar})^{\frac{n}{2}}\xi_k\\
&\mapsto&\hbar^{\frac{n}{2}}\sum_{i=1}^{\frac{n}{2}-1}c_i(\Xi_k) (\log\hbar+2\pi\sqrt{-1})^i
+\frac{\hbar^{\frac{n}{2}}}{2(-n)^{\frac{n}{2}}}
\sum_{j=0}^k\frac{(-2n\pi\sqrt{-1})^{k-j}}{(k-j)!}\xi_j,
\end{eqnarray*}
while, by (\ref{eq-monodromy-Xik}) and (\ref{eq-fix-Xi}),
\begin{eqnarray*}
\int_{\Xi_k}e^{\frac{W_0}{\hbar}}z_1\Omega &\mapsto&
\int_{\sum_{j=0}^k\frac{(-2n\pi\sqrt{-1})^{k-j}}{(k-j)!}\Xi_j+ b_k \Xi_{n+1}}e^{\frac{W_0}{\hbar}}z_1\Omega\\
&=&\sum_{j=0}^k\Big(\frac{(-2n\pi\sqrt{-1})^{k-j}}{(k-j)!}\hbar^{\frac{n}{2}}\sum_{i=1}^{\frac{n}{2}-1}c_i(\Xi_j) \log^i\hbar\Big)\\
&&+\frac{\hbar^{\frac{n}{2}}}{2(-n)^{\frac{n}{2}}}
\sum_{j=0}^k\frac{(-2n\pi\sqrt{-1})^{k-j}}{(k-j)!}\xi_j+b_k\hbar^{\frac{n}{2}}.
\end{eqnarray*}
It follows that
\begin{eqnarray}\label{lem-ciXik-consequenceOfMonodromy}
&&\sum_{i=1}^{\frac{n}{2}-1}c_i(\Xi_k) (\log\hbar+2\pi\sqrt{-1})^i\nn\\
&=&\sum_{j=0}^k\Big(\frac{(-2n\pi\sqrt{-1})^{k-j}}{(k-j)!}\sum_{i=1}^{\frac{n}{2}-1}c_i(\Xi_j) \log^i\hbar\Big)
+b_k.
\end{eqnarray}
For $0\leq k\leq n$, comparing the coefficients of $\log^{\frac{n}{2}-1}\hbar$ of both sides of (\ref{lem-ciXik-consequenceOfMonodromy}) yields linear equations for $c_{\frac{n}{2}-1}(\Xi_k)$. The matrix for this system of linear equations is an upper triangular matrix with a zero diagonal and all superdiagonal entries are nonzero. It follows that $c_{\frac{n}{2}-1}(\Xi_k)=0$ for $k\leq n-1$. This shows the $i=\frac{n}{2}-1$ case of this lemma.

Suppose we have shown our assertion for $l+1\leq i\leq \frac{n}{2}-1$. For $0\leq k\leq l+\frac{n}{2}+1$, comparing the coefficients of $\log^{l}\hbar$ of both sides of (\ref{lem-ciXik-consequenceOfMonodromy}) yields linear equations for $c_{l}(\Xi_k)$. Then the same argument as above yields $c_{l}(\Xi_k)=0$ for $k\leq l+\frac{n}{2}$. Hence a downward induction on $i$ concludes the proof.

\end{proof}

\begin{corollary}\label{cor-APrioriPicardFuchs-broad}
Let 
\begin{equation*}
      \theta_i=\int_{\Xi_i}e^{\frac{W_0}{\hbar}}z_1 \Omega.
\end{equation*}
Denote $D=\hbar\frac{\partial}{\partial \hbar}$. Then locally there exist $v_k(\hbar)\in \mathbb{C}[\hbar,\hbar^{-1}]$, $0\leq k \leq n+2$, such that 
\begin{equation}\label{eq-APrioriPicardFuchs-broad}
      \sum_{k=0}^{n+2}v_k(\hbar)D^k \theta_i=0,\ \mbox{for}\ 0\leq i\leq n+1.
\end{equation}
They satisfy
\begin{equation}\label{eq-equationsOfCoeff-PicardFuchs-broad-0}
      \sum_{k=0}^{n+2} (\frac{n}{2})^k v_k(\hbar)=0.
\end{equation} 
Moreover, let $(\mathcal{P}_{\mathrm{na}})$ be the left ideal generated by $\mathcal{P}_{\mathrm{na}}$ in the ring of differential operators  $\mathbb{C}(\hbar)[\partial_{\hbar}]$. Then 
\begin{equation}\label{eq-equationsOfCoeff-PicardFuchs-broad-1}
      \sum_{k=0}^{n+2} v_k(\hbar)(\frac{n}{2}+D)^k D^{\frac{n}{2}}\in (\mathcal{P}_{\mathrm{na}}).
\end{equation}
\end{corollary}
\begin{proof}
The existence of $v_k(\hbar)$ in our assertion follows from Lemma \ref{lem-algebraicity-restrictionOfcalE}.
Let $D^k$ act on (\ref{eq-thetak-xik-relation-monodromy}) (resp. (\ref{eq-broadPeriod-Xin+1})) from the left, for $1\leq i\leq \frac{n}{2}+1$ using Lemma \ref{lem-vanishing-ciXik-preliminary} we get
\begin{equation}\label{eq-thetak-xik-relation-monodromy-2}
D^k\theta_i=\frac{\hbar^{\frac{n}{2}}}{2(-n)^{\frac{n}{2}}}(\frac{n}{2}+D)^k D^{\frac{n}{2}} \xi_i.
\end{equation}
\[
\big(\mbox{resp.}\quad
D^k\theta_{n+1}=(\frac{n}{2})^k \hbar^{\frac{n}{2}}.\big)
\]
The latter equation yields (\ref{eq-equationsOfCoeff-PicardFuchs-broad-0}). 
Recall Proposition \ref{prop-PicardFuchs-Quadric}, for $0\leq i\leq n$, $\mathcal{P}_{\mathrm{na}}\xi_i=0$.
By Corollary \ref{cor-irreducibility-narrowPicardFuchsOperator}, $\mathcal{P}_{\mathrm{na}}$ is irreducible. So (\ref{eq-equationsOfCoeff-PicardFuchs-broad-1}) follows from  (\ref{eq-thetak-xik-relation-monodromy-2}).
\end{proof}

\begin{lemma}\label{lem-PicardFuchs-EqCoeff-rank}
The functions $(v_i)_{0\leq i\leq n+2}$ satisfying (\ref{eq-equationsOfCoeff-PicardFuchs-broad-0}), 
(\ref{eq-equationsOfCoeff-PicardFuchs-broad-1}) is unique up to a common factor.
\end{lemma}
\begin{proof}
By induction on $k$, we have
\begin{eqnarray}\label{eq-Dn+k-mod-narrowPicardFuchsOperator}
&& D^{n+k}\equiv 2 n^n\hbar^{-n}(-n+D)^{k-1}(2D-n),\mod (\mathcal{P}_{\mathrm{na}}).
\end{eqnarray}
For $0\leq k\leq n+2$, we set
\begin{equation*}
      E_k:=(\frac{n}{2}+D)^k D^{\frac{n}{2}}.
\end{equation*}
For  $\frac{n}{2}+1\leq k\leq n+2$, expanding $E_k$ and using (\ref{eq-Dn+k-mod-narrowPicardFuchsOperator}) we get
\begin{eqnarray}\label{eq-Ek-modulo-narrowPicardFuchsOperator}
E_k&\equiv& \sum_{i=0}^{\frac{n}{2}}\binom{k}{i}(\frac{n}{2})^{k-i}D^{\frac{n}{2}+i}\nn\\
&&+2 n^n\hbar^{-n}(2D-n)\sum_{i=\frac{n}{2}+1}^{k}\binom{k}{i}(\frac{n}{2})^{k-i}
(D-n)^{i-\frac{n}{2}-1}, \mod (\mathcal{P}_{\mathrm{na}}).
\end{eqnarray}
For $0\leq k\leq n+2$, let $G_k=$ RHS of (\ref{eq-Ek-modulo-narrowPicardFuchsOperator}). Then
\begin{equation}\label{eq-sum-fkGk}
      \sum_{k=0}^{n+2}v_k G_k=0.
\end{equation}
For each $0\leq j\leq n$, the coefficient of $D^j$ in the sum (\ref{eq-sum-fkGk}) gives a linear equation for $v_k$'s. Together with (\ref{eq-equationsOfCoeff-PicardFuchs-broad-0}), we have a system of $n+2$ linear equations for $v_k$'s. For our assertion we need to show that this system has a full rank $n+2$.

More precisely, for $0\leq k\leq n+2$ we define $a_{k,j}\in \mathbb{C}[q,\hbar^{-1}]$ by
\begin{equation*}
      G_k=\sum_{j=0}^{n}a_{k,j}D^j,
\end{equation*}
and let
\begin{equation*}
      a_{k,n+1}=(\frac{n}{2})^k.
\end{equation*}
We need to show that the rank of the matrix $M=(a_{k,j})_{0\leq k\leq n+2,0\leq j\leq n+1}$ is $n+2$, where we regard $k$ (resp. $j$) as the index of columns (resp. rows); in particular, both indices start from 0. For example, when $n=6$ the matrix $M$ is
{\scriptsize
\begin{equation*}
\begin{pNiceArray}{cccc:cccc:c}
 0 & 0 & 0 & 0 & -559872 \hbar^{-6} & -5038848 \hbar^{-6} & -35271936 \hbar^{-6} & -196515072 \hbar^{-6} & -997691904 \hbar^{-6} \\
 0 & 0 & 0 & 0 & 186624 \hbar^{-6} & 1119744 \hbar^{-6} & 8398080 \hbar^{-6} & 40310784 \hbar^{-6} & 211631616 \hbar^{-6} \\
 0 & 0 & 0 & 0 & 0 & 186624 \hbar^{-6} & 559872 \hbar^{-6} & 6718464 \hbar^{-6} & 20155392 \hbar^{-6} \\
 1 & 3 & 9 & 27 & 81 & 243 & 186624 \hbar^{-6}+729 & 2187 & 6718464 \hbar^{-6}+6561 \\
 \hdashline  
 0 & 1 & 6 & 27 & 108 & 405 & 1458 & 186624 \hbar^{-6}+5103 & 17496-559872 \hbar^{-6} \\
 0 & 0 & 1 & 9 & 54 & 270 & 1215 & 5103 & 186624 \hbar^{-6}+20412 \\
 0 & 0 & 0 & 1 & 12 & 90 & 540 & 2835 & 13608 \\
 1 & 3 & 9 & 27 & 81 & 243 & 729 & 2187 & 6561 \\
\end{pNiceArray}.
\end{equation*}}

One sees that
\begin{enumerate}
      \item[1)] The $(\frac{n}{2}+1)\times (\frac{n}{2}+1)$ submatrix of $M$ at the bottom-left corner is nonsingular. 
      \item[2)] For $0\leq k\leq n+2$, the entry in the $\frac{n}{2}$-th row
      \begin{equation*}
            a_{k,\frac{n}{2}}=(\frac{n}{2})^k+2 n^n\hbar^{-n}\mathrm{Coeff}_{D^j}\big\{(2D-n)\sum_{i=\frac{n}{2}+1}^{k}\binom{k}{i}(\frac{n}{2})^{k-i}
(D-n)^{i-\frac{n}{2}-1}\big\}
      \end{equation*}
      has the same constant term as the entry $a_{k,n+1}$ in the bottom row. 
\end{enumerate}
Denote by $\tilde{M}$ the matrix obtained from $M$ by subtracting the bottom row from the $\frac{n}{2}$-th row. Then one sees furthermore that
\begin{enumerate}
      \item[3)] The $(\frac{n}{2}+1)\times (\frac{n}{2}+1)$ submatrix of $\tilde{M}$ at the top-left  corner is 0.
\end{enumerate}
Suppose  $\tilde{M}=(\tilde{a}_{k,j})_{0\leq k\leq n+2,0\leq j\leq n+1}$. Then $\tilde{a}_{k,j}=2n^n  \hbar^{-n} e_{k,j}$, where
\begin{equation}\label{eq-def-ekj}
      e_{k,j}=\mathrm{Coeff}_{D^j}\big\{(2D-n)\sum_{i=\frac{n}{2}+1}^{k}\binom{k}{i}(\frac{n}{2})^{k-i}
(D-n)^{i-\frac{n}{2}-1}\big\}.
\end{equation}
From the above observations, to show $\mathrm{rank}\ M=n+2$, or equivalently $\mathrm{rank}\ \tilde{M}=n+2$, it suffices to show that the $(\frac{n}{2}+1)\times (\frac{n}{2}+1)$ submatrix 
$N:=(\tilde{a}_{k,j})_{\frac{n}{2}+1\leq k\leq n+1,0\leq j\leq \frac{n}{2}}$ (as indicated above in the case $n=6$) is nonsingular. By (\ref{eq-def-ekj}), we have
\begin{equation*}
      e_{k,j}=0\ \mbox{for}\ k-j<\frac{n}{2}
\end{equation*}
and
\begin{equation*}
      e_{\frac{n}{2}+j,j}=2.
\end{equation*}
This means that 
\begin{equation}\label{eq-lower-left-ofN-unipotent}
      \mbox{the lower-left $\frac{n}{2}\times \frac{n}{2}$ submatrix of $N$ is unipotent up to a scalar multiple.}
\end{equation}
If $\frac{n}{2}+1\leq k\leq n$, the highest power of $D$ in (\ref{eq-def-ekj}) is lower than or equal to $\frac{n}{2}$. So
 \begin{equation}\label{eq-sum-ekj-k<=n}
       \sum_{j=0}^{\frac{n}{2}}(\frac{n}{2})^j e_{k,j}= \sum_{j=0}^{\frac{n}{2}}(\frac{n}{2})^j \mathrm{Coeff}_{D^j}\big\{(2D-n)\sum_{i=\frac{n}{2}+1}^{k}\binom{k}{i}(\frac{n}{2})^{k-i}
(D-n)^{i-\frac{n}{2}-1}\big\}=0.
 \end{equation}
If $k=n+1$ we have
 \begin{eqnarray}\label{eq-sum-ekj-k=n+1}
&&\sum_{j=0}^{\frac{n}{2}}(\frac{n}{2})^j e_{n+1,j}= \sum_{j=0}^{\frac{n}{2}}(\frac{n}{2})^j \mathrm{Coeff}_{D^j}\big\{(2D-n)\sum_{i=\frac{n}{2}+1}^{n+1}\binom{n+1}{i}(\frac{n}{2})^{n+1-i}
(D-n)^{i-\frac{n}{2}-1}\big\}\nn\\
&=& -(\frac{n}{2})^{\frac{n}{2}+1} \mathrm{Coeff}_{D^{\frac{n}{2}+1}}\big\{(2D-n)\sum_{i=\frac{n}{2}+1}^{n+1}\binom{n+1}{i}(\frac{n}{2})^{n+1-i}
(D-n)^{i-\frac{n}{2}-1}\big\}\nn\\
&=&-2(\frac{n}{2})^{\frac{n}{2}+1}.
 \end{eqnarray}
Now we add $(\frac{n}{2})^j\times$ $j$-th row of $N$ to the $0$-th row, for $j=1,2,\dots,\frac{n}{2}$. Then  (\ref{eq-sum-ekj-k<=n}) and (\ref{eq-sum-ekj-k=n+1}) mean that the resulted matrix has the $0$-th row of the form
\begin{equation*}
      (0,0,\dots,0,-2(\frac{n}{2})^{\frac{n}{2}+1}).
\end{equation*}
Combined with (\ref{eq-lower-left-ofN-unipotent}) this yields that the matrix $N$ is nonsingular. Hence the proof is complete.
\end{proof}
\begin{proposition}\label{prop-PicardFuchsOperator-broad}
\begin{equation*}
      v_k=\begin{cases}
      6n^{n+1} \hbar^{-n}\delta_{1,k},& \mbox{for odd}\ k;\\
      (-1)^{\frac{n-k}{2}+1}(\frac{n}{2})^{n-k+2}\binom{\frac{n}{2}+1}{\frac{k}{2}}
      -\delta_{2,k}4n^n \hbar^{-n}-\delta_{0,k}2n^{n+2} \hbar^{-n},& \mbox{for even}\ k.
      \end{cases}
\end{equation*}
\end{proposition}
\begin{proof}
By Lemma \ref{lem-PicardFuchs-EqCoeff-rank}, it suffices to show that  $v_i$'s in the assertion satisfy (\ref{eq-equationsOfCoeff-PicardFuchs-broad-0}) and (\ref{eq-equationsOfCoeff-PicardFuchs-broad-1}). The former is easy. We only show (\ref{eq-equationsOfCoeff-PicardFuchs-broad-1}). 
By (\ref{eq-Ek-modulo-narrowPicardFuchsOperator}) we have
\begin{eqnarray*}
 -2n^{n+2} \hbar^{-n}E_0+6n^{n+1}\hbar^{-n}E_1-4n^n \hbar^{-n}E_2
=-2n^{n}\hbar^{-n}D^{\frac{n}{2}+1}(2 D-n).
\end{eqnarray*}
A direct computation yields
\begin{eqnarray*}
&&\sum_{k=0}^{n+2} \big(\delta_{k,\mathrm{even}}(-1)^{\frac{n-k}{2}+1}(\frac{n}{2})^{n-k+2}\binom{\frac{n}{2}+1}{\frac{k}{2}}\sum_{i=0}^{\frac{n}{2}}\binom{k}{i}(\frac{n}{2})^{k-i}D^{\frac{n}{2}+i}\big)\\
&=& \sum_{i=0}^{\frac{n}{2}}\Big((\frac{n}{2})^{n+2-i}D^{\frac{n}{2}+i}\sum_{l=0}^{\frac{n}{2}+1} (-1)^{\frac{n}{2}-l+1}\binom{\frac{n}{2}+1}{l}\binom{2l}{i}\Big)=0,
\end{eqnarray*}
where we have used that for $k<\frac{n}{2}+1$
\[
\sum_{l=0}^{\frac{n}{2}+1} (-1)^{\frac{n}{2}-l+1}\binom{\frac{n}{2}+1}{l}l^k=0.
\]
Moreover we have
\begin{eqnarray*}
&& \sum_{k=0}^{n+2} \left(\delta_{k,\mathrm{even}}(-1)^{\frac{n-k}{2}+1}(\frac{n}{2})^{n-k+2}\binom{\frac{n}{2}+1}{\frac{k}{2}}
\sum_{i=\frac{n}{2}+1}^{k}\binom{k}{i}(\frac{n}{2})^{k-i}
(D-n)^{i-\frac{n}{2}-1}\right)\\
&=&\sum_{i=0}^{\frac{n}{2}+1}\sum_{l=0}^{\frac{n}{2}+1} \left((-1)^{\frac{n}{2}-l+1}(\frac{n}{2})^{\frac{n}{2}+1-i}\binom{\frac{n}{2}+1}{l}\binom{2l}{i+\frac{n}{2}+1}
(D-n)^{i}\right).
\end{eqnarray*}
We are left to show that the last sum  is equal to $D^{\frac{n}{2}+1}$. It suffices to show
\begin{equation*}
      \sum_{l=0}^{\frac{n}{2}+1} (-1)^{\frac{n}{2}-l+1}\binom{\frac{n}{2}+1}{l}\binom{2l}{i+\frac{n}{2}+1}
      =2^{\frac{n}{2}+1-i}\binom{\frac{n}{2}+1}{i},
\end{equation*}
which follows from
\begin{equation*}
      \sum_{l=0}^{n} (-1)^{n-l}\binom{n}{l}\binom{2l}{m}
      =2^{2n-m}\binom{n}{m-n}.
\end{equation*}
\end{proof}

\begin{corollary}\label{cor-PicardFuchsOperator-broad}
Let $q=1$. Then
\begin{equation*}
\big((\nabla_{\hbar\partial_{\hbar}}^2-\frac{n^2}{4})^{\frac{n}{2}+1}-4n^n \hbar^{-n}(\nabla_{\hbar\partial_{\hbar}}-\frac{n}{2})(\nabla_{\hbar\partial_{\hbar}}-n)\big)[z_1\Omega]=0.
\end{equation*}
\end{corollary}
\begin{proof}
This is a  consequence of (\ref{eq-APrioriPicardFuchs-broad}) and Proposition \ref{prop-PicardFuchsOperator-broad}:
\begin{eqnarray*}
\sum_{k=0}^{n+2}v_k(\hbar)D^k
&=&\sum_{j=0}^{\frac{n}{2}+1}(-1)^{\frac{n}{2}-j+1}(\frac{n}{2})^{n-2j+2}\binom{\frac{n}{2}+1}{j}D^{2j}
      -2n^n \hbar^{-n}(2D-n)(D-n)\nn\\
&=& (D^2-\frac{n^2}{4})^{\frac{n}{2}+1}-4n^n \hbar^{-n}(D-\frac{n}{2})(D-n).
\end{eqnarray*}
\end{proof}

\begin{lemma}\label{lem-vanishingCondition-ci-under-broadPicardFuchsOperator}
Let $D=\hbar\frac{\partial}{\partial \hbar}$, and let
\begin{eqnarray}\label{eq-PicardFuchsOperator-broad-q=1}
\mathcal{P}_{\mathrm{br}}:=(D^2-\frac{n^2}{4})^{\frac{n}{2}+1}-4n^n \hbar^{-n}(D-\frac{n}{2})(D-n).
\end{eqnarray}
Suppose $c_i\in \mathbb{C}(\hbar)$ for $1\leq i\leq \frac{n}{2}-1$, and
\begin{equation}\label{eq-vanishingCondition-ci-under-broadPicardFuchsOperator}
     \mathcal{P}_{\mathrm{br}}(\hbar^{\frac{n}{2}}\sum_{i=1}^{\frac{n}{2}-1}c_i \log^i\hbar)=0.
\end{equation}
Then $c_i=0$ for $1\leq i\leq \frac{n}{2}-1$.
\end{lemma}
\begin{proof}
Decompose $D^2-\frac{n^2}{4}$ as $(D+\frac{n}{2})(D-\frac{n}{2})$, and note that
\[
(D-\frac{n}{2})\circ \hbar^{\frac{n}{2}}= \hbar^{\frac{n}{2}}D.
\]
Since $D^k(\log^i \hbar)=0$ for $i<k$, we have
\[
(D^2-\frac{n^2}{4})^{\frac{n}{2}+1}(\hbar^{\frac{n}{2}}\log^i\hbar)=0
\]
for $i<\frac{n}{2}+1$. Then from
\begin{eqnarray*}
(D-n)(D-\frac{n}{2})\hbar^{\frac{n}{2}}\log^i\hbar
=\hbar^{\frac{n}{2}}\big( i(i-1)\log^{i-2}\hbar-\frac{ni}{2}\log^{i-1}\hbar\big),
\end{eqnarray*}
the assertion follows.
\end{proof}

\begin{corollary}\label{cor-periods-final}
Locally there exists a unique basis $\Xi_0,\dots,\Xi_{n+1}$ of the local system $R$,  such that
\begin{equation*}
      \int_{\Xi_i}e^{\frac{W_0}{\hbar}} \Omega=\xi_i,\  0\leq i\leq n,
\end{equation*}
\begin{equation*}
      \int_{\Xi_{n+1}}e^{\frac{W_0}{\hbar}} \Omega=0,
\end{equation*}
\begin{equation*}
      \int_{\Xi_i}e^{\frac{W_0}{\hbar}}z_1\Omega
      =\frac{q^{-1}\hbar^{\frac{n}{2}}}{2(-n)^{\frac{n}{2}}}(\hbar\frac{\partial}{\partial \hbar})^{\frac{n}{2}}\xi_i,\  0\leq i\leq n,
\end{equation*}
\begin{equation*}
      \int_{\Xi_{n+1}}e^{\frac{W_0}{\hbar}}z_1 \Omega=q^{-1}\hbar^{\frac{n}{2}}.
\end{equation*}
\end{corollary}
\begin{proof}
First let $q=1$.
For each $0\leq k\leq n$, by (\ref{eq-thetak-xik-relation-monodromy}) and (\ref{eq-equationsOfCoeff-PicardFuchs-broad-1}), $c_i(\Xi_k)$'s satisfy (\ref{eq-vanishingCondition-ci-under-broadPicardFuchsOperator}). Lemma \ref{lem-vanishingCondition-ci-under-broadPicardFuchsOperator} then yields that all  $c_i(\Xi_k)$'s vanish, and  the assertion follows from Proposition \ref{prop-fix-homologyBasis-preliminary}.

By the change of variables (\ref{eq-changeOfVar-lambda}) one has
\begin{equation}\label{eq-changeOfVar-broadPeriod}
      \int_{\Xi_i}e^{\frac{W_0}{\hbar}}z_1\Omega=q^{-\frac{1}{2}}\int_{\Xi_i}e^{\frac{q^{\frac{1}{n}}W_0(q=1)}{\hbar}}z_1\Omega.
\end{equation}
The assertion for general $q$ then follows from the $q=1$ case, (\ref{eq-changeOfVar-narrowPeriod}), and (\ref{eq-changeOfVar-broadPeriod}).
\end{proof}

\section{Mirror symmetry in terms of Frobenius manifolds}\label{sec:mirrorSymmetryFrobeniusManifolds}
With the results in the previous sections we are ready to complete the construction of a Frobenius manifold associated with our Landau-Ginzburg model $(Z^n,W_0,\Omega)$. Then by establishing a reconstruction theorem we can obtain the mirror symmetry for quadric hypersurfaces.
\subsection{The opposite space and the Frobenius manifold}\label{sec:opposite space}
Let $\Xi_0,\dots,\Xi_{n+1}$ be a local basis in Corollary \ref{cor-periods-final}. Following \cite{Gro11}, we denote the dual basis by $\alpha^i$ with $0\leq i\leq n$, and $\beta$. Then
\begin{equation}\label{eq-def-dualBasis-alphaBeta}
      [f\Omega]=\sum_{i=0}^{n}\alpha^i\int_{\Xi_i}fe^{\frac{W}{\hbar}}\Omega+\beta\int_{\Xi_{n+1}}fe^{\frac{W}{\hbar}}\Omega.
\end{equation}
By (\ref{eq-monodromy-Xik}) and Corollary \ref{cor-periods-final}, when $\hbar$ undergoes a counterclockwise loop in the $\mathbb{C}^{\times}_{\hbar}$ plane around $0$, the effect of monodromy on this dual basis is
\begin{equation*}
      \alpha^k\mapsto \exp(2n\pi\sqrt{-1}\alpha)\alpha^k,\ \beta\mapsto \beta.
\end{equation*}
Thus 
\begin{equation*}
      \{\hbar^{-n\alpha}\alpha^k|0\leq k\leq n\}\cup\{\beta\}
\end{equation*}
is a basis of  $\mathcal{H}$ as a free $\mathbb{C}\{\hbar,\hbar^{-1}\}$-module.
We take $\mathcal{H}_{-}$ to be the $\mathcal{O}(\mathbb{P}^1\backslash \{0\})$-submodule of $\mathcal{H}$ generated by
\begin{equation*}
      \{(\hbar \alpha)^k \hbar^{-1}\hbar^{-n\alpha}|0\leq k\leq n\}\cup\{ \hbar^{\frac{n}{2}-1}\beta\}.
\end{equation*}
Let $\Omega_0$ to be the section of $\mathcal{H}$ whose value at $0\in \widetilde{M}$ is $\Omega$.

\begin{proposition}\label{prop-filtration-FrobeniusManifold}
\begin{enumerate}
      \item[(i)] $\mathcal{H}=\mathcal{E}_q\oplus \mathcal{H}_{-}$ for $q\in \wcM_{\mathrm{red}}$.
      \item[(ii)] $\mathcal{H}_{-}$ is isotropic with respect to $\overline{\Omega}$.
      \item[(iii)] $\mathrm{Gr}$ preserves $\mathcal{H}_{-}$, and $\Omega_{0}$ represents an eigenvector of $\mathrm{Gr}_0$ with eigenvalue $-1$.
      \item[(iv)] $s'_0=\tau(\Omega_0\otimes 1)$ yields  miniversality along $\wcM_{\mathrm{red}}$.
\end{enumerate}
      Thus we get a Frobenius manifold structure on $\wcM$.
\end{proposition}
\begin{proof} 
(i) By Corollary \ref{cor-periods-final},
\begin{eqnarray}\label{eq-W0iOmega-expansion}
      [W_0^i\Omega]&=&\sum_{i=0}^{n}\alpha^i\int_{\Xi_i}W_0^ie^{\frac{W_0}{\hbar}}\Omega+\beta\int_{\Xi_{n+1}}W_0^ie^{\frac{W_0}{\hbar}}\Omega\nn\\
      &=& (-\hbar^2\frac{\partial}{\partial \hbar})^i\xi
      = (-\hbar^2\frac{\partial}{\partial \hbar})^i \big(q^{\alpha}\hbar^{-n \alpha}\sum_{d=0}^{+\infty}\hbar^{-nd}\frac{(2d)!}{(d!)^{n+2}}\frac{q^d\binom{2 \alpha+2d}{2d}}{\binom{\alpha+d}{d}^2}\big)\nn\\
      &\equiv& (n\hbar \alpha)^i \hbar^{-n \alpha} \mod \mathcal{H}_{-}.
\end{eqnarray}
\begin{eqnarray}\label{eq-z1Omega-expansion}
      [z_1\Omega]&=&\sum_{i=0}^{n}\alpha^i\int_{\Xi_i}z_1e^{\frac{W_0}{\hbar}}\Omega+\beta\int_{\Xi_{n+1}}z_1e^{\frac{W_0}{\hbar}}\Omega\nn\\
      &=& \frac{\hbar^{\frac{n}{2}}}{2(-n)^{\frac{n}{2}}q}(\hbar\frac{\partial}{\partial \hbar})^{\frac{n}{2}}\xi
      +\hbar^{\frac{n}{2}}\beta\nn\\
      &=& \frac{q^{-\frac{1}{2}}\hbar^{\frac{n}{2}}}{2(-n)^{\frac{n}{2}}}(\hbar\frac{\partial}{\partial \hbar})^{\frac{n}{2}}\big(q^{\alpha}\hbar^{-n \alpha}\sum_{d=0}^{+\infty}\hbar^{      -nd}\frac{(2d)!}{(d!)^{n+2}}\frac{q^d\binom{2 \alpha+2d}{2d}}{\binom{\alpha+d}{d}^2}\big)
         +\hbar^{\frac{n}{2}}\beta\nn\\
      &\equiv& \frac{1}{2q} (\hbar \alpha)^{\frac{n}{2}} \hbar^{-n \alpha}
      +q^{-1}\hbar^{\frac{n}{2}}\beta\mod \mathcal{H}_{-}.
\end{eqnarray}
So (i) follows from the definition of $\mathcal{H}_{-}$.

(ii) Since $\hbar^{-n\alpha}\alpha^k$ and $\beta$ are single-valued sections of $\mathcal{R}^{\vee}$, 
\begin{equation}\label{eq-pairing-singleValued}
\big((-)^*\hbar^{-n\alpha}\alpha^i,\hbar^{-n\alpha}\alpha^{j} \big),\ 
\big((-)^*\hbar^{-n\alpha}\alpha^i,\beta \big)
\end{equation}
are single valued. As argued in \cite[page 85]{Gro11}, it follows that (\ref{eq-pairing-singleValued}) are independent of $\hbar$. Thus
\begin{eqnarray*}
 0=\hbar\partial_{\hbar}\big((-)^*\hbar^{-n\alpha}\alpha^i,\beta \big)
=\big(-(-)^*n\hbar^{-n\alpha}\alpha^{i+1},\beta \big).
\end{eqnarray*}
Hence
\begin{equation*}
      \big((-)^*\hbar^{-n\alpha}\alpha^i,\beta \big)=0.
\end{equation*}
Moreover since $\beta$ is a flat section of $\mathcal{R}^{\vee}$, $\big((-1)^*\beta,\beta\big)$ is a constant, so by definition
\begin{equation*}
      (\beta,\beta)_{\mathcal{E}}=\mathrm{const}\cdot \hbar^{-n}.
\end{equation*}
Then for $i,j\geq 1$,
\[
\overline{\Omega}(\hbar^{\frac{n}{2}-i}\beta,\hbar^{\frac{n}{2}-j}\beta)
=\mathrm{Res}_{\hbar=0}(\hbar^{\frac{n}{2}-i}\beta,\hbar^{\frac{n}{2}-j}\beta)_{\mathcal{E}} \mathrm{d}\hbar
=\mathrm{Res}_{\hbar=0}(\mathrm{const}\cdot \hbar^{-i-j})\mathrm{d}\hbar=0.
\]
So $\mathcal{H}_{-}$ is isotropic. 

(iii) 
For $j\geq 1$,
\begin{eqnarray*}
&& \mathrm{Gr}\big((\hbar \alpha)^k \hbar^{-j}\hbar^{-n\alpha}\big)
=\nabla_{\hbar\partial_{\hbar}+E}^{\mathrm{GM}}
\big((\hbar \alpha)^k \hbar^{-j}\hbar^{-n\alpha}\big)
-(\hbar \alpha)^k \hbar^{-j}\hbar^{-n\alpha}\\
&=&(k-j-n \alpha-1)(\hbar \alpha)^k \hbar^{-j}\hbar^{-n\alpha}\in \mathcal{H}_{-},
\end{eqnarray*}
\begin{eqnarray*}
&& \mathrm{Gr}(\hbar^{\frac{n}{2}-j}\beta)
=\nabla_{\hbar\partial_{\hbar}+E}^{\mathrm{GM}}(\hbar^{\frac{n}{2}-j}\beta)
-\hbar^{\frac{n}{2}-j}\beta\\
&=&(\frac{n}{2}-j-1)\hbar^{\frac{n}{2}-j}\beta\in \mathcal{H}_{-}.
\end{eqnarray*}
So $\mathrm{Gr}$ preserves $\mathcal{H}_{-}$.

(iv) By Corollary \ref{cor-periods-final}, $\Omega_0\equiv \hbar^{-n\alpha}\mod \mathcal{H}_{-}$. Suppose $\Omega_0\equiv[f\Omega] \mod \mathcal{H}_{-}$, then $f(0)=1$ and we can write
\begin{equation}\label{eq-fOmega-expansion}
      [f\Omega]=q^{\alpha}\hbar^{-n\alpha}\sum_{i=0}^{n}\varphi_i(\mathbf{t},\hbar^{-1})(\hbar \alpha)^i
      +\varphi_{n+1}(\mathbf{t},\hbar^{-1}) \hbar^{\frac{n}{2}}\beta,
\end{equation}
where 
\begin{equation*}
      \varphi_i(\mathbf{t},\hbar^{-1})=\delta_{0,i}+\sum_{j=1}^{\infty}\varphi_{i,j}(\mathbf{t})\hbar^{-j}.
\end{equation*}
$f$ is uniquely determined by these conditions. It follows from Proposition \ref{prop-PicardFuchs-Quadric-solution} and Corollary \ref{cor-periods-final} that the dependence of $f$ on $t_1$ is throught $q=e^{t_1}$.
Recall the Gauss-Manin connection
\begin{equation}\label{eq-GM-connection-X}
      \nabla_X^{\mathrm{GM}}[f \Omega]=[\big(X(f)+\hbar^{-1}X(W)f\big)\Omega],
\end{equation}

By (\ref{eq-GM-connection-X}), (\ref{eq-EulerField-liftingProperty-final}), and (\ref{eq-EulerField-final}),
\begin{equation}\label{eq-derivativeInCoordinateT-Omega0}
      \hbar\nabla_{\partial_{t_i}}[f\Omega]=\begin{cases}
      [(\hbar\frac{\partial f}{\partial t_i}+W_0^i f)\Omega],& i\neq 1,n,n+1;\\
      [(\hbar\frac{\partial f}{\partial t_1}+\frac{1}{n}(W_0+\sum_{i=2}^nit_iW_0^i-2ne^{t_1}t_n+\frac{n}{2}t_{n+1}e^{t_1}z_1) f)\Omega],& i=1;\\
      [(\hbar\frac{\partial f}{\partial t_i}+(W_0^i-2e^{t_1}) f)\Omega],& i=n;\\
      [(\hbar\frac{\partial f}{\partial t_{n+1}}+e^{t_1}z_1 f)\Omega],& i=n+1,
      \end{cases}
\end{equation}
where we have used
\begin{eqnarray*}
\partial_{t_1}=\frac{1}{n}\big(E+\sum_{i=0}^n(i-1)t_i\partial_{t_i}+(\frac{n}{2}-1)t_{n+1}\partial_{t_{n+1}}\big)
\end{eqnarray*}
and
\begin{eqnarray*}
&&\tilde{E}(W)+\big(\sum_{i=0}^n(i-1)t_i\partial_{t_i}+(\frac{n}{2}-1)t_{n+1}\partial_{t_{n+1}}\big)W\\
&=& W+\big(\sum_{i=0}^n(i-1)t_i W_0^i-2(n-1)e^{t_1}t_n+(\frac{n}{2}-1)t_{n+1}e^{t_1}z_1\big)\\
&=& W_0+\sum_{i=0}^nit_i W_0^i-2ne^{t_1}t_n+\frac{n}{2}t_{n+1}e^{t_1}z_1.
\end{eqnarray*}
So
\begin{equation*}
      (\hbar\nabla_{\partial_{t_i}}[f\Omega])|_{\mathbf{t}=0}\stackrel{\mod \hbar \mathcal{E}}{\equiv}\begin{cases}
      [W_0^i \Omega],& i\neq 1,n+1;\\
      [\frac{W_0}{n}\Omega],& i=1;\\
      [z_1 \Omega],& i=n+1,
      \end{cases}
\end{equation*}
\begin{equation*}
      (\hbar\nabla_{\partial_{t_i}}[f\Omega])|_{t_i=0\ \mathrm{for}\ i\neq 1}\stackrel{\mod \hbar \mathcal{E}}{\equiv}\begin{cases}
      [W_0^i f_{t_i=0\ \mathrm{for}\ i\neq 1}\Omega],& i\neq 1,n+1;\\
      [\frac{W_0}{n}f_{t_i=0\ \mathrm{for}\ i\neq 1}\Omega],& i=1;\\
      [e^{t_1}z_1f_{t_i=0\ \mathrm{for}\ i\neq 1} \Omega],& i=n+1.
      \end{cases}
\end{equation*}
Since $\{W_0^i\}_{0\leq i\leq n}\cup\{z_1\}$ form a basis of the Milnor ring $\mathrm{MR}_n$ (recall (\ref{eq-quadric-mirror-jac-relations}) and (\ref{eq-W0=nz2})), we get the miniversality  along $\wcM$.
\end{proof}

\begin{lemma}\label{lem-mirrorTransform-t-To-y}
Let 
\begin{equation*}
      y_i=\varphi_{i,1}(\mathbf{t}),\ 0\leq i\leq n+1.
\end{equation*}
Then
\begin{equation*}
      y_i=\begin{cases}
       t_1+O(\mathbf{t}^2),& i=1;\\
       n^{\frac{n}{2}}t_{\frac{n}{2}}+\frac{1}{2q} t_{n+1}+O(\mathbf{t}^2),& i=\frac{n}{2};\\
      n^i t_i+ O(\mathbf{t}^2),& 0\leq i\leq n,\ i\neq 1,\frac{n}{2};\\
      t_{n+1}+O(\mathbf{t}^2),& i=n+1.
       \end{cases}
\end{equation*}
\end{lemma}
\begin{proof}
By (\ref{eq-fOmega-expansion}),
\begin{equation*}
      [f\Omega]\equiv \hbar^{-n\alpha}+\hbar^{-n\alpha}\alpha\log(q)
      +\sum_{i=0}^n y_i \hbar^{-n\alpha-1}(\hbar \alpha)^i
      +y_{n+1}\hbar^{\frac{n}{2}-1}\beta\mod \hbar^{-1} \mathcal{H}_{-}.
\end{equation*}
So
\begin{equation*}
      \hbar\nabla_{\partial_{y_i}}[f\Omega]\stackrel{\mod  \mathcal{H}_{-}}{\equiv}\begin{cases}
      \hbar^{-n\alpha}(\hbar \alpha)^i,& 0\leq i\leq n;\\
      \hbar^{\frac{n}{2}}\beta,& i=n+1.
      \end{cases}
\end{equation*}
By (\ref{eq-W0iOmega-expansion}) and (\ref{eq-z1Omega-expansion}), under the canonical isomorphisms
\begin{equation*}
      \mathcal{E}_0/\hbar\mathcal{E}_0\xleftarrow{\sim} \mathcal{E}_0\cap \hbar \mathcal{H}_{-}\xrightarrow{\sim} \hbar \mathcal{H}_{-}/ \mathcal{H}_{-}
\end{equation*}
one has
\begin{equation}\label{eq-W0i&z1-mapsto-periods}
[W_0^i\Omega]\mapsto (n\hbar \alpha)^i \hbar^{-n \alpha},\ 
[z_1\Omega]\mapsto  \frac{1}{2q} (\hbar \alpha)^{\frac{n}{2}} \hbar^{-n \alpha}+q^{-1}\hbar^{\frac{n}{2}}\beta. 
\end{equation}
\end{proof}

\begin{corollary}\label{cor-frobeniusAlg-MilnorRing}
The Frobenius algebra  along $t_0=t_i=0$ for $0\leq i\leq n$ is isomorphic to the Milnor ring of $W_0$.
\end{corollary}
\begin{proof}
Recall (\ref{eq-mirrorLGofQuadric-mirrorBasisToSmallQH})
\[
f_1=z_2=\frac{W_0}{n}=\frac{\partial}{\partial t_1}=\frac{\partial}{\partial y_1}+O(\mathbf{y}),
\]
and
\begin{eqnarray*}
&&f_2=(\sqrt{-1})^{\frac{n}{2}}\sqrt{2}(q z_1-\frac{1}{2}z_2^{\frac{n}{2}})
=(\sqrt{-1})^{\frac{n}{2}}\sqrt{2}(q\frac{\partial}{\partial t_{n+1}}-\frac{1}{2 n^{\frac{n}{2}}}\frac{\partial}{\partial t_{\frac{n}{2}}})\\
&=&(\sqrt{-1})^{\frac{n}{2}}\sqrt{2}(q\frac{\partial y_{n+1}}{\partial t_{n+1}}\frac{\partial}{\partial y_{n+1}}
+q\frac{\partial y_{\frac{n}{2}}}{\partial t_{n+1}}\frac{\partial}{\partial y_{\frac{n}{2}}}
-\frac{1}{2 n^{\frac{n}{2}}}\frac{\partial y_{\frac{n}{2}}}{\partial t_{\frac{n}{2}}}\frac{\partial}{\partial y_{\frac{n}{2}}})+O(\mathbf{y})\\
&=&(\sqrt{-1})^{\frac{n}{2}}\sqrt{2}q\frac{\partial}{\partial y_{n+1}}+O(\mathbf{y}).
\end{eqnarray*}
By using (\ref{eq-derivativeInCoordinateT-Omega0}), the assertion follows from Corollary \ref{cor-milnorRing-smallQuantumCohOfQn}.
\end{proof}

\begin{lemma}\label{lem-EulerField-LGMirrorOfQuadric-inFlatCoordinates}
In the coordinates $\{y_i\}_{0\leq i\leq n+1}$ the Euler field is equal to
\begin{eqnarray}\label{eq-EulerField-coordY}
      E=n\partial_{y_1}+\sum_{i=0}^n(1-i)y_i\partial_{y_i}
      +(1-\frac{n}{2})y_{n+1}\partial_{y_{n+1}}.
\end{eqnarray}
\end{lemma}
\begin{proof}
By (\ref{lem-mirrorTransform-t-To-y}),
\begin{align}\label{eq-EulerField-coordY-proof}
      E=&n\partial_{t_1}+\sum_{i=0}^{n}(1-i)t_i\partial_{t_i}+(1-\frac{n}{2})t_{n+1}\partial_{t_{n+1}}\nn\\
      = & n\partial_{y_1}+\sum_{\begin{subarray}{c}0\leq i\leq n\\ i\neq \frac{n}{2}\end{subarray}}(1-i)y_i\partial_{y_i}
      +(1-\frac{n}{2})(y_{\frac{n}{2}}-\frac{1}{2q}y_{n+1})\partial_{y_{\frac{n}{2}}}\nn\\
      &+(1-\frac{n}{2})y_{n+1}(\frac{1}{2q}\partial_{y_{\frac{n}{2}}}+\partial_{y_{n+1}})+O(\mathbf{y}^2)\nn\\
      =& n\partial_{y_1}+\sum_{\begin{subarray}{c}0\leq i\leq n\\ i\neq \frac{n}{2}\end{subarray}}(1-i)y_i\partial_{y_i}
      +(1-\frac{n}{2})y_{\frac{n}{2}}\partial_{y_{\frac{n}{2}}}
      +(1-\frac{n}{2})y_{n+1}\partial_{y_{n+1}}+O(\mathbf{y}^2)\nn\\
      =& n\partial_{y_1}+\sum_{i=0}^n(1-i)y_i\partial_{y_i}
      +(1-\frac{n}{2})y_{n+1}\partial_{y_{n+1}}+O(\mathbf{y}^2).
\end{align}
By Lemma \ref{cor-milnorRing-smallQuantumCohOfQn} and Corollary \ref{cor-frobeniusAlg-MilnorRing}, the B-model Frobenius manifold associated to $W$ is generically semisimple. 
Then by \cite[\S 2.4]{Man99}, the terms of order $\geq 2$ in (\ref{eq-EulerField-coordY-proof}) vanish.
\end{proof}

\subsection{Pairings}\label{sec:pairings}
\begin{proposition}\label{prop-pairing-flatSections}
\begin{equation*}
 \big((\hbar \alpha)^i \hbar^{-n\alpha},(\hbar \alpha)^j \hbar^{-n\alpha}\big)_{\hbar\mathcal{H}_{-}/\mathcal{H}_{-}}
=\begin{cases}
2,& \mbox{if}\ i+j=n;\\
8q,& \mbox{if}\ i=j=n;\\
0,& \mbox{otherwise},
\end{cases}
\end{equation*}
\begin{equation*}
      \big((\hbar \alpha)^i \hbar^{-n\alpha},\hbar^{\frac{n}{2}}\beta\big)_{\hbar\mathcal{H}_{-}/\mathcal{H}_{-}}=0,\ \mbox{for}\ 0\leq i\leq n,
\end{equation*}
\begin{equation*}
      (\hbar^{\frac{n}{2}}\beta,\hbar^{\frac{n}{2}}\beta)_{\hbar\mathcal{H}_{-}/\mathcal{H}_{-}}
      =\frac{(-1)^{\frac{n}{2}}}{2}.
\end{equation*}
\end{proposition}
\begin{proof}
By (\ref{eq-W0iOmega-expansion}) and (\ref{eq-z1Omega-expansion}),
\begin{equation}\label{eq-prop-pairing-flatSections-proof-1}
\begin{split}
     & [W_0^i \Omega]=(n\hbar \alpha)^i \hbar^{-n\alpha}+\mbox{terms in}\ \mathcal{H}_{-},\\
& [z_1\Omega]= \frac{(\hbar\alpha)^{\frac{n}{2}}}{2q}\hbar^{-n\alpha}+\hbar^{\frac{n}{2}}\beta
+\mbox{terms in}\ \mathcal{H}_{-}.
\end{split}
\end{equation}
By the stationary phase approximation (e.g. \cite[Prop. 2.35]{Gro11}) and Lemma \ref{lem-hessianValue-quadric}, we have
\begin{eqnarray*}
&& ([W_0^i\Omega],[W_0^j\Omega])_{\mathcal{E}}
=\frac{(-1)^{n(n+1)/2}}{(2\pi\sqrt{-1}\hbar)^n}\sum_{p\in \mathrm{Crit}(W_0)}(\int_{\Delta_p^{-}}e^{-\frac{W_0}{\hbar}}W_0^i\Omega)(\int_{\Delta_p^{+}}e^{\frac{W_0}{\hbar}}W_0^i\Omega)\\
&=&\sum_{p\in \mathrm{Crit}(W_0)}\frac{W_0^i(p,0)W_0^j(p,0)}{\mathrm{Hess}(W_0(p))}+O(\hbar)\\
&=&  \big(\sum_{k=1}^{n}\frac{(n\zeta^k4^{\frac{1}{n}}q^{\frac{1}{n}})^i (n\zeta^k4^{\frac{1}{n}}q^{\frac{1}{n}})^j}{2n q}
-2\times\delta_{i,0}\delta_{j,0}\frac{1}{4q}
\big)+O(\hbar)\\
&=& \begin{cases}
2n^n+O(\hbar),& \mbox{if}\ i+j=n;\\
8n^{2n}q+O(\hbar),& \mbox{if}\ i=j=n;\\
0,& \mbox{otherwise},
\end{cases}
\end{eqnarray*}
and
\begin{eqnarray*}
&& ([z_1\Omega],[z_1\Omega])_{\mathcal{E}}
=\frac{(-1)^{n(n+1)/2}}{(2\pi\sqrt{-1}\hbar)^n}\sum_{p\in \mathrm{Crit}(W_0)}(\int_{\Delta_p^{-}}e^{-\frac{W_0}{\hbar}}W_0^i\Omega)(\int_{\Delta_p^{+}}e^{\frac{W_0}{\hbar}}z_1\Omega)\\
&=& \sum_{p\in \mathrm{Crit}(W_0)}\frac{z_1(p,0)^2}{\mathrm{Hess}(W_0(p))}+O(\hbar)\\
&=& \big(\sum_{k=1}^{n}\frac{(\zeta^{-\frac{kn}{2}}q^{-\frac{1}{2}})^2 }{2n q}
-\frac{\big((\sqrt{-1})^{\frac{n-2}{2}}q^{-\frac{1}{2}}\big)^2}{4q}-\frac{\big(-(\sqrt{-1})^{\frac{n-2}{2}}q^{-\frac{1}{2}}\big)^2}{4q}
\big)+O(\hbar)\\
&=& \frac{1}{2q^2}+\frac{(-1)^{\frac{n}{2}}}{2q^2}+O(\hbar),
\end{eqnarray*}
and
\begin{eqnarray*}
&& ([W_0^i\Omega],[z_1\Omega])_{\mathcal{E}}
=\frac{(-1)^{n(n+1)/2}}{(2\pi\sqrt{-1}\hbar)^n}\sum_{p\in \mathrm{Crit}(W_0)}(\int_{\Delta_p^{-}}e^{-\frac{W_0}{\hbar}}z_1\Omega)(\int_{\Delta_p^{+}}e^{\frac{W_0}{\hbar}}z_1\Omega)\\
&=& \sum_{p\in \mathrm{Crit}(W_0)}\frac{W_0(p,0)^i z_1(p,0)}{\mathrm{Hess}(W_0(p))}+O(\hbar)\\
&=& \big(\sum_{k=1}^{n}\frac{(n\zeta^k4^{\frac{1}{n}}q^{\frac{1}{n}})^i(\zeta^{-\frac{kn}{2}}q^{-\frac{1}{2}}) }{2n q}
-\frac{\delta_{i,0}\times(\sqrt{-1})^{\frac{n-2}{2}}q^{-\frac{1}{2}}}{4q}-\frac{\delta_{i,0}\times\big(-(\sqrt{-1})^{\frac{n-2}{2}}q^{-\frac{1}{2}}\big)}{4q}
\big)+O(\hbar)\\
&=& \delta_{i,\frac{n}{2}} \frac{n^{\frac{n}{2}}}{q}+O(\hbar).
\end{eqnarray*}
In view of (\ref{eq-prop-pairing-flatSections-proof-1}) and Construction \ref{cons-Barannikov} (vi), these pairings  yield
\begin{equation*}
 \big((\hbar \alpha)^i \hbar^{-n\alpha},(\hbar \alpha)^j \hbar^{-n\alpha}\big)_{\hbar\mathcal{H}_{-}/\mathcal{H}_{-}}
=\begin{cases}
2,& \mbox{if}\ i+j=n;\\
8q,& \mbox{if}\ i=j=n;\\
0,& \mbox{otherwise},
\end{cases}
\end{equation*}
\begin{eqnarray*}
  \big(\frac{1}{2q} (\hbar \alpha)^{\frac{n}{2}} \hbar^{-n \alpha}+q^{-1}\hbar^{\frac{n}{2}}\beta, 
  \frac{1}{2q} (\hbar \alpha)^{\frac{n}{2}} \hbar^{-n \alpha}+q^{-1}\hbar^{\frac{n}{2}}\beta\big)_{\hbar\mathcal{H}_{-}/\mathcal{H}_{-}}
= \frac{1}{2q^2}+\frac{(-1)^{\frac{n}{2}}}{2q^2},
\end{eqnarray*}
\begin{eqnarray*}
  \big((\hbar \alpha)^i \hbar^{-n \alpha}, \frac{1}{2q} (\hbar \alpha)^{\frac{n}{2}} \hbar^{-n \alpha}+q^{-1}\hbar^{\frac{n}{2}}\beta\big)_{\hbar\mathcal{H}_{-}/\mathcal{H}_{-}}
= \delta_{i,\frac{n}{2}} \frac{n^{\frac{n}{2}-i}}{q},
\end{eqnarray*}
and our assertion follows.
\end{proof}

\begin{remark}\label{rem:orientation-stationaryPhaseApproximation}
The factor $(-1)^{\frac{n(n+1)}{2}}$ is eliminated in the above computations due to the choice of the orientations of $\Delta_p^+$ and $\Delta_p^{-}$. Indeed, if $w_i=u_i+\sqrt{-1}v_i$ are  analytic local coordinates such that $W_0=\sum_i w_i^2$ at $p$, then $\Delta_p^{+}$ is the submanifold defined by $u_i=0$ for $1\leq i\leq n$, and  $\Delta_p^{-}$ is the submanifold defined by $v_i=0$ for $1\leq i\leq n$. The orientations of  $\Delta_p^+$ and $\Delta_p^{-}$ are chosen such that
\begin{equation}\label{eq-orientation-DeltaP-1}
      \Delta_p^+\cdot\Delta_p^{-}=1.
\end{equation}
On the other hand, in the application of the stationary phase approximation, the choice of the square root $A=\sqrt{\frac{\partial^2 W_0}{\partial z_i\partial z_j}}$ is made such that, with $\mathbf{w}=A \mathbf{z}$, the orientation of $\mathbb{R}^n$ in the integration
\begin{equation*}
      \int_{\mathbb{R}^n}e^{-v_1^2-\dots-v_n^2}\mathrm{d}v_1\cdots \mathrm{d}v_n
\end{equation*}
coincides with the orientation of $\Delta_p^+$, and similarly for $\Delta_p^-$. Note that 
\begin{equation}\label{eq-orientation-DeltaP-2}
\mathrm{d}v_1\cdots  \mathrm{d}v_n \mathrm{d}u_1\cdots \mathrm{d}u_n=(-1)^{\frac{n(n+1)}{2}}\mathrm{d}u_1 \mathrm{d}u_2\cdots \mathrm{d}u_n \mathrm{d}v_n,
\end{equation}
and $\mathrm{d}u_1 \mathrm{d}u_2\cdots \mathrm{d}u_n \mathrm{d}v_n$ is our standard orientation of $\mathbb{C}^n$. Then the elimination of the factor $(-1)^{\frac{n(n+1)}{2}}$ follows from (\ref{eq-orientation-DeltaP-1}) and (\ref{eq-orientation-DeltaP-2}).\pqed
\end{remark}

\begin{corollary}\label{cor-isoOfFrobAlg-milnorRing-smallQuantumCohOfQn}
The map in Corollary \ref{cor-milnorRing-smallQuantumCohOfQn} induces an isomorphism  from the Frobenius algebra at the origin of the Frobenius manifold obtained in Proposition \ref{prop-filtration-FrobeniusManifold}  to the Frobenius algebra of the Frobenius manifolds associated with the small quantum cohomology of $Q^n$.
\end{corollary}
\begin{proof}
Under (\ref{eq-W0i&z1-mapsto-periods}),
\begin{equation*}
\begin{split}
      f_1^i=z_2^i\mapsto (\hbar \alpha)^i \hbar^{-n \alpha},
      f_2=(\sqrt{-1})^{\frac{n}{2}}(2q z_1-z_2^{\frac{n}{2}})\mapsto (\sqrt{-1})^{\frac{n}{2}}2\hbar^{\frac{n}{2}}\beta
\end{split}
\end{equation*}
Then by Proposition \ref{prop-pairing-flatSections},
\begin{equation*}
      \begin{cases}
      (f_1^i,f_1^j)=2\delta_{i+j,n}+8q\delta_{i,n}\delta_{j,n},\\
      (f_1^i,f_2)=0,\\
      (f_2,f_2)=2.
      \end{cases}
\end{equation*}
Since (the small quantum power) $\sfh^{\bullet n}=\sfh^n+2q$, the pairings coincides with those in the quantum cohomology of $Q^n$.

\end{proof}

\subsection{Reconstruction of the Frobenius manifold associated with \texorpdfstring{$(Z^n,W_0)$}{(Zn,W0)}}
Denote the (germ of) Frobenius manifold obtained in Proposition \ref{prop-filtration-FrobeniusManifold} by $\mathrm{FM}_B(Z^n,W_0)$. 
In this section we show the reconstruction of the associated generating function and hence the mirror symmetry of quadric hypersurfaces. Let $y_i$, $0\leq i\leq n+1$ be the dual basis with respect to the basis
\begin{equation*}
      \{(\hbar \alpha)^i \hbar^{-n\alpha}\}_{0\leq i\leq n-1}\cup\{(\hbar \alpha)^n\hbar^{-n\alpha}-2q\hbar^{-n\alpha}\}\cup\{(\sqrt{-1})^{\frac{n}{2}}2\hbar^{\frac{n}{2}}\beta\}.
\end{equation*}
By Corollary \ref{cor-isoOfFrobAlg-milnorRing-smallQuantumCohOfQn} we have
\begin{equation}\label{eq-pairing-FrobManifold-LGMirro}
      g_{i,j}:=g(\partial_{y_i},\partial_{y_j})=\begin{cases}
      2 \delta_{n,i+j},& \mbox{if}\ 0\leq i,j\leq n;\\
      2,& \mbox{if}\ i=j=n+1;\\
      0,& \mbox{otherwise}.
      \end{cases}
\end{equation}
Let $(g^{i,j})_{0\leq i,j\leq n+1}$ be the inverse matrix to $(g_{i,j})_{0\leq i,j\leq n+1}$.
We have  be a formal power series $F(y_0,\dots,y_{n+1})\in \mathbb{C}[[y_0,\dots,y_{n+1}]]$ satisfying the WDVV equation
\begin{equation}\label{eq-WDVV-inTermsOfPartialDerivatives}
      \sum_{i=0}^{n+1}\sum_{j=0}^{n+1}(\partial_{y_a}\partial_{y_b}\partial_{y_i}F)g^{i,j}(\partial_{y_j}\partial_{y_c}\partial_{y_d}F)
      =\sum_{i=0}^{n+1}\sum_{j=0}^{n+1}(\partial_{y_a}\partial_{y_c}\partial_{y_i}F)g^{i,j}(\partial_{y_j}\partial_{y_b}\partial_{y_d}F),
\end{equation}
and
\begin{equation}\label{eq-EulerField}
  E(F)=(3-n)F+\mbox{quadratic terms},
\end{equation}
where $E$ is the Euler field (\ref{eq-EulerField-coordY}).
Denote
\begin{equation}\label{eq-correlators-BSide}
  \langle \gamma_{i_1},\dots,\gamma_{i_l}\rangle:=(\partial_{y_{i_1}}\dots\partial_{y_{i_l}}F)(0).
\end{equation}
We call $\langle \gamma_{i_1},\dots,\gamma_{i_l}\rangle$ a \emph{correlator of length $l$}.
For a multiset $S=\{i_1,\dots,i_l\}$ with elements in $[0,n+1]\cap \mathbb{Z}$, we denote by $\gamma_S$ the sequence $\gamma_{i_1},\dots,\gamma_{i_l}$. Then the WDVV equation can be written as, for any $S$,
\begin{eqnarray}\label{eq-WDVV}
      &&\sum_{i=0}^{n+1}\sum_{j=0}^{n+1}\sum_{S_1\sqcup S_2=S}\langle \gamma_a,\gamma_b,\gamma_{S_1},\gamma_i\rangle g^{i,j}\langle \gamma_j,\gamma_{S_2},\gamma_c,\gamma_d\rangle\nn\\
      &=&\sum_{i=0}^{n+1}\sum_{j=0}^{n+1}\sum_{S_1\sqcup S_2=S}\langle \gamma_a,\gamma_c,\gamma_{S_1},\gamma_i\rangle g^{i,j}\langle \gamma_j,\gamma_{S_2},\gamma_b,\gamma_d\rangle.
\end{eqnarray}
By construction, $F$ has no terms of degree $\leq 2$.

\begin{lemma}\label{lem-nonzerothreepointinvariants-LGMirrorOfQuadric}
Using the notation (\ref{eq-correlators-BSide}), the degree 3 terms of $F$ are given by
\begin{equation}\label{eq-nonzerothreepointinvariants-LGMirrorOfQuadric}
      \begin{cases}
      \langle \gamma_{n+1},\gamma_{n+1}, \gamma_0\rangle=2,\ \langle \gamma_{n+1},\gamma_{n+1}, \gamma_{n}\rangle=-4,\\
      \langle \gamma_{i},\gamma_{j},\gamma_{k}\rangle=2,\ \mbox{if}\ i+j+k=n,\\
      \langle \gamma_{i},\gamma_{n-i},\gamma_{n}\rangle=4,\ \mbox{if}\ 1\leq i\leq n-1,\\
      \langle \gamma_{i},\gamma_{j},\gamma_{k}\rangle=8,\ \mbox{if}\ i+j+k=2n,\ \mbox{and}\ i,j,k\neq n,\\
      \langle \gamma_{n},\gamma_{n},\gamma_{n}\rangle=8.
      \end{cases}
\end{equation}
\end{lemma}
\begin{proof}
Follows from Lemma \ref{lem-nonzerothreepointinvariants-quadric} and Corollary \ref{cor-isoOfFrobAlg-milnorRing-smallQuantumCohOfQn}.
\end{proof}

We say that $\gamma_i$ for $0\leq i\leq n$ are \emph{narrow} insertions, and $\gamma_{n+1}$ is a \emph{broad} insertion. Note that $\gamma_0=\mathbbm1$ is the flat identity.

\begin{lemma}\label{lem-quadric-vanishing-genus0GWInv-oddPrimInsertions}
Let $k$ be a non-negative integer. Let $S$ be a multiset with elements from $[0,n+1]\cap \mathbb{Z}$, and suppose that $S$ has exactly $2k+3$ copies of  $n+1$. Then $\langle \gamma_S\rangle=0$.
\end{lemma}
\begin{proof}
The assertion when $k=0$ follows from Lemma \ref{lem-nonzerothreepointinvariants-LGMirrorOfQuadric}. We show the general case by induction on $k$.
First we show $\langle \gamma_{n+1},\dots,\gamma_{n+1}\rangle_{0,2k+3}=0$. 
 In (\ref{eq-WDVV}), we take $S$ to be the multiset of $2k+1$ copies of $\gamma_{n+1}$, and $\gamma_a=\gamma_b=\gamma_{n+1}$, $\gamma_c=\gamma_1$, $\gamma_d=\gamma_n$. Then the induction assumption and (\ref{eq-pairing-FrobManifold-LGMirro}) yield
\begin{eqnarray*}
      &&(2k+1)\langle \underbrace{\gamma_{n+1},\dots,\gamma_{n+1}}_{2k+3}\rangle\times\frac{1}{2}\times\langle \gamma_{n+1},\gamma_{n+1},\gamma_1,\gamma_n\rangle
      +\langle \underbrace{\gamma_{n+1},\dots,\gamma_{n+1}}_{2k+3},\gamma_1\rangle\times\frac{1}{2}\times\langle \gamma_{n-1},\gamma_1,\gamma_n\rangle\nn\\
      &=&\langle \gamma_1,\underbrace{\gamma_{n+1},\dots,\gamma_{n+1}}_{2k+3}\rangle\times\frac{1}{2}\times\langle \gamma_{n+1},\gamma_{n+1},\gamma_n\rangle
      +\langle \gamma_1,\gamma_{n+1},\gamma_{n+1}\rangle\times\frac{1}{2}\times\langle \underbrace{\gamma_{n+1},\dots,\gamma_{n+1}}_{2k+3},\gamma_n\rangle.
\end{eqnarray*}
Thus from (\ref{eq-nonzerothreepointinvariants-quadric}) we get
\begin{eqnarray*}
      -2(2k+1)\langle\underbrace{\gamma_{n+1},\dots,\gamma_{n+1}}_{2k+3}\rangle
      +2\langle \underbrace{\gamma_{n+1},\dots,\gamma_{n+1}}_{2k+3},\gamma_1\rangle
      =-2\langle \gamma_1,\underbrace{\gamma_{n+1},\dots,\gamma_{n+1}}_{2k+3}\rangle.
\end{eqnarray*}
By the Euler field (Lemma \ref{lem-EulerField-LGMirrorOfQuadric-inFlatCoordinates}) we have
\begin{equation*}
      n\langle \underbrace{\gamma_{n+1},\dots,\gamma_{n+1}}_{2k+3},\gamma_1\rangle+(2k+3)(1-\frac{n}{2})\langle \underbrace{\gamma_{n+1},\dots,\gamma_{n+1}}_{2k+3}\rangle
      =(3-n)\langle \underbrace{\gamma_{n+1},\dots,\gamma_{n+1}}_{2k+3}\rangle,
\end{equation*}
i.e.
\begin{equation*}
      \langle \underbrace{\gamma_{n+1},\dots,\gamma_{n+1}}_{2k+3},\gamma_1\rangle
      =(k+\frac{1}{2}-\frac{2k}{n})\langle \underbrace{\gamma_{n+1},\dots,\gamma_{n+1}}_{2k+3}\rangle.
\end{equation*}
Since
\begin{equation*}
      4(k+\frac{1}{2}-\frac{2k}{n})-(4k+2)=-\frac{8k}{n}<0,
\end{equation*}
we obtain
\begin{equation*}
      \langle \underbrace{\gamma_{n+1},\dots,\gamma_{n+1}}_{2k+3}\rangle=0.
\end{equation*}

Now we show the general case inductively on the number of narrow insertions. Suppose the assertion holds for correlators with exactly $2k+3$ broad insertions and less than or equal to $m$ narrow insertions. In (\ref{eq-WDVV}), we take $S$ to be the union of the multiset of $2k+1$ copies of $\gamma_{n+1}$, and a multiset $S'$ of $m$ narrow insertions. Let $\gamma_a=\gamma_b=\gamma_{n+1}$, $\gamma_c=\gamma_1$, $\gamma_d=\gamma_r$ for some $0\leq r\leq n$. Then the inductive assumption and Euler field yields
\begin{eqnarray*}
      \langle \underbrace{\gamma_{n+1},\dots,\gamma_{n+1}}_{2k+3},\gamma_{S'},\gamma_i\rangle g^{i,j}\langle \gamma_j,\gamma_1,\gamma_r\rangle=0.
\end{eqnarray*}
If $0\leq r\leq n-2$, $\langle \gamma_j,\gamma_1,\gamma_r\rangle_{0,3}=2\delta_{i,n-1-r}$. Thus $\langle \underbrace{\gamma_{n+1},\dots,\gamma_{n+1}}_{2k+3},\gamma_{S'},\gamma_{r+1}\rangle=0$. If $r=n-1$, $\langle \gamma_j,\gamma_1,\gamma_r\rangle_{0,3}=2\delta_{i,0}+4\delta_{n,0}$. Since $\mathbbm1$ is the flat identity, we have (\cite[Corollary 2.1.1]{Man99}) $\langle \underbrace{\gamma_{n+1},\dots,\gamma_{n+1}}_{2k+3},\gamma_{S'},\mathbbm1\rangle=0$, and thus $\langle \underbrace{\gamma_{n+1},\dots,\gamma_{n+1}}_{2k+3},\gamma_{S'},\gamma_n\rangle=0$.

\end{proof}

\begin{proposition}\label{prop-reconstruction-FrobManifold}
The formal series $F$ is uniquely determined by the WDVV equation (\ref{eq-WDVV-inTermsOfPartialDerivatives}), the Euler field (\ref{eq-EulerField}), and the length 3 correlators (\ref{eq-nonzerothreepointinvariants-LGMirrorOfQuadric}).
\end{proposition}
\begin{proof}
We show the assertion by induction on the length $k$ of the correlators. We use $\sim$ to denote that the difference of the expressions on both sides is a combination of correlators that have been reconstructed in the induction.

 Firstly  suppose that we have reconstructed the correlators with less than $2k+2$ broad insertions. We are going to show $\langle \gamma_{n+1},\dots,\gamma_{n+1}\rangle_{0,2k+2}\sim 0$ for $k\geq 1$. 
 In (\ref{eq-WDVV}), we take $S$ to be the multiset of $2k$ copies of $\gamma_{n+1}$, and $\gamma_a=\gamma_b=\gamma_{n+1}$, $\gamma_c=\gamma_1$, $\gamma_d=\gamma_n$. Then the assumption yields
\begin{eqnarray*}
      &&2k\langle \underbrace{\gamma_{n+1},\dots,\gamma_{n+1}}_{2k+2}\rangle\times\frac{1}{2}\times\langle \gamma_{n+1},\gamma_{n+1},\gamma_1,\gamma_n\rangle
      +\langle \underbrace{\gamma_{n+1},\dots,\gamma_{n+1}}_{2k+2},\gamma_1\rangle\times\frac{1}{2}\times\langle \gamma_{n-1},\gamma_1,\gamma_n\rangle\nn\\
      &\sim &\langle \gamma_1,\underbrace{\gamma_{n+1},\dots,\gamma_{n+1}}_{2k+2}\rangle\times\frac{1}{2}\times\langle \gamma_{n+1},\gamma_{n+1},\gamma_n\rangle
      +\langle \gamma_1,\gamma_{n+1},\gamma_{n+1}\rangle\times\frac{1}{2}\times\langle \underbrace{\gamma_{n+1},\dots,\gamma_{n+1}}_{2k+2},\gamma_n\rangle.
\end{eqnarray*}
Thus from (\ref{eq-nonzerothreepointinvariants-quadric}) we get
\begin{eqnarray*}
      -4k\langle\underbrace{\gamma_{n+1},\dots,\gamma_{n+1}}_{2k+2}\rangle
      +2\langle \underbrace{\gamma_{n+1},\dots,\gamma_{n+1}}_{2k+2},\gamma_1\rangle
      \sim-2\langle \gamma_1,\underbrace{\gamma_{n+1},\dots,\gamma_{n+1}}_{2k+2}\rangle.
\end{eqnarray*}
By the Euler field (Lemma \ref{lem-EulerField-LGMirrorOfQuadric-inFlatCoordinates}) we have
\begin{equation*}
      n\langle \underbrace{\gamma_{n+1},\dots,\gamma_{n+1}}_{2k+2},\gamma_1\rangle+(2k+2)(1-\frac{n}{2})\langle \underbrace{\gamma_{n+1},\dots,\gamma_{n+1}}_{2k+2}\rangle
      =(3-n)\langle \underbrace{\gamma_{n+1},\dots,\gamma_{n+1}}_{2k+2}\rangle,
\end{equation*}
i.e.
\begin{equation*}
      \langle \underbrace{\gamma_{n+1},\dots,\gamma_{n+1}}_{2k+2},\gamma_1\rangle
      =(k-\frac{2k-1}{n})\langle \underbrace{\gamma_{n+1},\dots,\gamma_{n+1}}_{2k+2}\rangle.
\end{equation*}
Since
\begin{equation*}
      4(k-\frac{2k-1}{n})-4k=-\frac{8k-4}{n}<0,
\end{equation*}
we obtain
\begin{equation*}
      \langle \underbrace{\gamma_{n+1},\dots,\gamma_{n+1}}_{2k+2}\rangle\sim0.
\end{equation*}

Secondly we reconstruct the general case inductively on the number of narrow insertions. Fix $k\geq 0$. Suppose that we have computed correlators with less than $2k+2$ broad insertions, and those with exactly $2k+2$ broad insertions and less than or equal to $m$ narrow insertions. In (\ref{eq-WDVV}), we take $S$ to be the union of the multiset of $2k$ copies of $\gamma_{n+1}$, and a multiset $S'$ of $m$ narrow insertions, and let $\gamma_a=\gamma_b=\gamma_{n+1}$, $\gamma_c=\gamma_1$, $\gamma_d=\gamma_r$ for some $0\leq r\leq n$. Then 
\begin{eqnarray*}
      \langle \underbrace{\gamma_{n+1},\dots,\gamma_{n+1}}_{2k+2},\gamma_{S'},\gamma_i\rangle g^{i,j}\langle \gamma_j,\gamma_1,\gamma_r\rangle\sim0.
\end{eqnarray*}
If $0\leq r\leq n-2$, $\langle \gamma_j,\gamma_1,\gamma_r\rangle_{0,3}=2\delta_{i,n-1-r}$. Thus $\langle \underbrace{\gamma_{n+1},\dots,\gamma_{n+1}}_{2k+2},\gamma_{S'},\gamma_{r+1}\rangle=0$. If $r=n-1$, $\langle \gamma_j,\gamma_1,\gamma_r\rangle_{0,3}=2\delta_{i,0}+4\delta_{n,0}$. Since $\mathbbm1$ is the flat identity, we have (\cite[Corollary 2.1.1]{Man99}) $\langle \underbrace{\gamma_{n+1},\dots,\gamma_{n+1}}_{2k+2},\gamma_{S'},\mathbbm1\rangle=0$, and thus $\langle \underbrace{\gamma_{n+1},\dots,\gamma_{n+1}}_{2k+2},\gamma_{S'},\gamma_n\rangle=0$. 

Finally we reconstruct the correlators with only narrow insertions. By (\ref{eq-pairing-FrobManifold-LGMirro}) and Lemma \ref{lem-quadric-vanishing-genus0GWInv-oddPrimInsertions}, the correlators with only narrow insertions give a Frobenius (sub)manifold whose Frobenius algebra at the origin is generated by $\gamma_1$, and this reconstruction essentially is a consequence of \cite[Theorem 3.1]{KM94}. Here we sketch a proof for the reader's convenience. Suppose we have reconstructed the correlators of length $\leq k$ with only narrow insertions. Given a correlator of length $k+1$ with only narrow insertions. If it has an insertion $\gamma_1$, then we can eliminate it by the Euler field. In (\ref{eq-WDVV}), we take $\gamma_a=\gamma_b=\gamma_1$, and suitable $\gamma_c=\gamma_i$, $\gamma_d=\gamma_j$ and $S$, where $i,j\geq 2$. Then we can eliminate $\gamma_2$. By induction on $r$, suppose that we can eliminate $\gamma_r$.  In (\ref{eq-WDVV}), we take $\gamma_a=\gamma_1,\gamma_b=\gamma_r$, and suitable $\gamma_c=\gamma_i$, $\gamma_d=\gamma_j$ and $S$, where $i,j\geq r+1$. Then we can eliminate $\gamma_{r+1}$. 

\end{proof}

\begin{remark}\label{rem:oddBroadInsertions-vanishing}
The proof of Lemma \ref{lem-quadric-vanishing-genus0GWInv-oddPrimInsertions} can be united into the reconstruction in Proposition \ref{prop-reconstruction-FrobManifold}. We separate it in order to emphasize that this vanishing result follows from merely WDVV and the length 3 correlator. In \cite{HK21}, we deduced this vanishing in quantum cohomology of quadrics from the deformation invariance of Gromov-Witten invariants.
\end{remark}

\begin{theorem}\label{thm-mirrorSym-FrobManifolds-quadric}
Let $F_{Q^n}(y_0,\dots,y_{n+1})\in \mathbb{C}[[y_0,\dots,y_n]]$ be the generating function of genus 0 Gromov-Witten invariants of $Q^n$ in the basis $\mathbbm1=\sfh_0,\sfh_1,\dots,\sfh_n,\sfp$. Then $F_{Q^n}(y_0,\dots,y_{n+1})=F(y_0,\dots,y_{n+1})$. 
Consequently, there is an isomorphism of germs of Frobenius manifolds $\mathrm{FM}_A(Q^n)\cong \mathrm{FM}_B(Z^n,W_0)$.
\end{theorem}
\begin{proof}
By Lemma \ref{lem-nonzerothreepointinvariants-quadric}, the data
 (\ref{eq-pairing-FrobManifold-LGMirro}) and (\ref{eq-nonzerothreepointinvariants-LGMirrorOfQuadric}) match the corresponding data for $F_{Q^n}(y_0,\dots,y_{n+1})$ in the given basis. The Euler fields are also the same. Hence Proposition \ref{prop-reconstruction-FrobManifold} yields our assertion.
\end{proof}

\begin{corollary}
$ \mathrm{FM}_B(Z^n,W_0)$ is the germ of an analytic Frobenius manifold.
\end{corollary}
\begin{proof}
By \cite[Theorem 5.9]{HK21}, $F_{Q^n}(y_0,\dots,y_{n+1})$ has a positive convergence radius. The assertion then follows from Theorem \ref{thm-mirrorSym-FrobManifolds-quadric}.
\end{proof}



\end{document}